\documentclass[12pt]{amsart}
 \pdfoutput=1

\usepackage{rotating}
\usepackage{multirow}
\usepackage{amsmath}

\usepackage{amsfonts,amssymb,amscd,bbm,booktabs,color,enumerate,float,graphicx,latexsym, multirow,mathrsfs,psfrag}

\usepackage[dvipsnames]{xcolor}

\usepackage{mathtools}
\usepackage[all]{xy}
\usepackage[margin=1in,marginparwidth=0.8in, marginparsep=0.1in]{geometry}

\usepackage{comment}

\usepackage{booktabs,enumerate}
\usepackage{multirow}
\usepackage{mathrsfs}

\usepackage{tikz,tikz-3dplot}
\usetikzlibrary{matrix}

\usetikzlibrary{shapes.geometric,decorations.pathreplacing}

\usetikzlibrary{shapes.misc}

\tikzset{cross/.style={cross out, draw=black, minimum size=2*(#1-\pgflinewidth), inner sep=0pt, outer sep=0pt},
cross/.default={1pt}}

\usepackage{tikz-cd}

\usepackage{subcaption}

\usepackage{upgreek}

\usepackage[margin=1in,marginparwidth=0.8in, marginparsep=0.1in]{geometry}

\usepackage[bookmarks=true, bookmarksopen=true,%
bookmarksdepth=3,bookmarksopenlevel=2,%
colorlinks=true,%
linkcolor=blue,%
citecolor=blue,%
filecolor=blue,%
menucolor=blue,%
urlcolor=blue]{hyperref}

\newtheorem{dummy}{dummy}[section]
\newtheorem{lemma}[dummy]{Lemma}
\newtheorem{theorem}[dummy]{Theorem}

\newtheorem{conjecture}[dummy]{Conjecture}
\newtheorem{corollary}[dummy]{Corollary}
\newtheorem{proposition}[dummy]{Proposition}
\newtheorem{prop}[dummy]{Proposition}
\theoremstyle{definition}
\newtheorem{definition}[dummy]{Definition}

\newtheorem{example}[dummy]{Example}
\newtheorem{remark}[dummy]{Remark}

\newtheorem{Q}[dummy]{Question}

\newcommand{\bC}{\mathbb{C}}

\newcommand{\bN}{\mathbb{N}}
\newcommand{\bP}{\mathbb{P}}

\newcommand{\bR}{\mathbb{R}}
\newcommand{\bZ}{\mathbb{Z}}

\newcommand{\bc}{\mathbf{c}}

\newcommand{\cA}{\mathcal{A}}
\newcommand{\cC}{\mathcal{C}}
\newcommand{\cD}{\mathcal{D}}
\newcommand{\cE}{\mathcal{E}}

\newcommand{\cH}{\mathcal{H}}

\newcommand{\cL}{\mathcal{L}}
\newcommand{\cM}{\mathcal{M}}

\newcommand{\cO}{\mathcal{O}}
\newcommand{\cP}{\mathcal{P}}

\newcommand{\cT}{\mathcal{T}}

\newcommand{\cX}{\mathscr {X}}

\newcommand{\op}{\operatorname}

\newcommand{\N}{N}

\newcommand{\PGL}{\mathrm{PGL}}

\newcommand{\Hom}{\mathrm{Hom}}

\renewcommand{\log}{{\op{log}}}

\newcommand{\todo}[1]{{\marginpar{\tiny #1}}}

 \numberwithin{equation}{subsection}
\numberwithin{figure}{subsection}

\newcommand{\leg}{S}

\newcommand{\red}[1]{{\color{red}#1}}
\newcommand{\blue}[1]{{\color{blue}#1}}

\newcommand{\Exp}{\mathrm{Exp}}

\newcommand{\C}{\mathbb{C}}

\newcommand{\e}{\mathrm{exp}}
\newcommand{\E}{\mathrm{Exp}}
\renewcommand{\l}{\mathrm{log}}

\newcommand{\x}{\mathbf{x}}
\newcommand{\y}{\mathbf{y}}
\newcommand{\z}{\mathbf{z}}
\newcommand{\Ad}{\mathrm{Ad}}
\newcommand{\ad}{\mathrm{ad}}

\newcommand{\Sk}{\mathrm{Sk}}

\newcommand{\q}{q}

\newcommand{\fb}{\mathfrak{b}}
\newcommand{\fp}{\mathfrak{p}}
\renewcommand{\frame}[1]{\widetilde{P}_{#1}}
\newcommand{\skein}[1]{P_{#1}} 

\newcommand{\loo}[1]{{\tiny\text{\textcircled{$#1$}}}}

\newcommand{\qbinom}[2]{\begin{bsmallmatrix}{#1}\\ {#2}\end{bsmallmatrix}_q}

\newcommand{\seed}[1]{\underline{\mathbf{#1}}}

\newcommand{\link}{\begin{tikzpicture}[xscale=.1,yscale=.05]
\draw (0,-2.4)--(0,0);
\draw[ultra thick,white] (0,0) circle (1cm);
\draw (0,0) circle (1cm);
\draw[ultra thick,white] (0,0)--(0,2);
\draw (0,0)--(0,2);
\end{tikzpicture}}
 
\title[Skeins, Clusters and Wavefunctions]{Skeins, Clusters and Wavefunctions}

\author[Mingyuan Hu, Gus Schrader, and Eric Zaslow]{Mingyuan Hu, Gus Schrader, and Eric Zaslow\\
\\
{\tiny Department of Mathematics, Northwestern University}\\
\\
}

\begin{document}

\maketitle

\begin{abstract}
In \cite{SSZ}, the second- and third-named authors, with Linhui Shen, used cluster theory to
construct wavefunctions for branes in threespace and conjecturally
relate them to open
Gromov-Witten invariants.  
This was done by defining a quantum Lagrangian subvariety of a quantum cluster variety,
and mutating a simple solution to the defining equations from a distinguished seed.
In this paper,
we extend the construction to incorporate
the skein-theoretic approach to open Gromov-Witten theory of Ekholm-Shende \cite{ES19}.
In particular, we define a skein-theoretic version of cluster theory,
including the groupoid of seeds and mutations and a skein-theoretic version
of the quantum dilogarithm.  {We prove a pentagon relation in the skein of the closed torus in this context,
and give strong evidence that its analogue holds for arbitrary surfaces.}
We propose face relations satisfied by the skein-theoretic wavefunction, similar to \cite{SS},
prove their invariance under mutations, and show their solution is unique.
We define a skein version of framings in the story, and use the novel cluster
structure to compute wavefunctions in several examples.
The skein approach incorporates moduli spaces of sheaves of higher microlocal rank and their quantizations.
\end{abstract}

\setcounter{tocdepth}{1}
\tableofcontents

\section{Introduction and Summary}
\label{sec:intro}

 Several notable works have considered the three-dimensional physics
determined by a fivebrane in an M-theory compactification
that wraps a Lagrangian in a Calabi-Yau threefold.
The partition function of the theory defines a wavefunction
in an effective quantium-mechanical problem \cite{CEHRV,DGGo}.
In \cite{SSZ},
a class of theories defined by Lagrangians 
bounding Legendrian surfaces in the contact boundary of $\bC^3$ were considered.
Wavefunctions were constructed by exploiting cluster theory,
similar to \cite{CEHRV}.
The wavefunctions were shown to quantize the potential function defining
a classical Lagrangian moduli space inside a symplectic torus.
Following the pioneering work of Aganagic-Vafa \cite{AV},
the authors conjectured that the wavefunctions also
encode open Gromov-Witten (OGW) invariants with
one boundary component.

{In \cite{SSZ}, the wavefunction is locally a function  $\Psi(X_i;q)$ which is an element of a formal Laurent series ring acted upon by a quantum cluster torus. Its annihilator in this quantum torus quantizes the classical defining ideal of the moduli space of branes.

In \cite{ES19}, Ekholm and Shende construct a new
formalism for OGW theory in which
the boundary of a holomorphic map is recorded as an element of the HOMFLYPT skein module of $L$, rather than by the coarser information of its homology class.  
More specifically in our setting, if $L\subset \bC^3$ is Lagrangian brane, then one considers
a suitable moduli space $\cM$ of bare holomorphic maps $u: (\Sigma,\partial \Sigma)\to (\bC^3,L)$
and defines $\Psi(z,a) =  \sum_{u\in \cM} w(u) z^{-\chi(\Sigma)} a^{u\link L} u(\Sigma)\cdot [\partial u] \in \mathrm{Sk}(L),$
where $w$ is the weighted count of the point $u\in \cM$, $u\link L$ is a linking number, and $[\partial u]$
is an element of the skein module of $L$ --- see \cite[Definition 5.1]{ES19}
and subsequent sections for details.
It is tempting, therefore,
to extend the cluster-theoretic calculation schema of \cite{SSZ} to the new skein-valued formalism.
In \cite{SS}, the authors showed that the defining
equations for the moduli space in \cite{SSZ} had an evident extension
to the skein algebra of the Legendrian surface, providing a pathway for the sought-after extension.}

In this note, we merge these approaches
by considering 
a skein-theoretic approach to cluster theory.
{We introduce a skein-theoretic notion of mutation (Definition~\ref{def: Baxter operator}), defined using an analogue of the quantum dilogarithm function valued in the skein of the annulus.} The image of this skein dilogarithm under the rank-1 reduction to the linking skein of the annulus recovers the formal series of the standard $q$-dilogarithm.

We prove in Theorem~\ref{thm:face relation} that the defining relations for the skein-valued wavefunction
are compatible with skein mutation, and that (for the class of fillings considered in~\cite{SSZ}) they have a unique solution in the completed skein of $L$. 

The computational schema of \cite{SSZ} (including framings) is thereby extended to to the skein setting. We illustrate the utility of this skein cluster approach in several examples: 
in Proposition~\ref{prop:TV} we derive the topological vertex \cite{AKMV} solution for a single Aganagic-Vafa brane
in $\bC^3$, reproducing the analogous result from \cite{ES20} (while also incorporating
different framings), and in Example~\ref{eg:unknot-conormal} we show how to obtain the solution for the unknot-conormal using skein mutation. 

We also investigate some algebraic properties of the skein dilogarithm function. In Theorem~\ref{thm:baxter-pentagon} we show that the skein dilogarithms associated to the $(1,0)$- and $(0,1)$-curves on the closed torus $T^2$ satisfy a \emph{pentagon identity}, whose image in the linking skein recovers the celebrated one~\eqref{eq:ab-penta} satisfied by the quantum dilogarithm.

Counting curves in the HOMFLYPT skein  naturally incorporates
sheaves
of higher microlocal rank:
morally, the skein theory encodes
stacks of branes, meaning it should include
moduli of microlocal rank-$N$ sheaves for all
$N$. 
{In Section~\ref{sec:higher-rank}, we investigate this in  phenomenon in the fundamental case of a punctured genus-one surface $T^2\setminus D^2$, which captures the tubular neighborhood of two simple closed curves in a general surface intersecting at a single point. In this case we observe that the finite-rank skein algebra of $T^2\setminus D^2$ \emph{itself} admits a cluster structure in the standard (i.e., rank-1) sense.  Patching these local rank-1 cluster structures together, we show that the image of the skein dilogarithm under the reduction to the rank-$N$ skein factors into a product of standard quantum dilogarithm functions. We use this description to prove in Section~\ref{subsec: pentagon in finite rank} that the pentagon identity holds in all rank-$N$ specializations the skein of the punctured torus.}


\subsection*{Acknowledgements}

At an early stage of this project, we discussed several of these topics with Vivek Shende.
We also would like to thank Philippe di Francesco, Rinat Kedem, Peter Samuelson, Alexander Shapiro, and Boris Tsygan for helpful discussions.
G.S~has been supported by NSF grant
DMS-2302624.
E.Z.~has been supported by NSF grants
DMS-1708503 and DMS-2104087.

\section{Skeins}

Here we briefly recall some facts about the HOMFLYPT skein module, and define some notations which will be used later.

\subsection{The HOMFLYPT skein module}

The HOMFLYPT skein module $\Sk(M)$ of an oriented three-dimensional manifold with boundary $M$ is generated
by $R$-linear combinations of framed oriented links in $L$, up to isotopy,
modulo the relations

$$
\begin{tikzpicture}[scale=0.8]

\node (x) at (0,0){};
    \draw[thick,->] (x.45)-- (.5,.5);
    \draw[thick,->] (x.135) -- (-.5,.5);
    \draw[thick] (x.315) -- (.5,-.5);
    \draw[thick] (x.45) -- (-.5,-.5);

\node at (2,0) {$-$};

  \begin{scope}[xshift=4cm]
\node (y) at (0,0){};
  \draw[thick,->] (y.45)-- (.5,.5);
    \draw[thick,->] (y.135) -- (-.5,.5);
    \draw[thick] (y.225) -- (-.5,-.5);
    \draw[thick] (y.135) -- (.5,-.5);


\node at (2,0) {$=$};

  \end{scope}

  \begin{scope}[xshift=8cm]
  \node at (-.5,0) {$z$};
  \draw[thick, ->] (0,-.5) arc (-70:70:0.5);
  \draw[thick, ->] (1.2,-.5) arc (180+70:180-70:0.5);
  \end{scope}
\end{tikzpicture}
$$
\phantom{a}
$$
\begin{tikzpicture}

\node (x) at (0,0){};
\draw[thick] (x.45) to[out=45,in=180] (.3,.2);
\draw[thick] (.3,.2) to[out=0,in=0] (.3,-.2);
\draw[thick] (.3,-.2) to[out=180,in=315] (x.315);
\draw[thick] (x.45) -- (-.5,-.5);
\draw[thick,->] (x.135) -- (-.5,.5);
\node at (1.25,0) {$=$};
\node at (2.15,0) {$a$};
\draw[thick,->] (2.65,-.5)--(2.65,.5);
\begin{scope}[xshift = 5.5cm]
\node (y) at (0,0){};
\draw[thick] (y.45) to[out=45,in=180] (.3,.2);
\draw[thick] (.3,.2) to[out=0,in=0] (.3,-.2);
\draw[thick] (.3,-.2) to[out=180,in=315] (y.315);
\draw[thick,->] (y.315) -- (-.5,.5);
\draw[thick] (y.225) -- (-.5,-.5);
\node at (1.25,0) {$=$};
\node at (2.15,0) {$a^{-1}$};
\draw[thick,->] (2.75,-.5)--(2.75,.5);
    
\end{scope}

\begin{scope}[xshift=11cm]
    \draw[thick] (0,0) circle (.5);
\node at (1.25,0) {$=$};
\node at (2.5,0) {\Large $\frac{a-a^{-1}}{z}$};
\end{scope}

\end{tikzpicture}
$$
where $z = \q^{1/2} - \q^{-1/2}$
and $R = \bC[\q^{\pm 1/2},a^{\pm 1},(\q^{k/2}-\q^{-k/2})^{-1}]$ , where $k$ runs over $\bN$. 
\begin{example}
 It is well known that 
 $$
 \Sk(\bR^3) \simeq R \langle \emptyset \rangle. 
 $$
 A link is sent to its \emph{HOMFLYPT polynomial} under this identification. 
\end{example}

The skein algebra $\Sk(\Sigma)$ of an oriented surface $\Sigma$ is $\Sk(\Sigma\times I)$
endowed with the product $ab$ defined by placing the link $b$ above the link $a$:
embed $(0,1)\sqcup (0,1) \hookrightarrow (0,2)$ then rescale $(0,2) \cong (0,1)=:I.$

For a $3$-manifold $M$ with boundary $\partial M =  \Sigma$, the skein algebra $\Sk (\Sigma) $ acts naturally on the skein module $\Sk(M)$, given by gluing $\Sigma \times I$ to the boundary of $M$.

Following \cite{SS} and \cite{przytycki1998q}, if we furthermore quotient by the relation
$$
\begin{tikzpicture}
      \begin{scope}[xshift=4cm]
\node (y) at (0,0){};
  \draw[thick, ->] (y.45)-- (.5,.5);
    \draw[thick, ->] (y.135) -- (-.5,.5);
    \draw[thick] (y.225) -- (-.5,-.5);
    \draw[thick] (y.135) -- (.5,-.5);


\node at (2,0) {$=$};

  \end{scope}

  \begin{scope}[xshift=8cm]
  \node at (-.5,0) {$\q^{-1/2}$};
  \draw[thick, ->] (0,-.5) arc (-70:70:0.5);
  \draw[thick, ->] (1.2,-.5) arc (180+70:180-70:0.5);
  \end{scope}

\end{tikzpicture}
$$
and set
$$
a = \q^{1/2},
$$ 
we get the linking skein module $Lk (M)$, with a  natural quotient map:
\begin{equation}\label{map: link alg}
    \Sk(M) \longrightarrow Lk(M). 
\end{equation}
In the case where $M = \Sigma \times I$, the intersection pairing on $\Sigma$ defines a skew bilinear form on the lattice $H_1(M,\mathbb{Z})$, and the linking skein $Lk(M)$ is isomorphic to the quantum torus associated to this form (see \cite{przytycki1998q}).

\subsection{The skein module of the solid torus and the elliptic Hall algebra}
\label{sec:skt2}

Although a workable algebraic description of the HOMFLYPT skein of a general 3-manifold seems elusive, in the case $M\simeq D^2\times S^1 $ is a solid torus with $\partial M= T^2\coloneqq \bR^2/\bZ^2$ a very useful representation-theoretic interpretation has been found by Morton and Samuelson \cite{MS17}.
They prove that $\Sk(T^2)$ admits a presentation as an $R$-algebra
by certain elements $P_{\mathbf{x}},$ $\mathbf{x}\in \bZ^2$, with defining relations
\begin{align}
\label{eq:eha-rels}
    [P_\mathbf{x},P_\mathbf{y}] = \{d\}P_{\mathbf{x}+\mathbf{y}},
\end{align}    
where $d = \det(\mathbf{x}|\mathbf{y})$ and $\{d\} = \q^{d/2} - \q^{-d/2}.$
When the vector $\mathbf{x}\in\mathbb{Z}^2$ is primitive, the generator $P_\mathbf{x}$ corresponds to an embedded curve with homology class $\mathbf{x}\in H_1(T^2)\cong \bZ^2.$  
The relations \eqref{eq:eha-rels} imply that $\Sk(T^2)$ is generated over $R$ by the elements $P_{(\pm1,0)},P_{(0,\pm1)}$.

The above description identifies $\Sk(T^2)$ with the specialization of the \emph{elliptic Hall algebra} $\mathcal{E}_{\sigma,\overline\sigma}\otimes_\mathbb{C}\mathbb{C}[a^{\pm1}]$  at the parameter-values $\sigma=\overline\sigma=\q^{-\frac{1}{2}}$.

The structure of $\Sk(D^2\times S^1)$ as a module over $\Sk(T^2)$ depends on the choice of identification of $T^2$ with $\partial(D^2\times S^1)$; let us choose this identification so that the $(1,0)$-curve in $T^2$ is the isotopic to the boundary of a disk $D^2\times\{*\}$ in the solid torus. 

As explained in \cite{MS17}, in order to write explicit formulas for this action it is convenient to further choose an identification $D^2\simeq I\times I$, thereby fixing an isomorphism of the solid torus with the product of an interval $I$ and an annulus $I\times S^1$, so that the $(0,1)$-curve on the torus is isotopic to $\{*\}\times\{*\}\times S^1$  in the solid torus.

Under these identifications, the image of the skein of the annulus in $\Sk(T^2)$ is the commutative  `vertical' subalgebra in $\Sk(T^2)$ generated by the  $\{P_{(0,m)}\}_{m\in\mathbb{Z}}$ and $\Sk(D^2\times S^1)$ is a free module of rank 1 over this subalgebra. 

We may regard the vertical subalgebra as the tensor product $\mathcal{F}^+\otimes\mathcal{F}^-$ of two polynomial rings in infinitely many variables $\{P_{(0,\pm m )}\}_{m>0}$, and we will be mostly interested in the `positive' part of the skein module $\Sk^+(D^2\times S^1)$ given by the subspace $\mathcal{F}^+\otimes1$. In addition to being obviously preserved by the $P_{(0,m)}$ with $m>0$, the subspace $\Sk^+(D^2\times S^1)$ is also preserved by the action of all $P_{(m,0)}$ with $m>0$. The action of this positive part of the horizontal subalgebra is in fact simultaneously diagonalizable, with distinct eigenvalues: $\Sk^+(D^2\times S^1)$ admits a basis $W_\lambda$ indexed by integer partitions $\lambda$, on which the action of the algebra $\Sk^+(T^2)$ can be described combinatorially as follows. Recall that a \emph{border strip} is a connected skew partition containing no $2\times 2$ squares, and that the \emph{height} $\mathrm{ht}(\alpha)$ of such a border strip $\alpha$ is defined to be the number of rows in $\alpha$ plus 1. Then writing $\lambda+n$ for the set of all partitions obtained from $\lambda$ by adding a border strip with $n$ boxes, we have
\begin{align}
\label{eq:MN-rule}
P_{(0,n)}\cdot W_\lambda = \sum_{\mu\in \lambda+n} (-1)^{\mathrm{ht}(\mu-\lambda)}W_{\mu}.
\end{align}
Note that this formula takes the particularly simple form when $n=1$:
$$
P_{(0,1)}\cdot W_\lambda = \sum_{\mu= \lambda+\square} W_{\mu}.
$$
where the sum is taken over all partitions $\mu$ obtained from $\lambda$ by adding a single box.

To describe the action of the horizontal subalgebra, recall that the \emph{content} $c(x)$ of a box $x$ in row $i$ and column $j$ of a Young diagram is $c(x)=j-i$. Then the eigenvalues of the $P_{(m,0)}$ for $m>0$ are given by
\begin{align}
\label{eq:P-eigenvalues}
P_{(m,0)}\cdot W_\lambda = \left(\frac{a^m-a^{-m}}{\{m\}} + a^m\{m\}\sum_{x\in\lambda}\q^{mc(x)}\right)\cdot W_{\lambda}.
\end{align}
Note that in terms of the parts $\lambda = (\lambda_1\geq\lambda_2\geq\cdots\geq0)$ of the partition $\lambda$, we have
\begin{align}
\label{eq:content-sum}
    \sum_{x\in \lambda}\q^{c(x)} &= \sum_{k\geq1} \q^{-(k-1)}\frac{(1-\q^{\lambda_k})}{1-\q}
\end{align}

Introducing coordinates $l_k:=\lambda_k - k +1$ on the set of all partitions, we can write the latter formally as
$$
    \sum_{x\in \lambda}\q^{c(x)}= -\q^{-1}-\frac{1}{1-\q}\sum_{k\geq1} \q^{l_k}.
$$

When both $m>0$ and $n>0$, we have a hybrid formula
\begin{align}
P_{m,n}\cdot W_\lambda = a^{m}\frac{\{m\}}{\{mn\}}\sum_{\mu\in \lambda+n}(-1)^{\mathrm{ht}(\mu-\lambda)}\sum_{x\in {\mu - \lambda}}\q^{mc(x)}\cdot W_{\mu}.
\end{align}

The formulas for the action of the generators above can be interpreted in terms of symmetric functions as follows. Let us identify $\mathcal{F}^+=\Sk^+(I\times S^1)$ with the ring of symmetric functions in infinitely many variables $x_1,x_2,\ldots$ by identifying the basis element $W_\lambda$ with the Schur function $s_\lambda(x)$. Then under this identification, the element $P_{(0,n)}$ acts by multiplication by the $n$-th power sum symmetric function $p_n(x) = \sum_i x_i^n$, so that~\eqref{eq:MN-rule} corresponds to the Murnaghan-Nakayama rule.

In Theorem 3.1 of~\cite{SV11} it is shown that for each $N\geq 1$, the elliptic Hall algebra $\mathcal{E}_{q,t}$ has a surjective map to the \emph{spherical double affine Hecke algebra} $SH_{q,t}(\mathfrak{gl}_N)$ of $\mathfrak{gl}_N$ introduced by Cherednik~\cite{C05}. 
 This quotient will be relevant to our discussion of moduli spaces of sheaves of finite microlocal rank in Section~\ref{sec:higher-rank}.

This identification of the skein of the annulus with the ring of symmetric functions allows us to define \emph{Adams operations} at the level of skeins. Indeed, recall that the ring of symmetric functions is the free $\lambda$-ring generated by a single line element $e_1=p_1=\sum_i x_i$, which corresponds to the first elementary symmetric function; the other elementary symmetric functions are given by $e_n = \lambda^n(e_1)$. On the other hand, the $n$-th Adams operation $\psi^n$ takes the line element $p_1$ to the $n$-th power sum $p_n$, and thus at the level of skeins sends $P_{(0,1)}$ to $P_{(0,n)}$.



\subsection{Skein modules with framing lines and base points}\label{sec:skeins with b and f}
Following \cite{SS}, for a oriented $3$-manifold $M$ we can choose a \textit{framing line} $\mathfrak{l}$, which is a 1-chain in $M$. We write $Sk (M; \mathfrak{l})$ for the skein of $M \backslash \mathfrak{l}$, modulo an additional relation:
\begin{equation}\label{framing line}
\begin{tikzpicture}[scale = 1.5]

\node (x) at (0,0){};
    \draw[thick,->] (x.45)-- (.5,.5);
    \draw[thick,->] (x.135) -- (-.5,.5);
    \draw[thick] (x.315) -- (.5,-.5);
    \draw[thick] (x.45) -- (-.5,-.5);
    
\node at (-0.7, 0.5) {$\mathfrak{l}$};
\node at (0.7, 0.5) {$L$};
    
\node at (1.2,0) {$= $};

\node at (2, 0) {$ (-a)$};

  \begin{scope}[xshift= 3 cm]
\node (y) at (0,0){};
  \draw[thick,->] (y.45)-- (.5,.5);
    \draw[thick,->] (y.135) -- (-.5,.5);
    \draw[thick] (y.225) -- (-.5,-.5);
    \draw[thick] (y.135) -- (.5,-.5);

    \node at (-0.7, 0.5) {$\mathfrak{l}$};
    \node at (0.7, 0.5) {$L$};
  
  \end{scope}

\end{tikzpicture} 
\end{equation}
i.e. every time a link $L$ passes through $\mathfrak{l}$, it is multiplied by $(-a)^{\pm }$, where the sign depends on the orientation. 

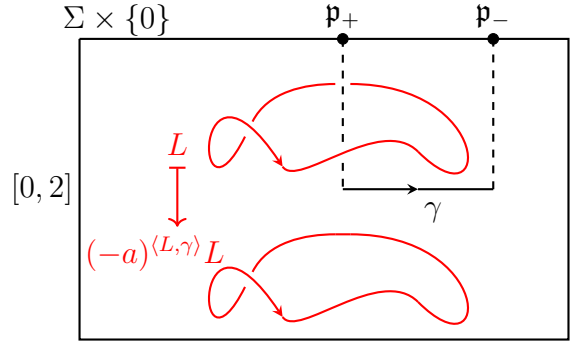
\begin{figure}[htpb]
	\centering
	\begin{tikzpicture}[scale=1]
		\filldraw (2, 0) circle(2pt) node[anchor = south] {$\fp_+$} ; 
		\filldraw (4, 0) circle(2pt) node[anchor = south] {$\fp_-$} ; 
		
		\node at (-1, .25) {$\Sigma \times \{ 0 \}$}; 
		\node at (-2, -2) {$[0, 2]$};
		
		\draw[thick] (-1.5, 0) -- (5, 0) -- (5, -4) -- (-1.5, -4) -- (-1.5, 0);
		\draw[thick, dashed] (2, 0) -- (2, -1.4);
		\draw[thick, dashed] (2, -1.65) -- (2, -2);
		\draw[thick, -stealth] (2, -2) -- (3, -2);
		\draw[thick] (3, -2) -- (4, -2);
		\draw[thick, dashed] (4, -2) -- (4, 0);
		
		\draw[red, thick] (1.9, -0.6) .. controls (1.5, -.6)and (1, -.7).. (.8, -1.1);
		\draw[red, thick, -stealth] (.7, -1.3) .. controls(0, -2.7) and (0, .1).. (1.2, -1.7);
		\draw[red, thick] (1.2, -1.7) .. controls (1.4, -2) and (2.5, -1).. (3, -1.5) ;
		
		\draw[red, thick] (3, -1.5) .. controls (4, -2.5) and (4, -0.6) .. (2.1, -.6);
		
		\begin{scope}[yshift= -2cm]
		
		\draw[red, thick] (1.9, -0.6) .. controls (1.5, -.6)and (1, -.7).. (.8, -1.1);
		\draw[red, thick, -stealth] (.7, -1.3) .. controls(0, -2.7) and (0, .1).. (1.2, -1.7);
		\draw[red, thick] (1.2, -1.7) .. controls (1.4, -2) and (2.5, -1).. (3, -1.5) ;
		
		\draw[red, thick] (3, -1.5) .. controls (4, -2.5) and (4, -0.6) .. (2.1, -.6);
		
		\draw[red, thick] (1.9, -0.6) -- (2.1, -.6);
		\end{scope}
		
		\node at (3.2, -2.3) {$\gamma$};
			
		\node at (-.2, -1.4) {\color{red}$L$};
		\draw[red, thick, |->] (-.2, -1.7) -- (-.2, -2.5);
		
		\node at (-.5, -2.85) {\color{red}$(-a)^{\langle L, \gamma \rangle}L$ }; 
		
	\end{tikzpicture}
	\caption{Pushing a link from top to bottom}
	\label{fig: downward}
\end{figure}

For a surface $\Sigma$ and a 0-chain $\fp$, we define $\Sk(\Sigma; \fp) \coloneqq \Sk(\Sigma \times [0, 1]; \fp \times [0,1])$. If $\fp$ is exact, we can choose $\fp = \partial \gamma$. Then there is an isomorphism of algebras:
    \begin{align}\label{iso of skeins}
        \Sk(\Sigma; \fp) &\longrightarrow \Sk(\Sigma)  \\
        L &\longmapsto (-a)^{\langle L, \gamma \rangle} L \notag
    \end{align}
Here $\langle L, \gamma \rangle $ means that we first push $L$ into $\Sigma$ to get an immersed curve, and then take the intersection number with $\gamma$. 
 Geometrically, if we consider $\Sigma \times [0, 2]$ with the $1$-chain $\fp \times [0, 1] \cup \gamma\times \{1\} $,  this map is given by pushing the tangles from the top to the bottom, and multiplying by $(-a)^{\pm}$ (the sign given by \eqref{framing line}) every time $L$ passes through $\gamma$. (See Figure \ref{fig: downward})

We will also need the notion of skein modules with base points. Let $\fb = \fb_+ - \fb_-$ be a framed $0$-chain on $\partial M$, i.e. each point is equipped with a choice of tangent vector. Define $\Sk(M, \fb)$ to be the skein module consisting of strands going from $\fb_+$ to $\fb_-$. If $\fb$ is exact and we choose a \emph{capping path} (or paths) $C$ on $\partial M$ connecting $\fb_-$ to $\fb_+$ (i.e.~$\partial C = -\fb)$, then we get a map of skein modules:
\begin{equation}\label{base points}
  \Sk (M, \fb) \longrightarrow \Sk(M) 
\end{equation}
 by capping a tangle with $C$.

For a surface $\Sigma$ with $0$-chain $\fb$, we write $\Sk (\Sigma, \fb)$ for $\Sk(\Sigma \times I, \fb \times \{0\})$. Let $M$ be a $3$-manifold with $\partial M = \Sigma$. 
Then the action depicted in Figure~\ref{fig: action} makes $\Sk(M, \fb)$ into a module over the algebra $\Sk (\Sigma, \fb)$:
\begin{equation}\label{action}
    \Sk(\Sigma,\fb) \otimes \Sk(M) \longrightarrow \Sk(M, \fb) 
\end{equation}

\begin{figure}[htpb]
	\centering
	\begin{subfigure}[b]{0.35\textwidth}
	\centering
	\begin{tikzpicture}[scale = 1]
		\node at (-.6, -.5) {$\Sigma \times I$}; 
		\node at (-.5, -3) {$M$};

		\draw[thick] (0, 0) -- (3, 0) -- (3, -1.2) -- (0, -1.2) -- (0, 0);
		\draw[thick] (0, -1.8) -- (3, -1.8);
		\draw[thick] (0, -1.8) .. controls (0, -5) and (3, -5) .. (3, -1.8);
		
		\draw[thick, red] (0.3, 0) .. controls (2, -1.7) and (-.7, -.6) .. (.7, -.6);
		\draw[thick, red, ->] (.9, -.6) .. controls (1.5, -.6) and (1.5, -1.6)..  (1.8, -.7);
		\draw[thick, red] (1.8, -0.7)  .. controls (2.1, 0.2) and (3.2, -1.5) .. (2.4, 0);
		
		\draw[thick] (1.4, -2.55) .. controls (1.1, -2.8) and (1.1, -3.4).. (1.4, -3.5);
		\draw[thick] (1.4, -2.7) .. controls (1.7, -2.9) and (1.6, -3.2) .. (1.3, -3.3);
		
		\draw[thick, red, rotate around={50: (1.5, -3)} ] (1.5, -3) ellipse (0.9 and 0.7);
		
		\filldraw (.3, 0) circle(1.5pt) node[anchor = south] {$\fb_+$}; 
		\filldraw (2.4, 0) circle(1.5pt) node[anchor = south] {$\fb_-$}; 
	\end{tikzpicture}
	\end{subfigure}
	\begin{tikzpicture}
		\node at (0,0) {$\xymatrix{{}\ar[rr]&&{} \\{} \\ {} \\ {}}$};
	\end{tikzpicture}
	\begin{subfigure}[b]{0.4 \textwidth}
		\begin{tikzpicture}	
			\node at (3.5, -1) {$M$};

			\draw[thick] (0, 0) -- (3, 0);
			\draw[thick] (0, 0) .. controls (0, -5) and (3, -5) .. (3, 0) ;
			
			\draw[thick] (1.4, -2.25) .. controls (1.1, -2.5) and (1.1, -3.1).. (1.4, -3.2);
			\draw[thick] (1.4, -2.4) .. controls (1.7, -2.6) and (1.6, -2.9) .. (1.3, -3);
			
			\draw[thick, red, rotate around={50: (1.5, -2.6)} ] (1.5, -2.6) ellipse (0.9 and 0.7);
			
			\draw[thick, red] (0.3, 0) .. controls (2, -1.7) and (-.7, -.6) .. (.7, -.6);
			\draw[thick, red, ->] (.9, -.6) .. controls (1.5, -.6) and (1.5, -1.6)..  (1.8, -.7);
			\draw[thick, red] (1.8, -0.7)  .. controls (2.1, 0.2) and (3.2, -1.5) .. (2.4, 0);
			
			\filldraw (.3, 0) circle(1.5pt) node[anchor = south] {$\fb_+$}; 
			\filldraw (2.4, 0) circle(1.5pt) node[anchor = south] {$\fb_-$}; 
		\end{tikzpicture}
	\end{subfigure}
	\caption{}
	\label{fig: action}
\end{figure}

For a surface $\Sigma$ with framing points $\fp$ and base points $\fb$, we can 
consider both of the above constructions together,
and define $\Sk(\Sigma, \fb; \fp)$. Fixing $\fp = \partial \gamma$ and capping path $C$ connecting $\fb_+$ to $\fb_-$ gives us a commutative diagram:
$$
\xymatrix{
\Sk(\Sigma, \fb; \fp) \ar[d]^\gamma \ar[r]^C & \Sk(\Sigma; \fp) \ar[d]^\gamma\\
\Sk(\Sigma, \fb) \ar[r]^C &\Sk(\Sigma)
}
$$


\subsection{Completion}
The skein module $\Sk (M)$ is graded
by the group $H_1 (M; \mathbb{Z})$ . If we choose a homomorphism 
$$
\delta: \quad H_1 (M; \bZ) \longrightarrow \bZ
$$
then we can define the completion:
$$
\widehat{\Sk} (M) \coloneqq \varprojlim_{n\ge 0} \bigoplus_{ t(\lambda) \le n}  \Sk (M)_\lambda .
$$

For a surface $\Sigma$, the completion $\widehat{\Sk}(\Sigma):=
\widehat{\Sk}(\Sigma \times I)$ retains the structure of an algebra.  Moreover, if $\partial M = \Sigma$ and we choose compatible homomorphisms to $\bZ$,
i.e.~commuting with the natural map $H_1(\Sigma; \bZ) \rightarrow H_1 (M; \bZ)$, then $\widehat{\Sk}(M)$ inherits
a module structure over $\widehat {\Sk} (\Sigma)$. 

\begin{example}
\label{eg:completion}
    We will focus on the genus-$g$ handlebody $\cH_g$ with boundary $\Sigma_g$ --- see Figure \ref{fig:handlebody}.
    After choosing a set of standard circles $a_i$ and $b_i$, we get a basis for $H_1 ( \Sigma_g ; \bZ) $. The images of the $[a_i]$ form a basis for $H_1 (\cH_g)$. We then define the homomorphism for $\cH_g$ to be
\begin{align*}
    \delta: H_1 (\cH_g; \bZ) &\longrightarrow \bZ \\
    \sum_{i=1}^g c_i [a_i] &\longmapsto \sum_{i=1}^g c_i,
\end{align*}  
    and for $\Sigma_g$ to be $\delta$ composed with the natural map $H_1(\Sigma_g; \bZ) \rightarrow H_1( \cH_g ; \bZ)$. The corresponding completion $\widehat{\Sk} (\cH_g) $ is a well-defined module over the algebra $\widehat{\Sk} (\Sigma_g)$.

\end{example}
We say a loop is \textit{positive}/\textit{negative} with respect to $\delta$ if the image of its homology class under $\delta$ is positive/negative.

\section{The Chromatic Lagrangian}
\label{sec:chromatic}

Our constructions center around
the geometric and algebraic objects
spawned by a cubic planar graph $\Gamma$ on the 2-sphere, or -- dually -- by a \emph{triangulation} of the sphere. 
In this paper, we will always consider cubic graphs up to the action of $\mathrm{Homeo}^+_0(S^2;\mathrm{Vert}(\Gamma))$, the identity component of the group of orientation-preserving homeomorphisms of $S^2$ which fix each vertex of $\Gamma$.
\subsection{Basic Constructions from a Cubic Planar Graph}

To each $\Gamma$ we can assign an integer $g$, the ``genus,'' such that $\Gamma$ has ${\sf v} = 2g+2$ vertices, ${\sf e} = 3g+3$ edges, and ${\sf f} = g+3$ faces.  As in \cite{TZ} and \cite{SSZ}, one may associate to
$\Gamma$ the following objects.  Below we  include some straightforward
generalizations of the
results of the cited papers to higher microlocal rank.  Other parts of this paper address areas where generalizations are subtle.

\begin{enumerate}
\item A Legendrian surface $\leg_{\Gamma} \subset T^{\infty} \bR^3 \subset S^5$ of genus $g$ \cite[Def. 2.1]{TZ}.
The surface $\leg_\Gamma$ is a branched double cover of $S^2,$ branched over the vertices of $\Gamma.$
It is defined by its front projection, which is taken to be a two-sheeted cover of $S^2$ with crossing
locus over the edges of $\Gamma$ and looking like the following near vertices:
\begin{figure}[H]
\includegraphics[scale = .25]{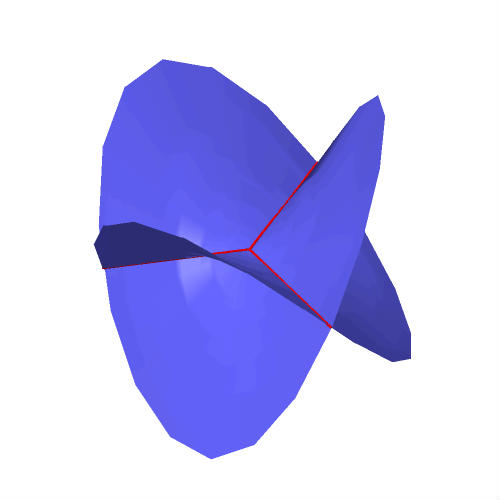}
\caption{The front projection of $\leg_\Gamma$ near a vertex.}
\label{fig:legprojs}
\end{figure}

\item \label{item:category} A category
$\cC_{\Gamma}$ of constructible sheaves on $\bR^3$ whose
singular support lies in $S_{\Gamma}$ \cite[\S 4.3]{TZ}.
By \cite{N,NZ} this is equivalent to a Fukaya category
whose geometric objects are Lagrangians
whose boundary-at-infinity is $S_\Gamma$.

\item A microlocal monodromy functor
$\cC_\Gamma \to \cL_\Gamma,$
where $\cL_\Gamma$ is
the category of local systems on $S_\Gamma$ minus the set ${\bf B}$ of branch points,
having monodromy $-\mathrm{id}$ around the points of ${\bf B}$.
Note that $\cL_\Gamma$ is a torsor for $\mathrm{Loc}(S_\Gamma),$ the
category of local systems on $S_\Gamma.$

\item \label{item:chromatic} 
For each
natural number $N$, a subcategory $\cC_{\Gamma,N}$ and 
moduli space $\cM_{\Gamma,N}$ of microlocal-rank-$N$ constructible sheaves on $\bR^3$ whose singular support lies in $S_{\Gamma}$  \cite[\S 4.3]{TZ},
and a moduli space $\cP_{\Gamma,N}$ of rank-$N$ local systems in $\cL_\Gamma$ which
we call the period domain.
More concretely, $\cM_{\Gamma,N}$ is the space of $\PGL_{2N}$-equivalence classes of colorings of the faces of $\Gamma$ by $N$-planes in $Gr(N,2N)$ such that neighboring faces are colored by transverse $N$-planes.  When $N=1,$
$\cM_{\Gamma,1}$ is the space of map colorings of $\Gamma$ by points in $\bP^1$
up to simultaneous action of $\PGL_2$, and
$\cP_{\Gamma,1}$ is non-canonically isomorphic to
$H^1(\leg_{\Gamma},\bC^*)$, which is an algebraic torus equipped with an algebraic symplectic form coming from the intersection pairing on $H_1(S_{\Gamma})$ \cite[\S 4.6]{TZ}. 
Note that $H^1(\leg_\Gamma, \bC^*)$ can be identified with the moduli space of flat line bundles over $S_\Gamma$. It acts on $\cP_\Gamma$ by taking the tensor product of corresponding line bundles, and this action equips $\cP_\Gamma$ with the structure of an $H^1(\leg_\Gamma, \bC^*)$ torsor. 

\item
\label{item-edge-fn} The microlocal monodromy
functor descends to a map $\cM_{\Gamma,N} \to \cP_{\Gamma,N}$ on moduli \cite[\S 4.7]{TZ}.
It can be described as follows.  
Every edge $e$ of $\Gamma$ connects branch points and therefore defines an element of $H_1(S_\Gamma).$ 
It gives rise to a holonomy $X_e: \cP_{\Gamma,N} \to GL_N/\mathrm{Ad}_{GL_N}$ taking values in the adjoint quotient for $GL_N$.
(Edges surround two points, so the holonomies in
$\leg \setminus {\bf B}$ are the same as those in $\leg$.)
To define $X_e$,
first label the four vector spaces around a face as such:

\begin{center}
\begin{tikzpicture}
\node at (0,.65) {$a$};
\node at (1.75,0) {$d$};
\node at (0,-.65) {$c$};
\node at (-1.75,0) {$b$};
\draw[very thick] (-2,1)--(-1,0)--(1,0)--(2,1);
\draw[very thick] (-2,-1)--(-1,0);
\draw[very thick] (2,-1)--(1,0);
\end{tikzpicture}
\end{center}
Neighboring faces are transverse in the vector space $V$
which is the stalk of the sheaf in the simply connected
region above the Legendrian front.  As a result, we have
equivalences
$$a\cong V/b \cong c \cong V/d \cong a.$$
Writing $\phi$ for the composition map $a\to a$,
we define $X_e = -\phi.$
\end{enumerate}


\subsection{Globalization via Clusters in
Microlocal Rank 1}
\label{subsec:globalization}

In microlocal rank 1, it was shown in \cite{SSZ} that the embeddings $\mathcal{M}_{\Gamma,1}\rightarrow\mathcal{P}_{\Gamma,1}$ described above associated to \emph{all} cubic graphs of a given genus $g$ assemble to define a global Lagrangian inside a holomorphic symplectic variety. This global symplectic variety is a symplectic leaf inside the Poisson moduli space $\cX_{PGL_2, S^2_{g+3}}$ of \emph{framed} $PGL_{2}$-\emph{local systems} on a sphere with $g+3$ punctures (which we regard as sitting inside the faces of $\Gamma$), where the framing data for such a local system consists of the choice near each puncture $p$ of a flat section of the associated $\mathbb{P}^1$-bundle on a small loop around $p$.  The moduli space $\cX_{PGL_2, S^2_{g+3}}$ was shown by Fock and Goncharov \cite{FG1} to carry a \emph{cluster Poisson structure}. Let us briefly recall the meaning of this structure in the present context, and its role in the constructions of \cite{SSZ}. 

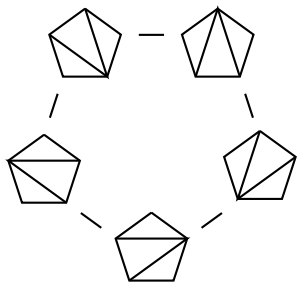
\begin{figure}[htpb]
    \centering
\begin{tikzpicture}[every node/.style={inner sep=0.5, thick, circle}, x=0.5cm, y=0.5cm, thick,scale=.5]

\def\n{2}

\begin{scope}[shift={(54:3*\n)}]
	\draw (90:\n) to (162:\n) to (-126:\n) to (-54:\n) to (18:\n) to (90:\n);
	\draw (-126:\n) to (90:\n) to (-54:\n);
\end{scope}

\begin{scope}[shift={(126:3*\n)}]
	\draw (90:\n) to (162:\n) to (-126:\n) to (-54:\n) to (18:\n) to (90:\n);
	\draw (90:\n) to (-54:\n) to (162:\n);
\end{scope}

\begin{scope}[shift={(198:3*\n)}]
	\draw (90:\n) to (162:\n) to (-126:\n) to (-54:\n) to (18:\n) to (90:\n);
	\draw (-54:\n) to (162:\n) to (18:\n);
\end{scope}

\begin{scope}[shift={(270:3*\n)}]
	\draw (90:\n) to (162:\n) to (-126:\n) to (-54:\n) to (18:\n) to (90:\n);
	\draw (162:\n) to (18:\n) to (-126:\n);
\end{scope}

\begin{scope}[shift={(344:3*\n)}]
	\draw (90:\n) to (162:\n) to (-126:\n) to (-54:\n) to (18:\n) to (90:\n);
	\draw (90:\n) to (-126:\n) to (18:\n);
\end{scope}

	\draw[-] (83:\n*11/4) to (97:\n*11/4);
	\draw[-] (155:\n*11/4) to (169:\n*11/4);
	\draw[-] (227:\n*11/4) to (241:\n*11/4);
	\draw[-] (299:\n*11/4) to (313:\n*11/4);
	\draw[-] (11:\n*11/4) to (25:\n*11/4);

\end{tikzpicture}

    \caption{Pentagon relation in the cluster modular groupoid for $\cX_{\PGL_2,S}$.}
    \label{fig:pentagon-rel}
\end{figure}

\begin{itemize}
    \item 
For each isotopy class $\Gamma$ of genus-$g$ cubic graph on $S^2$, we consider the lattice of edges $\Lambda_{\Gamma}\simeq \mathbb{Z}^{E_\Gamma}$, and the corresponding algebraic torus $\mathrm{Hom}_{\mathbb{Z}}(\Lambda_{\Gamma},\mathbb{C}^*)$. The lattice $\mathbb{Z}^{E_\Gamma}$ carries a skew-symmetric bilinear form whose value on pairs of basis vectors $e_i,e_j$ associated to edges is defined to be 
$$
(e_i,e_j)=\sum_{v\in\mathrm{Vertices}(\Gamma)}\epsilon_v(e_i,e_j),
$$
where $\epsilon_v(e_i,e_j)=1 $ (resp. $-1$) if edge $e_j$ is the predecessor (resp. successor) to edge $e_i$ with respect to the clockwise cyclic order on the set of edges incident to $v$, and zero otherwise. 
\item The skew bilinear form induces a log-canonical Poisson structure on 
the torus $\mathrm{Hom}_{\mathbb{Z}}(\Lambda_{\Gamma},\mathbb{C}^*)$, where the Poisson brackets between the canonical edge functions $X_e$ are defined by
$$
\{X_{e_i},X_{e_j}\} = (e_i,e_j)X_{e_i}X_{e_j}.
$$
\item For each isotopy class of genus-$g$ graph $\Gamma$, there is a birational Poisson map 
$$
\pi = (\pi^*X_{e}) ~:~\cX_{PGL_2, S^2_{g+3}} \longrightarrow (\mathbb{C}^*)^{E_\Gamma}
$$
whose components $\pi^*X_{e},e\in E_\Gamma$ are defined as follows. Dual to the graph $\Gamma$ is a triangulation of $S^2$ with vertices at the punctures, and each edge of $\Gamma$ is dual to the diagonal of a unique quadrilateral  in this triangulation. Since the quadrilateral is contractible, there is  well-defined rational functon $X_e$ obtained by taking the four flat sections $a,b,c,d$ near each corner, transporting to them to a common fiber of the $\mathbb{P}^1$-bundle, and computing their cross-ratio 
\begin{equation}
\label{eq:cross-ratio}
X_e = -\frac{a - b}{b-c}\cdot \frac{c-d}{d-a}.
\end{equation}
\item We say that two cubic graphs $\Gamma,\Gamma'$ are related by \emph{mutation at edge} $e_0$ if one is obtained from the other by performing the local operation shown in Figure~\ref{fig:q-mut}. At the level of the dual triangulations, mutation corresponds to flipping the diagonal in the quadrilateral containing the edge dual to $e_0$. The induced birational isomorphism between the tori $\mathrm{Hom}_{\mathbb{Z}}(\Lambda_{\Gamma'},\mathbb{C}^*)\dashrightarrow \mathrm{Hom}_{\mathbb{Z}}(\Lambda_{\Gamma},\mathbb{C}^*)$ is given by
\[
X_{e_i}'\longmapsto \left\{  \begin{array}{ll} 
      X_{e_0}^{-1} & \mbox{if } e_i=e_0, \\
      X_{e_i}(1+X_{e_0}^{-{\rm sgn}(e_i,e_0)})^{-(e_i,e_0)} & \mbox{if } e_i\neq e_0. \\
   \end{array} \right.
\]
This formula coincides with the $q=1$ specialization of that indicated in Figure~\eqref{fig:q-mut}.

\item Associated to a pair of graphs $\Gamma,\Gamma'$ related by the mutation at edge $e_0$ is a pair of lattice isometries $\nu_\pm:\Lambda_{\Gamma'}\rightarrow\Lambda_\Gamma$ called the \emph{positive/negative lattice mutations}, and defined by the formulas in Figure~\ref{eq:torusmutation}. 

\begin{equation}
\label{eq:torusmutation}
\begin{tikzpicture}
\draw[thick] (-1.5,-1)--(-1,0)--(-1.5,1);
\draw[thick] (-1,0)--(0,0);
\draw[thick] (.5,-1)--(0,0)--(.5,1);
\node at (-.5,1.6) {\small$\Gamma$};
\node at (-.5,.3) {\small$X_{e_0}$};
\node at (-1.6,.6) {\small$X_{e_1}$};
\node at (.7,.6) {\small$X_{e_4}$};
\node at (-1.7,-.6) {\small$X_{e_2}$};
\node at (.7,-.6){\small$X_{e_3}$};
\end{tikzpicture}
\qquad
\begin{tikzpicture}
\node at (0,0) {$\xymatrix{{}\\ {}&&\ar[ll]^{\nu^+_0}{}\\{}}$};
\end{tikzpicture}
\qquad
\begin{tikzpicture}

\node at (.15,1.6) {$ \Gamma'$};

\draw[thick] (-1,1)--(0,.5)--(1,1);
\draw[thick] (0,.5)--(0,-.5);
\draw[thick] (-1,-1)--(0,-.5)--(1,-1);
\node at (.55,0) {\small$ X_{-e_0}$};
\node at (-.8,.6) {\small$X_{e_1}$};
\node at (1.2,.6) {\small$X_{e_4+e_0}$};
\node at (-1.15,-.6) {\small$X_{e_2+e_0}$};
\node at (.8,-.6){\small$X_{e_3}$};
\end{tikzpicture}
\end{equation}

\begin{center}
\begin{tikzpicture}
\draw[thick] (-1.5,-1)--(-1,0)--(-1.5,1);
\draw[thick] (-1,0)--(0,0);
\draw[thick] (.5,-1)--(0,0)--(.5,1);
\node at (-.45,.25) {\small$X_{e_0}$};
\node at (-1.6,.6) {\small$X_{e_1}$};
\node at (.7,.6) {\small$X_{e_4}$};
\node at (-1.9,-.6) {\small$X_{e_2}$};
\node at (.7,-.6){\small$X_{e_3}$};
\end{tikzpicture}
\qquad
\begin{tikzpicture}
\node at (0,0) {$\xymatrix{{}\\ {}&&\ar[ll]^{\nu^-_0}{}\\{}}$};
\end{tikzpicture}
\qquad
\begin{tikzpicture}
\draw[thick] (-1,1)--(0,.5)--(1,1);
\draw[thick] (0,.5)--(0,-.5);
\draw[thick] (-1,-1)--(0,-.5)--(1,-1);
\node at (.55,0) {\small$X_{-e_0}$};
\node at (-1.05,.6) {\small$X_{e_1+e_0}$};
\node at (.9,.6) {\small$X_{e_4}$};
\node at (-.9,-.6) {\small$X_{e_2}$};
\node at (1,-.6){\small$X_{e_3+e_0}$};
\end{tikzpicture}
\end{center}

\item For each cubic graph $\Gamma$, the Poisson algebra $\mathcal{T}_\Gamma:= \Lambda_\Gamma\otimes_{\mathbb{Z}}\mathbb{C}^*$ of regular functions on the torus $\mathrm{Hom}_{\mathbb{Z}}(\Lambda_{\Gamma},\mathbb{C}^*)$ admits a canonical quantization to a \emph{quantum torus algebra} $\mathcal{T}^q_\Gamma$, the free $\mathbb{Z}[q^{\pm \frac{1}{2}}]$-module generated by $\{X_\lambda | \lambda\in \Lambda_\Gamma\}$ with multiplication given by\footnote{To make our conventions more compatible with the literature, we use $q^{1/2}$, whereas $q$ was used in \cite{SSZ}.}
$$
X_\lambda\cdot X_\mu = q^{\frac{1}{2}(\lambda,\mu)}X_{\lambda+\mu}.
$$
To a mutation at edge $e_0$ we associate the isomorphism $\mu^q_e$ between their noncommutative fraction fields indicated in Figure~\ref{fig:q-mut}.

\begin{center}
\begin{equation}
\label{fig:q-mut}
\begin{tikzpicture}

\draw[thick] (-1.5,-1)--(-1,0)--(-1.5,1);
\draw[thick] (-1,0)--(0,0);
\draw[thick] (.5,-1)--(0,0)--(.5,1);
\node at (-.45,.25) {\small$X_{e_0}$};
\node at (-1.6,.6) {\small$X_{e_1}$};
\node at (.7,.6) {\small$X_{e_4}$};
\node at (-1.9,-.6) {\small$X_{e_2}$};
\node at (.7,-.6){\small$X_{e_3}$};
\end{tikzpicture}
\qquad
\begin{tikzpicture}
\node at (0,0) {$\xymatrix{{}\\ {}&&\ar[ll]^{\mu_0}{}\\{}}$};
\end{tikzpicture}
\quad 
\begin{tikzpicture}
\draw[thick] (-2,1)--(-1,.5)--(0,1);
\draw[thick] (-1,.5)--(-1,-.5);
\draw[thick] (-2,-1)--(-1,-.5)--(0,-1);
\node at (-.9,1.6) {$\Gamma'$};
\node at (-.5,0) {\small$X_{-e_0}$};
\node at (-2.8,.6) {\small$X_{e_1}(1+q^{\frac{1}{2}}X_{e_0})$};
\node at (1.1,.6) {\small$X_{e_4}(1+q^\frac{1}{2}X_{-e_0})^{-1}$};
\node at (-2.85,-.6) {\small$X_{e_2}(1+q^\frac{1}{2}X_{-e_0})^{-1}$};
\node at (.95,-.65){\small$X_{e_3}(1+q^\frac{1}{2}X_{e_0})$};
\end{tikzpicture}
\end{equation}
\qquad
\end{center}

\item The quantum mutation automorphism $\mu^q_e: \mathcal{T}^q_{\Gamma'}\dashrightarrow \mathcal{T}^q_{\Gamma}$ associated to the flip at $e$ can be factored in two ways:
\begin{align*}
\mu_{e} &= \mathrm{Ad}_{\Phi(X_{e})}\circ\nu^+_e\\
&= \mathrm{Ad}_{\Phi(X_{-e})^{-1}}\circ\nu^-_e.
\end{align*}
Here $\nu^\pm_e:\mathcal{T}^q_{\Gamma'}\rightarrow \mathcal{T}^q_{\Gamma}$ refer to the isomorphisms of quantum tori induced by the corresponding lattice isometries, and $\mathrm{Ad}_{\Phi(X)}$ is the automorphism of $\mathcal{T}^q_{\Gamma}$ given by conjugation by the following element of the completed quantum torus:
\begin{align}
\label{eq:qdl-def}
\Phi(X)= \prod_{n=0}^\infty (1+q^{n+\frac{1}{2}}X)^{-1}
\end{align}
The formal power series $\Phi(X)$ is a version of the $q$-dilogarithm function. Indeed, we have
\begin{align*}
 \Phi(x) &= \exp\left(\sum_{m\geq1}\frac{(-q^\frac{1}{2}x)^m}{(1-q^{m})m}\right)\\
&=\exp\left(\frac{1}{q^\frac{1}{2}-q^{-\frac{1}{2}}}\sum_{m\geq1}\frac{-(-x)^m}{[m]_{q}m}\right)
\end{align*}
where $[m]_{q} = q^{\frac{m-1}{2}} + q^{\frac{m-3}{2}} + \cdots + q^{\frac{-m-1}{2}} = \frac{q^{m/2}-q^{-m/2}}{q^{1/2}-q^{-1/2}}$. It can be also interpreted via the \emph{plethystic exponential} function $\Exp$ associated to the (completed) representation ring of a rank-(1+1) torus $GL(1)_x\times(\mathbb{C}^*)_{q^{\frac{1}{2}}}$, where we identify the classes of the canonical 1-dimensional representations associated to each factor by $x,q^{\frac{1}{2}}$ respectively. The map $\Exp$ sends a class $[V]$ to the class $[S^\bullet V]$ of its symmetric algebra; for example, we have $\Exp(x) = 1/(1-x)$. In particular, note that $\Exp$ is multiplicative: $\Exp(f+g) = \Exp(f)\Exp(g)$. Now let $\mathcal{H} = \mathbb{C}[z]$ be a polynomial ring in a single variable $z$, where the monomial $z^k$ has weight $1$ with respect to $(\mathbb{C}^*)_x$ and weight $k+\frac{1}{2}$ with respect to $(\mathbb{C}^*)_{q}$. Then
$$
[\mathcal{H}] = \frac{q^{\frac{1}{2}}x}{1-q},
$$
and we have
\begin{align*}
\Phi(-x) &= \Exp(\mathcal{H})\\
&=\Exp\left(\frac{x}{q^{-1/2}-q^{1/2}}\right)
\end{align*}
\end{itemize}
The combinatorics of the cluster atlas for $\cX_{\PGL_2,S}$ is governed by the \emph{(cluster) modular groupoid} $\mathrm{M}_{\cX_{\PGL_2,S}}$ \cite[Section 6.1]{FG2}. The objects in this groupoid are \emph{isomorphism classes} of ideal triangulations (i.e. considered up to all homeomorphisms of $S$, not just those homotopic to the identity) on $S$, while its morphisms are generated by involutive ones corresponding to flips of diagonals (dually, graph mutations). The relations among the morphisms are generated by those of two kinds: the `rectangles', which assert that flips in disjoint quadrilaterals commute, and the `pentagons' which assert that the composite of five flips shown in Figure~\ref{fig:pentagon-rel} is the identity. The set of all cluster seeds corresponds to that of all morphisms in $\mathrm{M}_{\cX_{\PGL_2,S}}$. 

The space $\cX_{\PGL_2,S}$ is then obtained by gluing together a collection of Poisson tori, one for each cluster seed, where for each morphism $\mu_{e_{i_1}}\circ\cdots\circ\mu_{e_{i_k}}:\Gamma_1\rightarrow\Gamma_2$ the gluing map between the tori $\mathcal{T}_{\Gamma_1},\mathcal{T}_{\Gamma_2}$  is defined to be the composite of the birational transformations given by the $q^{\frac{1}{2}}=1$ specialization of Figure~\eqref{fig:q-mut}. That this prescription for the gluing data does not depend on the choice of representative of the equivalence class of a morphism follows from the \emph{pentagon identity} for the quantum dilogarithm:
\begin{align}
\label{eq:ab-penta}
(u,v)=1 \implies \Phi(X_{u})\Phi(X_{v}) = \Phi(X_{v})\Phi(X_{u+v})\Phi(X_{v}).
\end{align}

In the quantum setting, one can use a similar gluing  construction of to define a `quantum category of sheaves' on $\cX_{\PGL_2,S}$ as a category of \emph{descent data}, with the gluing pattern again determined by the groupoid $\mathrm{M}_{\cX_{\PGL_2,S}}$; see Remark 4.4 in ~\cite{SSZ} for further details. 

\subsection{Face Relations}
\label{subsec:face-relations}

In Section \ref{subsec:globalization} above, we defined the Fock-Goncharov quantum cluster Poisson variety
$\cX_{\PGL_2,S}$
of Borel-decorated $\PGL_2$ local systems on a punctured surface, $S.$  Cluster charts are labeled by cubic
graphs in $S$ with punctures in faces, or dually ideal triangulations.  
The subvariety of local systems with unipotent monodromy around punctures is a
quantized symplectic leaf, $\cX^{\mathrm{un}}_{\PGL_2,S},$ whose underlying (commuting) variety has
two components (see \cite[Proposition 2.10]{SSZ}).

We now fix $S = S^2\setminus \{p_1,...,p_{g+3}\}$ and write
simply $\cP$ for a particular \emph{component} of the unipotent symplectic leaf, $\cX^{\mathrm{un}}_{\PGL_2,S}$
which we will single out by equations below.
We write $\cP_\Gamma$ for the chart of $\cP$ obtained by a cubic planar graph, $\Gamma$.
As shown in \cite{SSZ}, there is a quantum Lagrangian subvariety, $\cM\subset \cP$,
and $\cM_\Gamma := \cM\cap \cP_\Gamma$
is cut out in the chart $\cP_\Gamma$ by explicit equations, one for each face of $\Gamma$.
This \emph{chromatic Lagrangian}
quantizes the Lagrangian subvariety of decorated local systems with trivial monodromy,
which by Item \eqref{item:chromatic} are nothing other than colorings of the faces of $\Gamma$ by elements of $\bP^1,$
modulo $\PGL_2$.
The ideals defined by these equations are compatible with the cluster structure:  the mutated face relations are
the face relations of the mutated graph.  This is the content of 
\cite[Theorem 4.3]{SSZ}, which we recall below.

Fix a cubic graph $\Gamma$ with $g+3$ faces and let $\mathcal{T}^q_{\Gamma}$ be the associated quantum torus
\begin{equation}
\label{eq:quantum-torus}
X_e X_{e'} = q^{\frac{1}{2}(e,e')}X_{e+e'} = q^{(e,e')} X_{e'}X_e.
\end{equation}
Let $f$ be a face of $\Gamma,$
with edges $e_1,...,e_n$ listed in counterclockwise order around the face,
with indices taken modulo $n$.
If $n>2$ then $e_{i+1}$ \emph{precedes} $e_i$ in counterclockwise order around their unique shared vertex,
and so the Poisson form is $(e_{i+1},e_i) = 1.$  If $n=2$ then the sum of the intersection forms at the two
half-edges is zero, and $(e_1,e_2) = 0.$
For $n>2,$ then, the corresponding cluster coordinates $X_{e_i}$ obey
$$X_{e_{i+1}}X_{e_i} = q^\frac{1}{2} X_{e_i + e_{i+1}} = q X_{e_i}X_{e_{i+1}},$$
with all other pairs commuting.

Then the relation which imposes unipotency of the
monodromy around the puncture corresponding to $f$ is
quantized as the multiplicative face relation
\begin{equation}
\label{eq:mult-face-relation}
X_{e_1}\cdots X_{e_n} = q^{-n/2},\quad\text{or equivalently}\;\quad X_{e_1+\dots +e_n} = q^{-1}
\end{equation}
while the distinguished component $\cP$
is cut out 
by the additional global relation 
\begin{equation}
\label{eq:global-relation}
    X_s = (-q)^{\frac{g+3}{2}},
\end{equation}
where $s = \sum_{e\in E}X_e$ and $g+3$ is the number of faces of $\Gamma.$
After quotienting $\mathcal{T}^q_{\Gamma}$ by these relations, we obtain a symplectic quantum torus algebra $\mathcal{T}^q_{\underline\Gamma}$. 

The \emph{additive} face relation $R_f$ (or $R_{\Gamma,f}$ if we want to show dependence on $\Gamma$)
quantizes the condition imposing triviality of the underlying unipotent
local system around the puncture inside $f$.
To this end,  set
\begin{align}
\label{eq:quantum.face.2h}
R_{f} &= {q^{-1/2}}+X_{e_1} + X_{e_1+e_2}+\ldots + X_{e_1+e_2+\cdots e_{n-1}}\\
\nonumber &= q^{-1/2} + X_{e_1} + q^{1/2}X_{e_1}X_{e_2} + qX_{e_1}X_{e_2}X_{e_3}+\cdots + q^{(n-2)/2}X_{e_1}X_{e_2}\cdots X_{e_{n-1}}.
\end{align}

\begin{remark}
\label{rmk:cyclic-order}
From the multiplicative face relation~\eqref{eq:mult-face-relation}, we see
that multiplying ~\eqref{eq:quantum.face.2h} by $q^{1/2}X_{e_n}$ yields
$$X_{e_n}+X_{e_n+e_1} + X_{e_n+e_1+e_2}+\ldots + q^{-1/2},$$
so we see that the ideal in the quantum torus $\mathcal{T}^q_{\underline{\Gamma}}$ generated by $R_f$ is independent of our arbitrary linearization of the cyclic order on the edges around the face $f$ implicit in~\eqref{eq:quantum.face.2h}.
\end{remark}

Let $\mathcal{I}_{\Gamma}$ be the left ideal in $\mathcal{T}_{\Gamma}^q$ generated by all~\eqref{eq:mult-face-relation}
along with the global relation~\eqref{eq:global-relation} and the relations $R_f$ for all faces $f$. As the quantization of a Lagrangian subvariety, the D-module $\mathcal{V}_{\Gamma}:=\mathcal{D}_{\Gamma}/ \mathcal{I}_{\Gamma}$ is holonomic.

Now suppose that two regular cubic graphs $\Gamma$ and $\Gamma'$ are related by mutation at  edge $e_0$. Let us write $\mathcal{T}_{\Gamma,\Gamma'}$ for the localization of the quantum torus $\mathcal{T}_{\Gamma}$ at the Ore set $\{(1+q^{k+1/2}X_{e_0})\}_{k\in\mathbb{Z}}$, and similarly write $\mathcal{T}_{\Gamma',\Gamma}$ for the localization of $\mathcal{T}_{\Gamma}$ at $\{(1+q^{k+1/2}X'_{e_0})\}_{k\in\mathbb{Z}}$. Then the quantum mutation map $\mu_0$ in~\eqref{fig:q-mut} defines an isomorphism
{$ \mu_0~:~\mathcal{T}_{\Gamma',\Gamma}\rightarrow \mathcal{T}_{\Gamma,\Gamma'}$. Let us write $\mathcal{I}_{\Gamma,\Gamma'}$ for the ideal in $\mathcal{T}_{\Gamma,\Gamma'}$ generated by the quantized chromatic ideal $\mathcal{I}_\Gamma$, and $\mathcal{I}_{\Gamma',\Gamma}$ for the ideal in $\mathcal{T}_{\Gamma',\Gamma}$ generated by $\mathcal{I}_{\Gamma'}$.}

\begin{theorem}[\cite{SSZ}]
\label{thm:d-mod} The system of quantized chromatic ideals $\{\mathcal{I}_\Gamma\}$ is compatible with quantum cluster mutations: if $\Gamma,\Gamma'$ are regular cubic graphs related by a flip at  edge $e_0$ as in Figure~\ref{fig:q-mut}, then we have $\mu_0(\mathcal{I}_{\Gamma',\Gamma}) = \mathcal{I}_{\Gamma,\Gamma'}$.
\end{theorem}

\begin{proof}
Consider the generator $R_{\Gamma',f}$ of $\mathcal{I}_{\Gamma',\Gamma}$ associated to the left face of the graph $\Gamma'$ in Figure~\ref{fig:q-mut}, as defined in ~\eqref{eq:quantum.face.2h}. 
(Note the $e_i$ do not correspond to consecutive edges around a face here.)
We show that it is mapped to the corresponding 
to a generator $R_{\Gamma,f}$ of $\mathcal{I}_{\Gamma,\Gamma'}$ under $\mu_0$. As explained in Remark~\ref{rmk:cyclic-order},
by multiplying $R_{\Gamma',f}$ by a unit in $\mathcal{T}^q_{\underline{\Gamma'}}$ we may assume that the edge $e_0$ at which we mutate is neither $e_1$ nor $e_{n-1}$ in the notations of~\eqref{eq:quantum.face.2h}, since we do
not consider mutations along edges of bigons.
Then reading counterclockwise around the left face of the right graph in Figure~\ref{fig:q-mut}, we see that
\begin{align*}
\mu_0(X'_{e_2} +q^{1/2}X'_{e_2}X'_{e_0} + q  X'_{e_2}X'_{e_0}X'_{e_1}) &= X_{e_2}(1+q^{1/2} X_{-e_0})^{-1} + q^{1/2} X_{e_2}X_{-e_0}(1+q^{1/2} X_{-e_0})^{-1} \\
&\phantom{=}\; +  qX_{e_2}X_{-e_0}(1+q^{1/2} X_{-e_0})^{-1}X_{e_1}(1+q^{1/2} X_{e_0})\\
&= X_{e_2}+ q^{1/2} X_{e_2}X_{e_1},
\end{align*}
where we used that $X_{e_0}X_{e_1} = q^2X_{e_1}X_{e_0}$ by the  relation~\eqref{eq:quantum-torus} applied to the graph on the left of Figure~\ref{fig:q-mut}. From this computation, we see that $\mu_0(R_{\Gamma',f}) = R_{\Gamma,f}$. The intertwining of the generators of the form~\eqref{eq:mult-face-relation} and~\eqref{eq:global-relation} follows in exactly the same way.
\end{proof}

\section{Skein-theoretic framed seeds}
\label{sec:groupoid}
In this section we propose a skein-theoretic extension of the notion of a framed seed  from \cite{SSZ}.   In the present context, the framing data will furnish an isomorphism of a given Lagrangian filling $L$ with the standard handlebody; such an isomorphism determines a basis for $\mathrm{Sk}(L)$ into which wavefunctions can be expanded, and numerical invariants extracted. 

\subsection{Fillings and framings}
Let us briefly recall some of the combinatorial and geometric setup employed in \cite{TZ,SSZ} to describe Lagrangian fillings of the Legendrian surfaces $S_\Gamma$. 
Such a Lagrangian is encoded by a \emph{smoothed ideal foam} filling the graph $\Gamma$, and admits a topological description as the double cover of the 3-ball branched over a tangle $\mathcal{T}$ whose endpoints lie the vertices of $\Gamma$. In order to ensure the resulting filling $L$ has the topology of a handlebody, we assume the tangle $\mathcal{T}$ to be \emph{rational} in the sense of Conway: that is, there exists a homeomorphism of $B^3$ with  $\mathbb{D}^2\times [-1,1]$ and a partition of the vertices of $\Gamma$ into two disjoint sets $V_\pm$, such that $\mathcal{T}$ is identified with the trivial tangle connecting vertices $V_-=\{v_1,\ldots, v_{g+1}\}\times \{-1\} $
with $V_+=\{v_1,\ldots, v_{g+1}\}\times \{1\}$.


Now recall that in the geometric setting of \cite[Section 5 and Definition 5.6]{SSZ},
a framing is an isotropic splitting of the short exact sequence $K\to H_1(S_\Gamma)\to H_1(L)$, and different choices of splitting of the same exact sequence are said to differ by a `change of framing'.

\begin{example}
When $S_\Gamma$ has genus $g$,
changes of framings as defined in \cite{SSZ} are parametrized by $g\times g$ symmetric,
integer matrices.
In the case $g=1$ familiar from knot theory, we have a single integer $p$ which can be regarded as the winding around the meridian
of a path transverse to the longitude of the torus.
\end{example}

\begin{figure}[htpb]
    \centering
    \begin{tikzpicture}
    \pgfmathsetmacro{\g}{4}
    \pgfmathsetmacro{\gg}{\g+1}
        \draw[blue,very thick] (0,-1)--(\gg,-1)--(\gg,1)--(0,1)--(0,-1);
        \foreach \i in {1,...,\g}{\draw[blue,very thick] (\i,0) circle (1/4);}
    \end{tikzpicture}
    \qquad
    \begin{tikzpicture}[yscale=.4]
    \pgfmathsetmacro{\g}{4}
    \pgfmathsetmacro{\gg}{\g+1}
    \pgfmathsetmacro{\vert}{3}
        \draw[blue,very thick] (-1/2,-1)--(-1/2+\gg,-1)--(1/2+\gg,1);
        \draw[blue,dashed] (1/2+1/4,1)--(1/2+0,1)--(-1/2+0,-1);
        \draw[blue,dashed] (\g+1/4,1)--(\g+3/2,1);
        \draw[blue,very thick] (-1/2,-1+\vert)--(-1/2+\gg,-1+\vert)--(1/2+\gg,1+\vert);
        \draw[blue,very thick] (1/2+\gg,1+\vert)--(1/2+0,1+\vert)--(-1/2+0,-1+\vert);
        \draw[blue,very thick] (-1/2,-1)--(-1/2,-1+\vert);
        \draw[blue,very thick] (-1/2+\gg,-1)--(-1/2+\gg,-1+\vert);
        \draw[blue,very thick] (1/2+\gg,1)--(1/2+\gg,1+\vert);
        \draw[blue,dashed] (1/2,1)--(1/2,1+\vert);

        \foreach \i in {1,...,\g} {\draw[blue, thick, dashed](\i+1/4,0)--(\i+1/4,\vert);
        \draw[blue, thick, dashed](\i-1/4,0)--(\i-1/4,\vert);
        \draw[red, thick, dashed](\i,0)--(\i,\vert-2/10);
        \draw[red,thick](\i-1/3,\vert-1)--(\i-1/3,-1);
        \draw[red, thick] (\i-1/3,\vert-1)--(\i-1/3+2/5,\vert-1+4/5);
        \draw[red,thick,dashed] (\i-1/3,-1)--(\i-1/3+2/5,-1+4/5);
        \draw[blue,dashed] (\i+1/4,1)--(\i+3/4,1);
        \draw[blue,thick,dashed] (\i,\vert*0) circle (1/4);\draw[blue,very thick] (\i,\vert*1) circle (1/4);
        }
    \end{tikzpicture}
    \caption{$D_g$ (left) and the handlebody $\cH_g$ (right), for $g=4$}
    \label{fig:handlebody}
\end{figure}
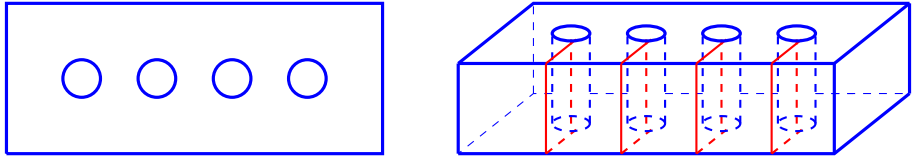

To generalize to skeins, we consider the following nonlinear counterpart of this definition. Write $\cH_g$ for the standard genus-$g$ handlebody,
which we parameterize explicitly as follows.
Let $D = [0,g+1]\times [-1,1]\subset \bR^2$, write $B_{\frac{1}{4}}(k,0)$ for disks of radius $\frac{1}{4}$ centered at $(k,0),$
and let $D_g := D\setminus \bigsqcup_{k=1}^g B_{\frac{1}{4}}(k,0)$ be the
disk minus these $g$ disjoint smaller disks --- \blue{see Figure \ref{fig:handlebody}}.
Then set $\cH_g = D_g \times I,$
an oriented three-manifold, and put $\Sigma_g = \partial \cH_g.$
We equip the fundamental group $\pi_1(\Sigma_g)$ with a distinguished set of generators $(a_1,b_1,\ldots, a_g,b_g)$ as indicated at right in Figure~\ref{fig:handlebody}: each $a_i$ is represented by a loop going {clockwise} around the boundary of the $i$th small disk, and each each $b_i$ (shown in red in the figure) is an essential curve on $\Sigma_g$ bounding a disk in $\cH_g $  (i.e.~a \emph{meridian}). 
In particular, our parametrization fixes an isomorphism
$$
\pi_1(L) \simeq \langle a_1,\ldots, a_g\rangle \subset \pi_1(\Sigma_g).
$$

Now recall the definition of the handlebody group,
or mapping class group $\mathrm{MCG}(\cH_g) := \pi_0\mathrm{Homeo}^+(\cH_g,\cH_g)$,
i.e.~isotopy classes of orientation-preserving
self-homeomorphisms of $\cH_g$.  Equivalently,
it is the subgroup of the mapping class group of the boundary $\mathrm{MCG}(\Sigma_g)$
consisting of homeomorphisms of which
extend to the interior. 
Closer to the spirit of Definition 5.6 of~\cite{SSZ}, the subgroup $\mathrm{MCG}(\cH_g)\subset \mathrm{MCG}(\Sigma_g)$ admits an alternative characterization as those mapping classes preserving
the kernel of the map $\pi_1(\Sigma_g)\to \pi_1(\cH_g)$
--- see \cite{H}.

\begin{definition}
Let $L$ be a filling of $S_\Gamma$.
    A \emph{framing} is a connected component of $\mathrm{Homeo}^+(L, \cH_g)$.
\end{definition}

\begin{remark}
The set of framings $\pi_0 \mathrm{Homeo}^+(L, \cH_g)$ is a torsor over the handlebody group
$\mathrm{MCG}(\cH_g) := \pi_0\mathrm{Homeo}^+(\cH_g,\cH_g).$  We call an element of the handlebody
group a \emph{change of framing}.
\end{remark}

\begin{remark}
    Taking the action on homology recovers the previous definitions from \cite{SSZ}
    of geometric framings and changes of framing.  In particular, the retraction $\cH_g \sim S_\Gamma$
    gives the isotropic splitting $H_1(L)\to H_1(S_\Gamma)$.
\end{remark}

\begin{example}
    When $g=1,$ nonlinear framings are still a torsor over the additive group $\bZ.$
\end{example}





\begin{definition} A \emph{geometric framed seed} of genus $g$ consists of the following data:
\begin{enumerate}[(i)]
\item A genus-$g$ cubic planar graph $\Gamma$ on $S^2=\partial B^3$;
\item A rational $(g+1)$-tangle $\mathcal{T}$ in $B^3$ with endpoints at the $2g+2$ vertices of $\Gamma$;
\item A framing $\mathbf{f}\in \pi_0\mathrm{Homeo}^+(L_{\mathcal{T}},\mathcal{H}_g)$ of the 3-manifold $L_{\mathcal{T}}$. 
\item For each face $f\in F_\Gamma$, a distinguished point $\fp^f$, and choice of \emph{framing path} $\gamma^f$ from $\fp^f$ to some vertex $v\in V_\Gamma$
that intersects no other vertices.
\end{enumerate}
\end{definition}
As usual, we consider the data above modulo the action of the identity component of the group of orientation preserving homeomorphisms of the ball which fix each vertex of $\Gamma$. 

We say that two geometric framed seeds $(\Gamma_1,\mathcal{T}_1,\mathbf{f}_1,\{\gamma_1^f\}), (\Gamma_2,\mathcal{T}_2,\mathbf{f}_2,\{\gamma_2^f\})$ are related by a change of framing if $\Gamma_1=\Gamma_2,\mathcal{T}_1=\mathcal{T}_2$. In this case the fillings $L_1\simeq L_2$ are \emph{canonically} homeomorphic (as the $(\Gamma_i,\mathcal{T}_i)$ are isotopic), and $\mathbf{f}_2\circ\mathbf{f}_1^{-1}$ is an element of $\mathrm{MCG}(\mathcal{H}_g)$.

\subsection{Mutation of framed seeds}
Suppose that $(\Gamma,\mathcal{T},\mathbf{f},\{\gamma^f\})$ is a geometric framed seed of genus $g$, and that $e$ is an edge of $\Gamma$. The edge $e$ corresponds the diagonal of a quadrilateral in the triangulation of $S^2$ dual to $\Gamma$. 

\begin{figure}[ht]
    \centering
    \begin{tikzpicture}[scale=.5]
    \pgfmathsetmacro{\a}{3}
    \pgfmathsetmacro{\b}{2.5}
        \draw[thick,black] (-\a,0)--(0,-\b)--(\a,0)--(0,\b)--(-\a,0);
        \draw[thick,black] (0,-\b)--(0,\b);
        \draw[thick,black,dashed] (-\a,0)--(\a,0);
        \filldraw[red] (\a/4,-\b/4) circle (.2cm) node[below]{$t^*_1$};
        \filldraw[red] (-\a/4,\b/4) circle (.2cm) node[above]{$t^*_2$};
        \filldraw[blue] (\a/4,\b/4) circle (.2cm) node[above]{$t_1$};
        \filldraw[blue] (-\a/4,-\b/4) circle (.2cm) node[below]{$t_2$};
    \end{tikzpicture}
    \qquad
        \begin{tikzpicture}[scale=.5]
    \pgfmathsetmacro{\a}{3}
    \pgfmathsetmacro{\b}{2.5}
        \draw[thick,black] (-\a,0)--(0,-\b)--(\a,0)--(0,\b)--(-\a,0);
        \draw[thick,black] (0,-\b)--(0,\b);
        \draw[thick,black,dashed] (-\a,0)--(\a,0);
        \draw[very thick,brown,-Stealth] (\a/4,-\b/4) .. controls (\a/3,0) .. (\a/4,\b/4);
        \draw[very thick,brown,-Stealth] (-\a/4,\b/4) .. controls (-\a/3,0) .. (-\a/4,-\b/4);

        \filldraw[red] (\a/4,-\b/4) circle (.2cm);
        \filldraw[red] (-\a/4,\b/4) circle (.2cm);
        \filldraw[blue] (\a/4,\b/4) circle (.2cm);
        \filldraw[blue] (-\a/4,-\b/4) circle (.2cm);
    \end{tikzpicture}
    \qquad
    \begin{tikzpicture}[scale=.5]
    \pgfmathsetmacro{\a}{3}
    \pgfmathsetmacro{\b}{2.5}
        \draw[thick,black] (-\a,0)--(0,-\b)--(\a,0)--(0,\b)--(-\a,0);
        \draw[thick,black] (0,-\b)--(0,\b);
        \draw[thick,black,dashed] (-\a,0)--(\a,0);
        \draw[very thick,brown,-Stealth] (\a/4,-\b/4) .. controls (0,-\b/3) .. (-\a/4,-\b/4);
        \draw[very thick,brown,-Stealth] (-\a/4,\b/4) .. controls (0,\b/3) .. (\a/4,\b/4);
        \filldraw[red] (\a/4,-\b/4) circle (.2cm);
        \filldraw[red] (-\a/4,\b/4) circle (.2cm);
        \filldraw[blue] (\a/4,\b/4) circle (.2cm);
        \filldraw[blue] (-\a/4,-\b/4) circle (.2cm);
    \end{tikzpicture}

    \caption{Tetrahedron (left) and tangles of positive (middle) and negative (right) mutations.}
    \label{fig:tetra}
\end{figure}
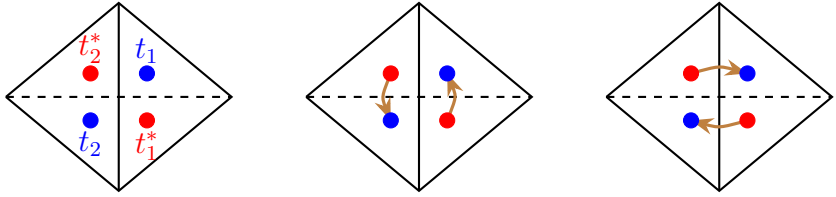

Consider a tetrahedron $T$ with a distinguished pair of faces $t^*_1,t^*_2$. We denote the other two faces by $t_1,t_2$ as shown in Figure~\ref{fig:tetra}. Recall from Section 5.7 of \cite{SSZ} the pair of smoothed Harvey-Lawson Lagrangians $\mathrm{HL}_{\pm}$ associated to such a tetrahedron: $\mathrm{HL}_{+}$, which corresponds to the tangle matching $t_1^*\leftrightarrow t_1$ and $t_2^*\leftrightarrow t_2$; and $\mathrm{HL}_{-}$ corresponding to the opposite matching.

For each edge $e$ of $\Gamma$ and choice of sign $\pm$, we define a 3-manifold $M_{\Gamma, e,\pm}$ as follows. Form  $S^2\times I$, and write $\Gamma_0,\Gamma_1$ for the graphs on $S^2\times\{0\}$ and  $S^2\times\{1\}$ respectively. For each vertex $v$ of $\Gamma$, we regard $v\times I$ as a strand of a trivial tangle in $S^2\times I$, and we consider 3-manifold obtained as the branched double cover of $S^2\times I$ along this tangle. Its boundary is identified with $S_{\Gamma_0}\sqcup S_{\Gamma_1}$. The 3-manifold $M_{\Gamma,e,\pm}$ is then obtained by attaching a tetrahedron $T_\pm$ to the quadrilateral on $S^2\times\{1\}$ containing $e$, thereby gluing a smoothed Harvey-Lawson Lagrangian $\mathrm{HL}_{\pm}$ to $S_{\Gamma_1}$ as shown in Figure~\ref{fig:tetra}. 

The 3-manifold $M_{\Gamma,e,\pm}$ is a (nonexact) \emph{Lagrangian concordance} between the Legendrians $S_{\Gamma}$ and $S_{\Gamma'}$, where $\Gamma'$ is the cubic graph obtained from $\Gamma$ by mutation at edge $e$. 

Now suppose that $L_\mathcal{T}$ is a Lagrangian filling of $S_\Gamma$ obtained as the branched double cover of $B^3$ along a tangle obtained from some smoothed ideal foam. Then gluing over $S_\Gamma$ produces a Lagrangian filling $\mu_{e,\pm}\left(L_{\mathcal{T}}\right) = L_{\mathcal{T}}\cup_{S_\Gamma}M_{\Gamma,e,\pm}$ of $S_{\Gamma'}$, which is the branched double cover over the tangle $\mu_{e,\pm}(\mathcal{T})$ obtained by propagating $\mathcal{T}$ to a new rational tangle with endpoints at vertices of $\Gamma'$ in accordance with the sign choice $\pm$.

As explained in Section 5.7 of~\cite{SSZ} (see, in particular, the proof of Proposition 5.22 of {\it loc.~cit.}), the concordance $M_{\Gamma,e,\pm}$ determines a homeomorphism $\varphi_\pm: L_{\mathcal{T}'} \rightarrow L_{\mathcal{T}}$.
Hence a framing $\mathbf{f}:{L}_{\mathcal{T}}\rightarrow \mathcal{H}_g$ gives rise to a framing $\mu_{e,\pm}(\mathbf{f}) = \mathbf{f}\circ\varphi_{\pm}$ for the mutated filling $\mathcal{L}_{\mathcal{T}'}$. Similarly, we transport the framing paths $\{\gamma^f\}$ for $L_{\mathcal{T}}$ to ones for $L_{\mathcal{T}'}$ by means of the isomorphism $(\mathbf{f}')^{-1}\circ\mathbf{f}$.

When $(\Gamma',\mathcal{T}',\mathbf{f}',\{\gamma^{f'}\})=\mu_{e,\pm}(\Gamma,\mathcal{T},\mathbf{f},\{\gamma^f\})$, we say that the two framed seeds are related by \emph{mutation at edge $e$ with sign }$\pm$. 

\subsection{{Skein Cluster Groupoid}}
\label{sec:skein-cluster-groupoid}

Consider the groupoid $\mathbb{G}$ whose objects are geometric framed seeds, and whose morphisms are formal composites of mutations and changes of framing.  

Just as in~\cite{SSZ}, we will also need a smaller groupoid $\mathbb{G}_{\mathrm{ad}}$ in which we allow as morphisms not all signed mutations, but only the \emph{admissible} ones.  In the present context, we define a mutation $\mu_{e,\pm}:(\Gamma',\mathcal{T}',\mathbf{f}',\{\gamma^{f'}\})\rightarrow (\Gamma,\mathcal{T},\mathbf{f},\{\gamma^f\})$ to be admissible if $\delta(\pm\mathbf{f}_*[e])>0$, where $[e]\in H_1(S_\Gamma,\mathbb{Z})$ is the homology class of the loop on $S_\Gamma$ corresponding to edge $e$, and $\delta:H_1(\Sigma_g,\mathbb{Z})\rightarrow \mathbb{Z}$ is the standard functional defined in Example~\ref{eg:completion} given by $\delta([a_i])=1, \delta([b_i])=0$. In particular, the notion of admissibility is independent of the choice of framing and framing path in a given framed seed.

\subsection{Algebraic data associated to framed seeds and their mutations}
To a geometric framed seed $(\Gamma,\mathcal{T},\mathbf{f},\{\gamma^f\})$, we can associate the following algebraic data:
\begin{enumerate}
\item A module $\Sk(L_{\mathcal{T}})$ over the algebra $\Sk(S_\Gamma)$.

\item An algebra isomorphism $\Sk(S_\Gamma)\rightarrow\Sk(\Sigma_g)$, and a compatible module isomorphism $~\Sk(\mathcal{L}_{\mathcal{T}})\rightarrow\Sk(\mathcal{H}_g).$
\item An algebra isomorphism $\Sk(S_\Gamma,\fp^\pm)\rightarrow \Sk(S_\Gamma)$, determined by the framing paths.
\end{enumerate}
We can be particularly concrete about datum (2) in the case $g=1$: there it consists of an isomorphism of the algebra $\Sk(S_\Gamma)$ with the abstract $R$-algebra~\eqref{eq:eha-rels}, together with an isomorphism of $\Sk(L_\mathcal{T})$ with the ring of symmetric functions, compatible with the action defined by~\eqref{eq:MN-rule},~\eqref{eq:P-eigenvalues}. 

When $g>1$, it is argued in~\cite{P1} that the skein $\mathrm{Sk}(\mathcal{H}_g)$ coincides as an $R$-module with the symmetric algebra on the span of all nontrivial conjugacy classes in $\pi_1(\mathcal{H}_g)$. Thus, as per Figure 1.2 in ~\cite{P1}, a geometric framed seed again determines a basis in  $\mathrm{Sk}(\mathcal{H}_g)$. 

{ 
Again following Section 3.4 of \cite{SSZ}, to any formal composite of framing shifts and admissible mutations we can associate an element of the semidirect product $R$-algebra $\widehat{\Sk} (\Sigma_g)\rtimes R\mathrm{MCG}(\cH_g)$ of the $\delta$-completed skein algebra of $\Sigma_g$ with the group algebra of the handlebody group. To a signed mutation at edge $e$ we assign the \emph{Baxter operator} $Q_{\pm \mathbf{f}(e)}^{\pm1}\in \widehat{\Sk} (\Sigma_g)$ {(see Definition~\ref{def: Baxter operator})}, while the framing shifts are simply mapped to the corresponding elements of $\mathrm{MCG}(\cH_g)$. 

In particular, using the natural action of $\widehat{\Sk} (\Sigma_g)\rtimes R\mathrm{MCG}(\cH_g)$ on the completed handlebody skein module $\widehat{\Sk}(\cH_g)$, any composite of admissible mutations and changes of framing data gives rise to an $R$-automorphism of $\widehat{\Sk}(\cH_g)$.

As in~\cite{SSZ}, we may impose a relation on morphisms in the admissible groupoid $\mathbb{G}_{\mathrm{ad}}$ by identifying any two morphisms with the same source and target for which the corresponding elements of 
$\widehat{\Sk} (\Sigma_g)\rtimes R\mathrm{MCG}(\cH_g)$ are \emph{equal}.

It would be desirable to describe this equivalence relation in terms of the Lagrangian concordances defining mutations of fillings. We save the problem of developing such a description for a future occasion. 
}

\section{Skein Face Relations and Mutations}
{We are now ready to define skein mutation, and thereby upgrade Theorem~\ref{thm:d-mod} to the skein setting.}

\label{sec:face-mutations}

\subsection{Skein valued face relation}
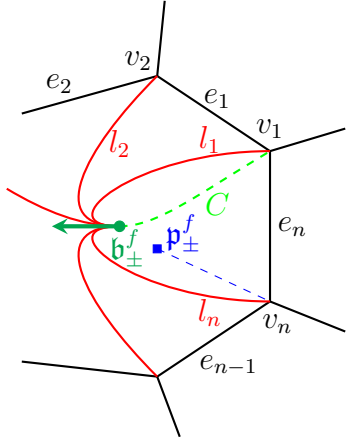
\begin{figure}[htpb]
	\centering
	\begin{tikzpicture}
		\draw[thick] (-3.3, 1.5) -- (-1.5, 2)--(0,1)--(0, -1) -- (-1.5, -2) -- (-3.3, -1.8);
		\draw[thick] (-1.5, 2)--(-1.4, 3) ;
		\draw[thick] (0, 1) -- (1, 1.3);
		\draw[thick] (0, -1) -- (1, -1.4);
		\draw[thick] (-1.5, -2) -- (-1.2, -2.8);
		
		\node at (0, 1.3) {$v_1$};  
		\node at (0.1, -1.3) {$v_n$};
		\node at (-1.75, 2.2) {$v_2$};
		\node at (-0.7, 1.7) {$e_1$};
		\node at (-0.55, -1.85) {$e_{n-1}$};
		\node at (0.3, 0) {$e_n$}; 
		\node at (-2.8, 1.9) {$e_2$};
		
		\draw[red, thick] (-2, 0) .. controls (-3, 0) and (-2.5, 1) .. (-1.5, 2);
		\node at (-2, 1.1) {\color{red} $l_2$};
		
		\draw[red, thick] (-2, 0) .. controls (-3, 0) and (-2.5, -1).. (-1.5, -2);
		\node at (-0.8, 1.15) {\color{red} $l_1$};
		
		\draw[red, thick] (-2, 0) .. controls (-3, 0) and (-1.8, 1).. (0, 1);

		\draw[red, thick] (-2, 0) .. controls (-3, 0) and (-1.8, -1).. (0, -1);
		\node at (-0.8, -1.15) {\color{red} $l_n$};
		
		\draw[red, thick] (-2, 0) .. controls (-2.5, 0) and (-3, 0.2) .. (-3.5, 0.5);

		\draw[thick, green, dashed] (-2, 0) .. controls (-1.5, 0) and (-1, 0.4) .. (0, 1);
		\node at (-0.7, 0.3) {\color{green}$C$};

		\draw[blue, dashed] (-1.5, -0.3) -- (0, -1) ;
		
		\filldraw[blue] (-1.5, -0.3) [yshift = -1.5pt, xshift = -1.5pt] rectangle ++(3pt, 3pt) ;
            \node at (-1.15, -0.1) {\color{blue} $\fp_\pm^f$};
  
		\filldraw[Green] (-2, 0) circle(2pt);
            \node at (-1.9, -0.3) {\color{Green} $\fb_\pm^f$};
		
		\draw[Green, ultra thick, ->, > = stealth] (-2, 0) -- (-2.9, 0); 
		
	\end{tikzpicture}
	\caption{A face $f$ of $\Gamma$ with base points $\fb^f_\pm$, framing points $\fp^f_\pm$, capping path $C$
 and skeins $l_i$.}
	\label{fig: face relation}
\end{figure}

Let $\Gamma$ be a cubic graph on $S^2$, $S_\Gamma$ its associated Legendrian surface, and $L$ a Lagrangian filling. For each face $f$, there is Reeb chord; we denote its endpoints by $\fb^f_+$ (on the upper sheet) and $\fb^f_-$ (on the lower sheet), and take the $0$-chain $\fb^f_\pm \coloneqq  \fb^f_+ - \fb^f_-$.  As is shown in Figure \eqref{fig: face relation}, we first choose a distinguished edge (here it is $e_n$), and the framing of $\fb_\pm$ to be pointing away from the edge. Then we choose another pair of points $\fp^f_\pm$, one on each sheet, lying in the triangle formed by the vertices of the distinguished edge and $\fb^f_\pm$. Globally, we take $\fb_\pm \coloneqq  \sum_{f \in F} \fb_\pm^f$, and $\fp_\pm = \sum_{f \in F} \fp_\pm^f $. 

Recall from section \ref{sec:skeins with b and f} we have the skein module $\Sk( S_\Gamma, \fb_\pm; \fp_\pm)$. For a loop or path $l$, we will denote the element it represents in the skein module with framing line by $\frame{l}$, and the element it represents in the normal skein module by $\skein{l}$. For an edge $e$ of the graph $\Gamma$, we write $\loo{e}$ for the loop it determines in $S_\Gamma$, and take $\skein{e} \coloneqq  \skein{\loo{e}}$.

We define the skein-valued face relation to be 
\begin{equation}
    \label{eq:skein-face-relation}
		A_{\Gamma, f} = \frame{l_1} + \frame{l_2} + \cdots + \frame{l_n} \in \Sk (S_\Gamma, \fb_\pm; \fp_\pm)
    \end{equation}
where $l_i$ is the path connecting $\fb^f_-$ to $\fb^f_+$ passing through
the vertex $v_i$, as is shown in figure \eqref{fig: face relation}. (The real $l_i$ as actually a double cover of the red line, branching at $v_i$.) The $\fp_\pm^f$ lie in the region adjacent to the edge $e_n$. 

Later on we will see that, as in \cite{SS}, this face relation annihilates the 
skein-valued wavefunction of $L$ defined by \cite{ES19}, under the map \eqref{action}. 

\begin{example}\label{eg: equivalent to SS face relation}
    If we choose a capping path $C$ going from $\fb_+^f$ to $\fb_-^f$ as in figure \eqref{fig: face relation}, then after capping with $C$ the path $l_1$ becomes an unknot with a twist, and $l_i$ becomes a loop isotopic to $\loo{e}_1 + \cdots +\loo{e}_{i-1}$. Here ``$+$'' means connected sum, about which we will say more in Section \ref{subsec:mutation}. Thus under the map \eqref{base points}, the face relation is sent to 
    \begin{equation*}
        a^{-1} \bigcirc + \frame{\loo{e}_1} + \frame{\loo{e}_1 + \loo{e}_2} + \cdots + \frame{\loo{e}_1 + \cdots + \loo{e}_{n-1}}, 
    \end{equation*}
    which is the face relation in \cite{SS}, and which recovers the face relation \eqref{eq:quantum.face.2h}
    upon mapping to the linking skein \eqref{map: link alg}.
\end{example}

Summarizing, if we choose a capping path, the element $A_{\Gamma,f}$ is mapped to the face relation defined in \cite{SS}.

For different choices of the framing of $\fb_\pm^f$, or equivalently, choices of the distinguished edge (remember that these choices are
one and the same by the standard form of loops in Figure \ref{fig: face relation}), there is an isomorphism between the corresponding skein
modules, given by ``rotating continuously''.  The next lemma establishes independence of the face relations on this choice.

\begin{lemma}
\label{lem: framing-independence}
    The element $A_{\Gamma,f}$ does not depend on the framing of $\fb_\pm^f$, in the sense that the rotation isomorphism identifies the two face relations. 
\end{lemma}

\begin{proof}

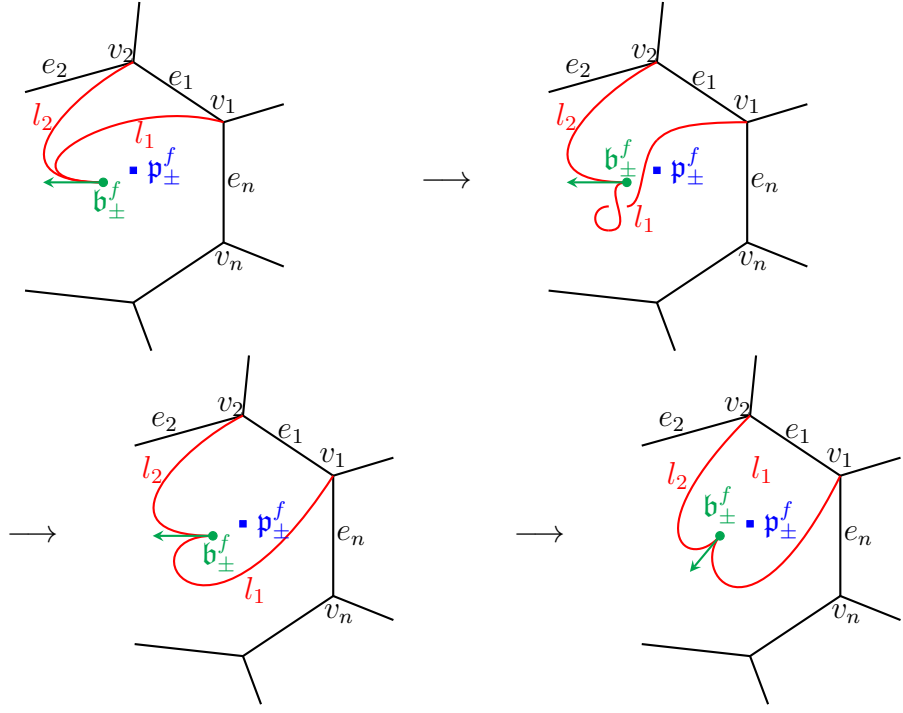
\begin{figure}[htbp]
	\centering
	\begin{subfigure}{0.25\textwidth} 
		\centering

	\begin{tikzpicture}[scale = 0.8]
		\draw[thick] (-3.3, 1.5) -- (-1.5, 2)--(0,1)--(0, -1) -- (-1.5, -2) -- (-3.3, -1.8);
		\draw[thick] (-1.5, 2)--(-1.4, 3) ;
		\draw[thick] (0, 1) -- (1, 1.3);
		\draw[thick] (0, -1) -- (1, -1.4);
		\draw[thick] (-1.5, -2) -- (-1.2, -2.8);
		
		 \node at (0, 1.3) {$v_1$};  
		\node at (0.1, -1.3) {$v_n$};
		\node at (-1.7, 2.15) {$v_2$};
		\node at (-0.7, 1.7) {$e_1$};
		 
		\node at (0.3, 0) {$e_n$}; 
		\node at (-2.8, 1.9) {$e_2$};
				
		\draw[red, thick] (-2, 0) .. controls (-4, 0) and (-2.5, 1.5) .. (-1.5, 2);
		\node at (-3, 1.1) {\color{red} $l_2$};

		\draw[red, thick] (-2, 0) .. controls (-4, 0) and (-1.8, 1.5).. (0, 1);
		\node at (-1.3, .8) {\color{red} $l_1$};

		\filldraw[blue] (-1.5, 0.2) [yshift = -1.5pt, xshift = -1.5pt] rectangle ++(3pt, 3pt) node[anchor = west] {\color{blue} $\fp_\pm^f$}; ;

		\filldraw[Green] (-2, 0) circle(2pt);
		\node at (-1.9, -0.3) {\color{Green} $\fb_\pm^f$};
				
		 \draw[Green, thick, ->, > = stealth] (-2, 0) -- (-3, 0);

	\end{tikzpicture}
	\end{subfigure} 
	\begin{subfigure}{0.5\textwidth} 
		\centering
	\begin{tikzpicture}[scale = 0.8]
		\node at (-5, 0) {$\longrightarrow$}; 
		
		\draw[thick] (-3.3, 1.5) -- (-1.5, 2)--(0,1)--(0, -1) -- (-1.5, -2) -- (-3.3, -1.8);
		\draw[thick] (-1.5, 2)--(-1.4, 3) ;
		\draw[thick] (0, 1) -- (1, 1.3);
		\draw[thick] (0, -1) -- (1, -1.4);
		\draw[thick] (-1.5, -2) -- (-1.2, -2.8);
		
		 \node at (0, 1.3) {$v_1$};  
		\node at (0.1, -1.3) {$v_n$};
		\node at (-1.7, 2.15) {$v_2$};
		\node at (-0.7, 1.7) {$e_1$};
	
		\node at (0.3, 0) {$e_n$}; 
		\node at (-2.8, 1.9) {$e_2$};
		
		\draw[red, thick] (-2, 0) .. controls (-4, 0) and (-2.5, 1.5) .. (-1.5, 2);
		\node at (-3, 1.1) {\color{red} $l_2$};

		\draw[red, thick] (-2, 0) .. controls (-2.5, 0) and (-1.8, -.8) .. (-2.3, -0.8);
		\draw[red, thick] (-2.3, -0.8) .. controls (-2.6, - 0.8) and (-2.6, -0.4) .. (-2.3, -0.4) ;
		\draw[red, thick] (-2, -0.4) .. controls (-1.7, -0.4) and (-1.8, 0.4).. (-1.5, 0.7) ;
		\draw[red, thick] (-1.5, 0.7) .. controls (-1.3, 0.9) and (-1, 1) .. (0, 1);
		
		\node at (-1.7, -0.6) {\color{red} $l_1$};

		\filldraw[blue] (-1.5, 0.2) [yshift = -1.5pt, xshift = -1.5pt] rectangle ++(3pt, 3pt) node[anchor = west] {\color{blue} $\fp_\pm^f$}; ;

		\filldraw[Green] (-2, 0) circle(2pt);
		\node at (-2.1, 0.4) {\color{Green} $\fb_\pm^f$};
		
		\draw[Green, thick, ->, > = stealth] (-2, 0) -- (-3, 0);

	\end{tikzpicture}
	\end{subfigure}

\begin{subfigure}{0.4\textwidth} 
	\centering
	\begin{tikzpicture}[scale = 0.8]
		\node at (-5, 0) {$\longrightarrow$}; 
		
		\draw[thick] (-3.3, 1.5) -- (-1.5, 2)--(0,1)--(0, -1) -- (-1.5, -2) -- (-3.3, -1.8);
		\draw[thick] (-1.5, 2)--(-1.4, 3) ;
		\draw[thick] (0, 1) -- (1, 1.3);
		\draw[thick] (0, -1) -- (1, -1.4);
		\draw[thick] (-1.5, -2) -- (-1.2, -2.8);
		
		 \node at (0, 1.3) {$v_1$};  
		\node at (0.1, -1.3) {$v_n$};
		\node at (-1.7, 2.15) {$v_2$};
		\node at (-0.7, 1.7) {$e_1$};
		 
		\node at (0.3, 0) {$e_n$}; 
		\node at (-2.8, 1.9) {$e_2$};
		
		\draw[red, thick] (-2, 0) .. controls (-4, 0) and (-2.5, 1.5) .. (-1.5, 2);
		\node at (-3, 1.1) {\color{red} $l_2$};

		\draw[red, thick] (-2, 0) .. controls (-3.5, 0) and (-2.2, -2.4).. (0, 1);
		\node at (-1.3, -.9) {\color{red} $l_1$};

		\filldraw[blue] (-1.5, 0.2) [yshift = -1.5pt, xshift = -1.5pt] rectangle ++(3pt, 3pt) node[anchor = west] {\color{blue} $\fp_\pm^f$}; ;

		\filldraw[Green] (-2, 0) circle(2pt);
		\node at (-1.9, -0.3) {\color{Green} $\fb_\pm^f$};
		
		 \draw[Green, thick, ->, > = stealth] (-2, 0) -- (-3, 0); 
		
	\end{tikzpicture}
\end{subfigure}
\begin{subfigure}{0.4\textwidth} 
	\centering
	\begin{tikzpicture}[scale = 0.8]
		\node at (-5, 0) {$\longrightarrow$}; 
		
		\draw[thick] (-3.3, 1.5) -- (-1.5, 2)--(0,1)--(0, -1) -- (-1.5, -2) -- (-3.3, -1.8);
		\draw[thick] (-1.5, 2)--(-1.4, 3) ;
		\draw[thick] (0, 1) -- (1, 1.3);
		\draw[thick] (0, -1) -- (1, -1.4);
		\draw[thick] (-1.5, -2) -- (-1.2, -2.8);
		
		 \node at (0, 1.3) {$v_1$};  
		\node at (0.1, -1.3) {$v_n$};
		\node at (-1.7, 2.15) {$v_2$};
		\node at (-0.7, 1.7) {$e_1$};
		 
		\node at (0.3, 0) {$e_n$}; 
		\node at (-2.8, 1.9) {$e_2$};
		
		\draw[red, thick] (-2, 0) .. controls (-2.5, -.6) and (-3.5, 0) .. (-1.5, 2);
		\node at (-2.7, 1) {\color{red} $l_2$};

		\draw[red, thick] (-2, 0) .. controls (-2.5, -.6) and (-1.5, -2).. (0, 1);
		\node at (-1.3, 1.1) {\color{red} $l_1$};

		\filldraw[blue] (-1.5, 0.2) [yshift = -1.5pt, xshift = -1.5pt] rectangle ++(3pt, 3pt) node[anchor = west] {\color{blue} $\fp_\pm^f$}; ;

		\filldraw[Green] (-2, 0) circle(2pt) node[anchor = south] {\color{Green} $\fb_\pm^f$};
		
		 \draw[Green, thick, ->, > = stealth] (-2, 0) -- (-2.5, -0.6); 
		
	\end{tikzpicture}
\end{subfigure}

	\caption{Steps in the proof of Lemma \ref{lem: framing-independence} showing independence of framing points.}
	\label{fig: framing of base pt}
\end{figure}

         It suffices to show independence if we change the distinguished edge from $e_n$ to $e_1$ (or \emph{vice versa}). 
         For the proof, we refer to Figure \eqref{fig: framing of base pt}.      
         The map taking the first picture to the second is an isotopy. From the second picture to the third one, we get a factor of $(-a)^{-2}$ when passing through $\fp_\pm^f$,  by Equation \eqref{iso of skeins}. This factor cancels out with the twists we get when passing through $\fb_\pm^f$, and the passage from the third to final picture is the rotation isotopy.
\end{proof}

\subsection{Mutation}\label{subsec:mutation}

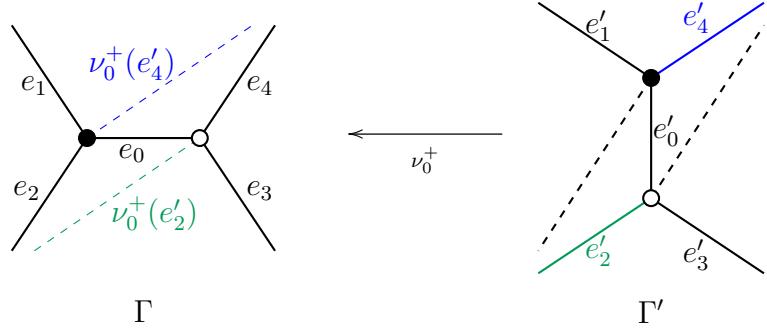
\begin{figure}[htbp]
        \centering
    \begin{tikzpicture}[scale = 1]
        
        \draw[thick] (-2.5,-1.5)--(-1.5,0)--(-2.5,1.5);
	\draw[thick] (-1.5,0)--(0,0);
	\draw[thick] (1,-1.5)--(0,0)--(1,1.5);

        \draw[blue, dashed] (-1.5, 0) -- (0.7, 1.5); 
        \draw[ForestGreen, dashed] (-2.2, -1.5) -- (0, 0);
        
        \filldraw[thick, black ] (-1.5, 0) circle (3pt); 
        \filldraw[thick, black, fill = white] (0, 0) circle (3pt); 

        \node at (-2.2, 0.7) {$e_1$};
        \node at (-2.3, -0.7) {$e_2$}; 
        \node at (-0.9, -0.2) {$e_0$};
        \node at (0.8, 0.7) {$e_4$}; 
        \node at (0.8, -0.7) {$e_3$};
        \node at (-0.9, 1)  {\color{blue} $\nu_0^+(e_4')$};
        \node at (-0.6, -1) {\color{ForestGreen} $\nu_0^+ (e_2')$};

        \node at (3, 0) {$\xymatrix{{}\\ {}&&\ar[ll]^{\nu^+_0}{}\\{}}$};

        \draw[thick] (4.5, 1.8) -- (6, 0.8) -- (6, - 0.8) -- (7.5, -1.8) ;
        \draw[thick, blue] (6, 0.8) -- (7.5, 1.8); 
        \draw[thick, ForestGreen] (4.5, -1.8)-- (6, -0.8); 

        \draw[thick, dashed] (6, 0.8) -- (4.5, -1.5); 
        \draw[thick, dashed] (6, -0.8) -- (7.5, 1.5); 
        
        \filldraw[thick, black ] (6, 0.8) circle (3pt); 
        \filldraw[thick, black, fill = white] (6, -0.8) circle (3pt); 

        \node at (5.3, -1.5) {\color{ForestGreen}$e_2'$}; 
        \node at (6.6, -1.55) {$e_3'$};
        \node at (5.3, 1.5) {$e_1'$};
        \node at (6.6, 1.6) {\color{blue} $e_4'$}; 
        \node at (6.2, .1) {$e_0'$};

        \node at (- 0. 75, -2.3) {$\Gamma$};
        \node at (6, -2.3) {$\Gamma'$};
        
    \end{tikzpicture}
        \caption{Positive mutation homeomorphism}
        \label{fig:skein mutation}
    \end{figure}

Let $\Gamma$ and $\Gamma'$ be two graphs which differ by a mutation, as in Figure \ref{fig:skein mutation}.  In this section we will define a notion of ``skein mutation'' relating the skein modules $\widehat{\Sk} (S_{\Gamma'}, \fb'_\pm; \fp'_\pm) $ and $\widehat{\Sk}(S_{\Gamma}, \fb_\pm; \fp_\pm)$, and show that it preserves the face relations. 

We first define a skein-theoretic analog of the quantum dilogarithm of a cluster monomial.  


For a knot $l$ in a $3$-manifold, its tubular neighborhood is a solid torus. We can use the decoration in \cite[Sections 2.4 and 2.5]{MS17} to define skein algebra elements $\skein{nl}$. By \cite[Theorem 3.2]{MS17} we know that this ``coloring'' defines Adams operations on $\Sk(T^2)$. However, for a general surface this may not be an algebra homomorphism\footnote{ For example, in the skein algebra of the punctured torus, it is not true that $[P_{(2, 0)}, P_{(0,2)}] = {\{2\}} P_{(2,2)}$. }. 

\begin{definition} \label{def: Baxter operator}
    Let $\Sigma$ be an oriented surface (possibly with framing lines) and let $l \in \Sk(\Sigma)$ be a positive simple closed curve. We define the \emph{Baxter operator} to be the following element of the completed skein algebra $\widehat\Sk(\Sigma):$
	\begin{equation*}
		Q_l (t)  \coloneqq  \e \left(\sum_{n = 1}^{\infty}\frac{ -(-t)^{n} \skein{nl}}{n (q^{n/2} - q^{-n/2} ) } \right) \quad\text{or} \quad Q_l (-t) =  \E \left(\frac{-t \skein{l}}{q^{1/2} - q^{-1/2}} \right), 
	\end{equation*}
where $t$ is viewed as a line element in the plethystic exponential. 
\end{definition}
We write $Q_l$ for $Q_l(1)$. 

\begin{remark}
Our use of the term `Baxter operator' is based on the fact that $Q_l(t)$ may be characterized as the solution of the $q$-difference equation
\begin{align}
\label{eq:bax-qde}
Q_{l}(qt) = \Gamma_+(t)Q_{l}(t),
\end{align}
where $\Gamma_+(t)$ is the vertex operator
$$
\Gamma_+(t) = \exp\left(-\sum_{n>0} \frac{\left(-q^{\frac{1}{2}}t\right)^n}{n}P_{nl}\right).
$$
\end{remark}
We now establish some relations satisfied by Baxter operators in the skein algebra of the closed torus $T^2$. 
For $\x$, $\y \in \bZ^2$, denote the elements they represent in $\Sk(T^2)$ by $P_\x$ and $P_\y$.
Set $d:= \det(\x | \y) \in \bZ$, and put $s = q^{1/2}$, so that $\{k\} \coloneqq  s^k - s^{-k}$. 
Recall the (symmetric) quantum numbers $$
[n]_q \coloneqq \frac{s^n-s^{-n}}{s-s^{-1}} = s^{n-1}+s^{n-3}+\ldots+s^{-{n-1}},
$$
as well as the quantum factorials $[n]_q \coloneqq \prod_{m=1}^n [m]_q$ and binomial coefficients
$$
\ \qbinom{n}{k} \coloneqq \frac{\prod_{m=0}^{k-1} [n-m]_q}{[k]_q!}
$$
obeying the $q$-binomial theorem:  with $n\geq 0$ we have
    \begin{align*}
        \sum_{k\ge 0}  \qbinom{n}{k} x^k 
        &= (1 + x s^{n-1}) (1 + x s^{n-3}) (1 + x s^{n-5}) \cdots (1 + x s^{-(n-1)}).
    \end{align*}
    We extend their definition to $n<0$ by setting
    \begin{align*}
    \sum_{n \ge 0}  \qbinom{-n}{k} x^k 
        &= \frac{1}{ (1 + x s^{n-1}) (1 + x s^{n-3}) \cdots (1 + x s^{-(n-1)}) }. 
    \end{align*}
Then we have:
\begin{lemma} The following identity holds in $\widehat\Sk(T^2)$:
\label{lem:Ad-action}
    \begin{align}\label{adjoint action}
        \Ad_{Q_\x(t)} \skein{\y} = \sum_{n\ge 0} \ \qbinom{d}{n}   t^n     
        \skein{\y+n\x} 
    \end{align}
  %

In particular, we have:
\begin{align}
    \Ad_{Q_\x} \skein{\y} 
    &= \skein{\y} + \skein{\x + \y} , \quad \text{if} \ d = 1;\\
    \Ad_{Q_\x} \skein{\y} 
    &=  \sum_n (-1)^n \skein{n\x + \y} \quad \text{if} \ d = -1. 
\end{align}
\end{lemma}

\begin{proof}

Set $\Theta_{\x}(t) \coloneqq  \sum_{k\geq1} \frac{(-1)^{k+1} t^k }{k\{k\}} P_{k\x}$, so  that $Q_\x (t) = \e (\Theta_\x(t))$. We have:

$$
\ad_{\Theta_\x(t)} P_\y = \sum_{k \geq 1} \frac{\{kd\}}{k\{k\}} (-1)^{k+1} t^k P_{k \x + \y }
$$
and 
\begin{align*}
	\Ad_{Q_\x(t)} \cdot P_\y 
	&= e^{\ad_{\Theta_\x}(t)} \cdot P_\y \\
	&= \sum_N \frac{1}{N!} (\ad_{\Theta_\x}(t))^N P_\y \\
	&= \sum_n \sum_{k_1 + \cdots +k_N = n} \frac{1}{N!} (-1)^{N+n} t^n \prod_i \frac{\{k_i d\}}{k_i \{k_i\}} P_{n\x + \y}                
\end{align*}

We show that 
$$\sum_{k_1 + \cdots +k_N = n} \frac{1}{N!} (-1)^{N+n} \prod_i \frac{\{k_i d\}}{k_i \{k_i\}} =  \qbinom{d}{n}.
$$
If $d \ge 0$, we have:
\begin{align*}
	\sum_k \frac{\{kd\}}{k\{k\}} (-1)^{k+1} x^k
	&= \sum_k \frac{1}{k} (s^{(d-1)k } + s^{(d-3)k} + \cdots s^{-(d-1)k} ) (-1)^{k+1} x^k \\
	&= \l(1+ x s^{d-1}) + \cdots + \l(1+ x s^{-(d-1)})	
\end{align*}
Thus 
\begin{align*}
        \sum_n \sum_{k_1 + \cdots +k_N = n} \frac{1}{N!} (-1)^{N+n} \prod_i \frac{\{k_i d\}}{k_i \{k_i\}}   x^n 
	&= \e(	\sum_k \frac{\{kd\}}{k\{k\}} (-1)^{k+1} x^k)\\
	&= (1 + xs^{d-1}) (1 + xs^{d-3} ) \cdots (1 + s^{-(d-1)}),
\end{align*}
which proves the claim. 
The case $d < 0$ is similar. 
\end{proof}
\begin{remark}
   By a similar computation, we have:
\begin{align*}
    \Ad_{Q^{-1}_\x(t)} P_\y = \sum_{n\ge 0} \qbinom{-d}{n} t^n P_{\y + n \x}
\end{align*}
    
\end{remark}

Now we focus on two simple loops $\x$ and $\y$ in a genus-$g$ surface, intersecting positively in a single point ---
so the symplectic pairing is $\det(\x|\y)=1$ ---
whose tubular neighborhood is a punctured torus $T^2\setminus D^2$.
We use the notation $\x+\y$ to indicate the loop defined by the connect sum
\begin{center}
\begin{tikzpicture}
    \draw[thick,black] (0,0)--(1,0)--(1,1)--(0,1)--(0,0);
    \filldraw[thick,black,fill=white] (0,0) circle (2pt);
    \filldraw[thick,black,fill=white] (1,0) circle (2pt);
    \filldraw[thick,black,fill=white] (1,1) circle (2pt);
    \filldraw[thick,black,fill=white] (0,1) circle (2pt);
    \draw[thick,blue,->,>=stealth'] (0,1/2) .. controls (1/2,1/2) .. (1/2,1);
    \draw[thick,blue,->,>=stealth'] (1/2,0) .. controls (1/2,1/2) .. (1,1/2);
\end{tikzpicture}
\end{center}
and similarly for $a\x + b\y,$ which in the
non-primitive case needs to be decorated
as in \cite[Sections 2.4 and 2.5]{MS17}.
Note that $\x+\y$ are $\y+\x$ are the same skein,
and are not to be confused
with the
compositions in the nonabelian fundamental group.
Note, too, that the defining relations of
the skein algebra are local,
so any skein relations involving operators made from $\x$
and $\y$ that hold in the punctured torus
will hold in $\Sk(\Sigma)$ as well.

An element ${\bf a} = (a_1, a_2) \in \bZ^2$ determines an element in the skein module of this punctured torus. 
As is pointed out in \cite[Remark 4.3]{MS21}, if $\bf a$ is primitive and ${\bf b} = n \bf b_0$ where $\det ({\bf a |  b_0}) = \pm 1$, then 
the following  relations still holds in the punctured torus skein:
\begin{equation}
\label{eq:ms-skein-equation 1}
    [P_{\bf a},P_{\bf b}] = \{\omega({\bf a},{\bf b})\}\, P_{\bf a+b},
\end{equation}
where as above $\{k\}:= s^k-s^{-k},$ and $s = q^\frac{1}{2}$.

More generally, if $\y$ is a path in the skein module with base points $\Sk(S_\Gamma, \fb)$, transversely intersecting the loop $\x$ at a single point, 
the skein $\y + n\x$ is still defined, and \eqref{eq:ms-skein-equation 1} 
is still true in the following sense (\cite[Theorem 4.2]{MR1915490})
\[
[\skein{n\x}, \skein{\y}] = \{n \omega(\x, \y)\} \skein{n\x + \y}, \quad n \in \bZ .
\]
In this case, note that $\widehat{\Sk} (S_\Gamma, \fb)$ can be viewed as a $\widehat{\Sk} (S_\Gamma - \fb)-\widehat{\Sk}(S_\Gamma)$ bimodule. Thus $\Ad_{Q_\x} \skein{y}$ still makes sense, and the conclusion of Lemma \ref{lem:Ad-action} holds for such $\x,\y$. 
 

Now we define the skein version of the ``signed mutation'' $\nu^\pm$, which are isomorphisms between $\Sk (S_{\Gamma'}, \fb'_\pm; \fp'_\pm) $ and $\Sk(S_{\Gamma}, \fb_\pm; \fp_\pm)$, induced by homeomorphisms between the two surfaces. By abuse of notation, we denote both the isomorphism between skein modules and the homeomorphism between surfaces by $\nu^\pm$.

We first consider a piecewise linear map between the neighborhoods of edge $e_0'$ and edge $e_0$. As in Figure \ref{fig:skein mutation}, each neighborhood is divided into six regions. The homeomorphism $\nu^+$ is defined by  ``counterclockwise twisting''. Each region is mapped linearly onto the corresponding one. This induces a homeomorphism on the branched double cover, i.e.~the cylinder neighborhood of the loop $\loo{e}_0'$ and the cylinder neighborhood of the loop $\loo{e}_0$. Hence we get a homeomorphism between $S_{\Gamma'}$ and $S_\Gamma$. 

Under this homeomorphism $\nu^+$, the loops are mapped to:
\begin{align*}
    \skein{e_0'} \longmapsto \skein{-e_0}\\
   \skein{e_1'} \longmapsto P_{e_1} & & \skein{e_4'} \longmapsto \skein{e_0+e_4}\\
   \skein{e_2'} \longmapsto \skein{e_2 + e_0} & &\skein{e_3'} \longmapsto \skein{e_3}
\end{align*}
where $\skein{e_0 + e_2}$ means the skein represented by $\loo{e}_0 + \loo{e}_2$. Thus upon projecting to the linking skein $Lk(S_\Gamma)$, this definition of $\nu^+$ recovers the one in \eqref{eq:torusmutation}. 

\begin{figure}[htbp]
        \centering
    \begin{tikzpicture}[scale = 1]
        
        \draw[thick] (-2.5,-1.5)--(-1.5,0)--(-2.5,1.5);
	\draw[thick] (-1.5,0)--(0,0);
	\draw[thick] (1,-1.5)--(0,0)--(1,1.5);

        \draw[blue, dashed] (0, 0) -- (-2.2, 1.5); 
        \draw[ForestGreen, dashed] (0.7, -1.5) -- (-1.5, 0);
        
        \filldraw[thick, black, fill = white ] (-1.5, 0) circle (3pt); 
        \filldraw[thick, black] (0, 0) circle (3pt); 

        \node at (-2.2, 0.7) {$e_1$};
        \node at (-2.3, -0.7) {$e_2$}; 
        \node at (-0.9, -0.2) {$e_0$};
        \node at (0.8, 0.7) {$e_4$}; 
        \node at (0.8, -0.7) {$e_3$};
        \node at (-0.7, 1.1)  {\color{blue} $\nu_0^-(e_1')$};
        \node at (-0.6, -1.2) {\color{ForestGreen} $\nu_0^- (e_3')$};

        \node at (3, 0) {$\xymatrix{{}\\ {}&&\ar[ll]^{\nu^-_0}{}\\{}}$};

        \draw[thick] (7.5, 1.8) -- (6, 0.8) -- (6, - 0.8) -- (4.5, -1.8) ;
        \draw[thick, blue] (6, 0.8) -- (4.5, 1.8); 
        \draw[thick, ForestGreen] (7.5, -1.8)-- (6, -0.8); 

        \draw[thick, dashed] (6, -0.8) -- (4.5, 1.5); 
        \draw[thick, dashed] (6, 0.8) -- (7.5, -1.5); 
        
        \filldraw[thick, black ] (6, 0.8) circle (3pt); 
        \filldraw[thick, black, fill = white] (6, -0.8) circle (3pt); 

        \node at (5.3, -1.55) {$e_2'$}; 
        \node at (6.6, -1.55) {\color{ForestGreen}$e_3'$};
        \node at (5.3, 1.55) {\color{blue}$e_1'$};
        \node at (6.6, 1.6) { $e_4'$}; 
        \node at (6.2, -.1) {$e_0'$};
        
        \node at (- 0. 75, -2.3) {$\Gamma$};
        \node at (6, -2.3) {$\Gamma'$};
 
    \end{tikzpicture}
        \caption{Negative mutation homeomorphism}
        \label{fig:skein negative mutation}
    \end{figure}
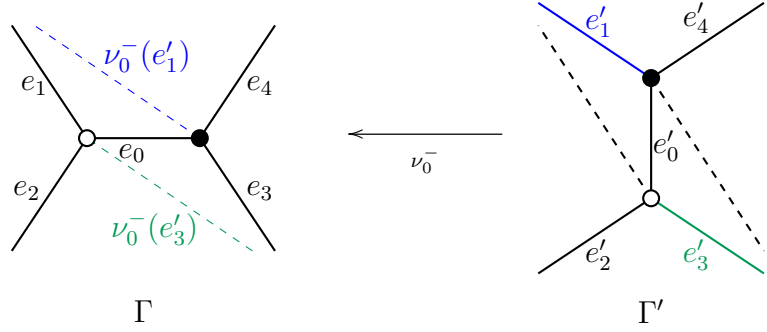

Similarly, we can take a ``clockwise twist'', as in figure \eqref{fig:skein negative mutation}, to define $\nu^-$.

    
We are now ready to define a skein version of cluster transformations. Let $\Gamma$ and $\Gamma'$ again be the two graphs in Figure \ref{fig:skein mutation}. 
\begin{definition}
    \label{eq:mutation-def}
We define the mutation at the edge $e_0$ to be the homomorphism  $\mu_0:\widehat{\Sk}(S_{\Gamma'}, \fb'_\pm; \fp'_\pm) \longrightarrow \widehat{\Sk}(S_\Gamma, \fb_\pm; \fp_\pm)$ given by:
    \begin{align} \label{\skein mutation}
		\mu_0 &\coloneqq \Ad_{Q_{e_0}} \circ \nu^+  \quad \text{if } e_0\  \text{is ``positive'', or}\\
     \mu_0 &\coloneqq \Ad_{Q_{-e_0} ^{-1}} \circ \nu^-  \quad \text{if } e_0\  \text{is ``negative''.} \notag
    \end{align}
\end{definition}
We can now prove the following skein analog of the globalization result Theorem~\ref{thm:d-mod}:
\begin{theorem}\label{thm:face relation}
    The face relation \eqref{eq:skein-face-relation} is preserved under the cluster mutation \eqref{\skein mutation}. 
\end{theorem}

\begin{proof}
     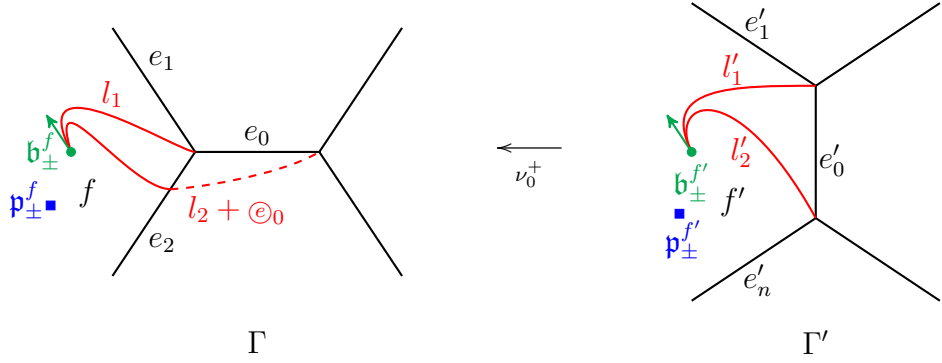
\begin{figure}
	\centering
	\begin{tikzpicture}[scale = 1.1]
		\draw[thick] (-2.5,-1.5)--(-1.5,0)--(-2.5,1.5);
		\draw[thick] (-1.5,0)--(0,0);
		\draw[thick] (1,-1.5)--(0,0)--(1,1.5); 
		
		\node at (-1.9, 1.1) {$e_1$};
		\node at (-1.9, -1.1) {$e_2$}; 
		\node at (-0.75, 0.2) {$e_0$};
		
		\draw[thick, red] (-3, 0) .. controls (-3.6, 1.2) and (-1.7, 0) .. (-1.5, 0); 
		\node at (-2.5, 0.7) {\color{red}$ l_1$}; 
	
		\draw[red, thick] (-3, 0) .. controls (-3.3, 1) and (-2.3, - 0.45) .. (-1.8, -0.45); 
		\draw[red, thick, dashed] (-1.8, -0.45) .. controls (-1.6, -0.45) and (-0.5, -0.3) .. (0, 0); 
		\node at (-1, -0.7) {\color{red} $l_2 + \loo{e}_0$};
		
		\filldraw[Green] (-3, 0) circle(1.5pt) node[anchor = east]{$\fb_\pm^f$};
		\draw[Green, thick, ->, > = stealth'] (-3, 0) -- (-3.3, 0.45) ; 
		
		\filldraw[blue] (-3.3, -.7)  rectangle ++(3pt, 3pt) node[anchor=east] {\color{blue} $\fp_\pm^f$}; 
		\node at (-2.8, -0.5) {$f$};

		\node at (2, 0) {$\xymatrix{{}\\ {}&&\ar[l]^{\nu^+_0}{}\\{}}$};
		
		\draw[thick] (7.5, 1.8) -- (6, 0.8) -- (6, - 0.8) -- (4.5, -1.8) ;
		\draw[thick] (6, 0.8) -- (4.5, 1.8); 
		\draw[thick] (7.5, -1.8)-- (6, -0.8); 
		
		\node at (5.3, -1.55) {$e_n'$}; 
		\node at (5.3, 1.55) {$e_1'$};
		\node at (6.2, -.1) {$e_0'$};
		
		\draw[red, thick] (4.5, 0) .. controls (4.1, 0.8) and (5, 0.8) .. (6, 0.8);
		\node at (5, 1) { \color{red}$l_1'$}; 
		\draw[red, thick] (4.5, 0) .. controls (4.2, 0.6) and (5, 1) .. (6, -0.8); 
		\node at (5.1, -0) {\color{red} $l_2'$}; 
		
		\node at (5, -.6) {$f'$};
		
		\filldraw[blue] (4.3, -.8) rectangle ++(3pt, 3pt) node[anchor=north]{$\fp_\pm^{f'} $};
		\filldraw[Green] (4.5, 0) circle(1.5pt) node[anchor = north] {$\fb_\pm^{f'}$};

		\draw[thick, Green, ->, >= stealth'] (4.5, 0) -- (4.2, 0.45);

		\node at (- 0. 75, -2.3) {$\Gamma$};
		\node at (6, -2.3) {$\Gamma'$};
		
	\end{tikzpicture}
	\caption{Face relations before and after mutation} 
	\label{fig:mutation of face rel}
\end{figure}

    Consider the face relations for $f'$ and $f$ as in Figure \ref{fig:mutation of face rel}. Assume that $\fb_\pm^{f'}$ and $\fb_\pm^f$ are framed as shown, and that $e_0'$ is a positive edge. The map $\nu_0^+$ sends all $\frame{l_i'}$ to $\frame{l_i}$, with the exception of $\frame{l_0'}$ which is sent to $\frame{l_1 + \loo{e}_0} $. So we have:
		\begin{align*}
			\mu_0 (A_{\Gamma', f'})
			& = \Ad_{Q_{e_0}} \circ \nu_0^+ (\frame{l'_0} + \frame{l'_1} + \cdots + \frame{l'_n}) \\
			& = \Ad_{Q_{e_0}} (\frame{l_1 + e_0 } + \frame{l_1} + \cdots + \frame{l_n})\\
			& = \Ad_{Q_{e_0}}(\frame{l_1 + e_0} + \frame{l_1}) + \frame{l_2} + \cdots + \frame{l_n}\\
			& \overset{Lemma \eqref{lem:Ad-action}}{=} \frame{l_1} + \cdots + \frame{l_n}  \\
			& = A_{\Gamma, f} 
 		\end{align*}
    The case where $e_0'$ is negative is similar. 

\end{proof}

\section{Skein Scattering and the Pentagon Relation}
\label{sec:skein-scattering}

Much of the magic of cluster theory can be attributed to the
pentagon relation obeyed by the quantum dilogarithm:
indeed, the noncommutative cluster mutation is effected by the adjoint action
of the dilogarithm.

At the level of wavefunctions of branes --- or modules for the
quantum torus algebra --- the quantum dilogarithm acts multiplicatively, \`a la Schr\"odinger representation,
to effect mutation.
As shown in \cite{SSZ}, mutation of the Chekanov torus to the
Clifford torus takes us from the constant wavefunction $1$ to the
wavefunction of the Aganagic-Vafa brane.  In \cite{ES19} this wavefunction
was generalized to skein modules.  

In Section \ref{sec:skein-pentagon} we discuss the operator on the
skein algebra analogous to the dilogarithm, then show 
that it obeys the pentagon relation.  First we briefly recall
in Section \ref{sec:scattering} that
the cluster mutation is but one form of wall crossing encountered in
a scattering diagram, and that scattering diagrams make sense
for the skein algebra.  This leads to a natural guess for
the skein version of the dilogarithm.

\subsection{Scattering Diagrams}
\label{sec:scattering}

This section is motivational and can be read independently, or not at all.
\vskip0.1in
Scattering diagrams and their wall functions
encode the combinatorics of automorphisms which
glue  charts of a (not necessarily commutative) space.
The wall functions are created from counts of objects, depending
on the setting --- examples are holomorphic maps, tropical curves, counts of objects
in a category.
Scattering diagrams in all of these different settings are constructed from
a lattice-graded Lie algebra, with the automorphisms constructed from a
corresponding pro-nilopotent Lie algebra --- see \cite{Ma21} for a clear exposition,
with examples.

For our purposes, it suffices to consider the following setting.
Let $N = \bZ^2$, equipped with the
standard antisymmetric pairing $\omega(a,b)=\det(a|b)$.

This leads to three classes of examples, in increasing generality:
\subsubsection{Symplectic torus}

The associated complex torus $T := \Hom(N,\bC^*) \cong (\bC^*)^2$
is symplectic, so $\mathfrak{g} := \cO_{T}$ has a Lie algebra structure defined
by the Poisson bracket $\{z^a,z^b\}=d\cdot z^{a+b},$
where $d = \omega(a,b).$  Walls are defined by
primitive ray vectors $m\in \bZ^2$ and wall functions
are of the form $f = 1 + \sum_{k>0}c_k \cdot z^{km}.$ 
The corresponding automorphism is $z^n\mapsto f^d \cdot z^n,$ $d = \omega(m,n)$.
We define for this wall $H = \sum H_k z^{km}$ by $\log f = \sum k H_k z^{km}$.  Note $\{H,z^n\} = \sum H_k kd z^{n+km}$ so
the time-one flow is
$$\mathrm{Ad}_{e^H}\cdot z^n = e^{\{H,-\}}\cdot z^n = f^d\cdot  z^n.$$

\subsubsection{Quantum torus and
\texorpdfstring{$q$}{q}-deformation}

We keep the notations as above, including the wall function $f$.
Recall the defining relations for the quantum torus
$$z^m z^n = q^{d} z^n z^m = s^{d} z^{n+m}, \qquad d = \omega(m,n),\qquad s = q^{\frac{1}{2}}.$$
In particular,
$$\qquad ad_{z^m}\cdot z^n = [z^m,z^n] = (q^d - 1)z^n\cdot z^m = (s^{d}-s^{-d})z^{m+n} =
\{d\}z^{m+n} = \{d\}T_m \cdot z^n,\qquad.$$
where we have defined $\{d\} := s^d - s^{-d}$ and the operator $T_a: z^b\mapsto z^{a+b}.$

Again for a wall with primitive ray vector $m$, we define $H = \sum H_k z^{km},$ where $H_k$ are defined by $\log f(z) = \sum H_k(q^k-1)z^k.$
If we write $q = e^{i\hbar}$, then up to $O(\hbar^2)$ this is indeed $i\hbar$ times the above classical
Hamiltonian.

Define $f_j$ by $f_j(z^m) = f(q^jz^m) = \sum c_k q^{jk}z^{km}.$
Then if $d>0$ we have
$$Ad_{e^H}\cdot z^n = \Bigl(\,\prod_{j=0}^{d-1} f_j\, \Bigr) \cdot z^n$$
(there is a similar formula for $d<0$).  Following \cite[Lemma 3.3]{B18} we prove this
as follows:
$$
\begin{array}{rcl}
Ad_{e^H}\cdot z^n &=& e^{ad_H}\cdot z^n = \exp\left(\sum_k H_k ad_{z^{km}}\right)\cdot z^n\\ &=&
 \exp\left(\sum_k H_k (q^{kd}-1)z^{km}\right)\cdot z^n\\
&=& \exp\left(\sum_k H_k \frac{(q^{kd}-1)}{(q^k-1)}(q^k-1) z^{km}\right)\cdot z^n\\
&=& \bigl( \,\prod_{j=0}^{d-1}f(q^j z^m)\, \bigr)\cdot z^n = \bigl(\, \prod_{j=0}^{d-1}f_j\,\bigr) \cdot z^n\\
\end{array}$$
where in the last line we use $\frac{(q^{kd}-1)}{(q^d-1)} = \sum_{j=0}^{d-1}q^{kj}$
and the relation between $H$ and $\log f.$  Note this reduces to the commutative case when $q\to 1.$

Let's write this in another form:
$$
\begin{array}{rcl}
Ad_{e^H}\cdot z^n &=& e^{ad_H}\cdot z^n = \exp\left(\sum_k H_k ad_{z^{km}}\right)\cdot z^n\\ &=&
\exp\left(\sum_k H_k \{kd\}T_{km}\right)\cdot z^n\\
&=& \exp\left(\sum_k H_k \frac{(s^{kd}-s^{-kd})}{(s^{2k}-1)}(q^k-1) (T_m)^k \right)\cdot z^n\\
&=& \prod_{j=1}^{d}f(q^{-j} T_m)\cdot z^n\\
\end{array}$$

\begin{example}
Suppose $d=1$ and $f = 1 + sz^m.$  Then $H = \sum_k H_k z^{km} = \frac{(-1)^{k+1}}{k(q^k-1)}(sz^m)^{k},$
so $$e^H = Q_m := Q_m(1), \qquad Q_m(-a) := \mathrm{Exp}\left(\frac{-az^m}{s-s^{-1}}\right).$$
In writing this, we recall that $\mathrm{Exp}$ is defined in the setting of a $\lambda$-ring with $n$th
Adams operation $\psi^n$, and above we have taken $\psi^n: z \mapsto z^n$, $s \mapsto s^n$,
i.e. $z$ and $s$ are line elements.  

\end{example}

\subsubsection{Skein Algebra}
\label{sec:skein-baxter}

Recall the Morton-Samuelson relation in the skein algebra of a torus:
$$[P_m,P_n] = \{d\}P_{m+n} = \{d\}T_m \cdot P_n,$$
where the operator $T_m$ is defined similarly as above: $T_m: P_n \mapsto P_{m+n}.$
We see that the skein algebra of a torus is thus naturally $\bZ^2$-graded.
So from a wall function $f$ as above, we can construct the quantum Hamiltonian $H$ as above
then graduate to a skein Hamiltonian simply by the replacement $z^k \rightsquigarrow P_k$.

\begin{example}
Let $m$ be a primitive vector along a wall.  Let $f = 1 + sP_m.$  The
Hamiltonian is then $H = \sum_k H_k P_{km} = \frac{(-1)^{k+1}s^k}{k(q^k-1)}P_{km},$
and therefore
\begin{equation}
\label{eq:plethbaxter}
    e^H = Q_m := Q_m(1),\qquad Q_m(-a) = \mathrm{Exp}\left(\frac{-aP_m}{s-s^{-1}}\right),
\end{equation}
where we have defined the Adams operation $\psi^n(P_m) = P_{nm}$ for skeins on a torus, as in Section \ref{sec:skt2}. 
This definition agrees with Definition \ref{def: Baxter operator},
though elsewhere in the paper when we discuss the skein of a general surface, the
Adams operation is not clear.
In Section \ref{sec:skein-pentagon} we prove this operator obeys the pentagon relation ---
see Remark \ref{rmk:pentagonrelation}.
\end{example}

\subsubsection{Consistent scattering diagrams}

A scattering diagram is a collection of rational polyhedral walls, each equipped with
a wall function and corresponding automorphism, as above.
To each oriented loop $\gamma$ in $M_\bR$ one associates an automorphism $A_\gamma$ obtained by composing, in order, those
of each wall met by the loop.
A diagram is said to be
\emph{consistent} if for all loops $\gamma$, the automorphism $A_\gamma$ is the identity.

Consistent diagrams are constructed by beginning with an initial diagram and
adding rays ``at higher order'' to account for any discrepancies in the automorphisms associated
to two paths with the same endpoints.  This introduces new discrepances, which are fixed at higher
and higher order.  

\begin{remark}

It would be interesting to implement the above procedure to a complex
scattering diagram, such as the one associated to a smooth cubic curve $E\subset \bP^2$ in \cite{GRZ}.
Work in progress shows that quantum wall crossings encode higher-genus log invariants of $(\bP^2,E)$
with one boundary component.
We would hope that the skein version should encode multiple boundary components, as well.
\end{remark}

\begin{example}
\label{ex:basicdiagram}
$$
\begin{tikzpicture}
\draw[thick] (-1,0)--(1,0);
\draw[thick] (0,-1)--(0,1);
\draw[thick] (0,0)--(1,1);
\node[left] at (-1,0) {$1+z^{(1,0)}$};
\node[below] at (0,-1) {$1+z^{(0,1)}$};
\node[above right] at (1,1) {$1+z^{(1,1)}$};
\draw[->,>=stealth',blue, thick] (305:.8cm) arc[radius=.8, start angle=305, end angle=145];
\draw[->,>=stealth',blue, thick] (325:.8cm) arc[radius=.8, start angle=-35, end angle=125];

\end{tikzpicture}
$$
As explained above, the indicated wall functions generate scattering diagrams in three contexts.
Consistency of the diagrams for the classical and quantum torus is ensured by the pentagon
relation for dilogarithms: $\Phi(z^{(0,1)})\Phi(z^{(1,0)})=\Phi(z^{(1,0)})\Phi(z^{(1,1)})\Phi(z^{(0,1)}).$  We now prove consistency for the skein algebra diagram.

\end{example}
%

\subsection{The Pentagon Relation for Skeins}
\label{sec:skein-pentagon}


    Let $\Sigma$ be an oriented, genus-$g$ surface (possibly with framing lines), for example $S_\Gamma$.  Let $\x$ be an oriented, simple closed curve.
    Let $\Theta_{\x}(t) \coloneqq \sum_{k\geq1} \frac{(-1)^{k+1} t^k }{k\{k\}} P_{k\x}$.
   Recall in \eqref{def: Baxter operator} we defined 
    \begin{equation} Q_\x(t) = \e (\Theta_\x(t)),\qquad\text{or\qquad} Q_\x(-t) = \mathrm{Exp}\left(\frac{-t}{s-s^{-1}}P_\x\right)
    \end{equation}
    and set $Q_\x = Q_\x(1).$

To investigate the pentagon relation, as in section \ref{subsec:mutation}, take two simple loops $\x$ and $\y$ in $\Sigma$, intersecting positively in a single point ---
so the symplectic pairing is $\omega(\x,\y)=1$ ---
whose tubular neighborhood is a punctured torus. Then $m\x + n\y$ is defined. 

In Section \ref{sec:skein-punctured-torus} and \ref{subsec: pentagon in finite rank} below, we offer
evidence for a conjectural pentagon relation on the
punctured torus.  Before that, we consider the unpunctured
case.

\subsubsection{Pentagon Relation for the Torus}
\label{sec:skein-torus}
Recall that in \cite{MS17}, the authors established an equivalence between $\Sk(T^2)$ and the specialized elliptic Hall algebra $\cE_{s, s}$. Hence the pentagon relation for $\Sk(T^2)$ can be obtained from \cite{Ze},  which proved the pentagon relation for the ellipical Hall (or DIM) algebra, extending work of Garsia-Mellit \cite{GM}.
Below we give an alternate, self-contained proof.


To begin, we recall from \cite{MS17} that the elements $P_\z$, $\z\in H_1(T^2)\cong \bZ^2,$ generate $Sk(T^2)$ as an algebra and obey
\begin{equation}
\label{eq:ms-skein-equation}
    [P_{\bf a},P_{\bf b}] = \{\omega({\bf a},{\bf b})\}\, P_{\bf a+b},
\end{equation}
where $\{k\}:= s^k-s^{-k},$ and $s = q^\frac{1}{2}.$

Let us define $\cA$ to be the abstract algebra generated by
the $P_{\x}$ with the relations of Equation \eqref{eq:ms-skein-equation}, so $\cA \cong Sk(T^2).$ We take the completion $\widehat{\cA}$ as in Example \ref{eg:completion}. 


\begin{proposition}
\label{prop:adjoint-pentagon}
Suppose $\x$ and $\y$ in $\bZ^2$ are primitive vectors, with $\omega(\x,\y)=1$. 
Then as operators on $\widehat{\cA}$, we have 
\begin{equation}
\label{eq:skein-pentagon}
\mathrm{Ad}_{Q_\x (v)}\cdot \mathrm{Ad}_{Q_\y(w) } = \mathrm{Ad}_{Q_\y (w)}\cdot \mathrm{Ad}_{Q_{\x+\y} (wv)}\cdot \mathrm{Ad}_{Q_\x (v)},
\end{equation}
where the Baxter operators $Q_\x$, $Q_\y$ are as in Definition \ref{def: Baxter operator}. 
\end{proposition}


\begin{proof}

First, recall from Lemma \ref{lem:Ad-action} that
$\Ad_{Q_{\bf a}}(v) \skein{\bf b} = \sum_{n\ge 0} \qbinom{d}{n} v^n \skein{n{\bf a}+{\bf b}},$
where $d = \omega({\bf a},{\bf b}),$
and in particular
\begin{align*}
    \Ad_{Q_{\bf a} (v) } \skein{\bf b} 
    &= \skein{\bf b} + v \skein{{\bf a} + {\bf b}} ,\\
    \Ad_{Q_{\bf a}^{-1} (v) } \skein{\bf b} 
    &= \sum_{ n\ge 0} (-1)^n v^n P_{n {\bf a} + \bf b}
    \quad \text{if} \ d = 1;\\
   \Ad_{Q_{\bf a} (v) } \skein{\bf b} 
    &=  \sum_{n \ge 0} (-1)^n v^n P_{n{\bf a} + {\bf b}}, \\
    \Ad_{Q_{\bf a} (v)^{-1} } \skein{\bf b} 
    &= \skein{\bf b} + v \skein{{\bf a} + {\bf b}} \quad \text{if} \ d = -1.
\end{align*}

Since conjugation is an algebra isomorphism
of $\cA$,
it suffices to prove it for the generating set
$P_\x, P_{-\x}, P_\y$ and $P_{-\y}$. 

The cases $P_{-\x}$ and $P_{\y}$ are immediate from the equations above, so it remains to treat those of $P_\x$ and $P_{-\y}$.

\underline{Case $P_\x$}: we compute
\begin{align*}
	\Ad_{Q_\x(v)} \Ad_{Q_\y(w)} \cdot P_\x 
	& = \sum_{k\geq 0} (-1)^k \Ad_{Q_\x(v)} w^k P_{\x + k\y } \\
	& = \sum_{k\geq 0, l\geq 0} (-1)^k \qbinom{k}{l}   v^l w^k P_{l\x + k \y + \x} ,
\end{align*}
while
\begin{align*}
	\Ad_{Q_\y(w)} \Ad_{Q_{\x + \y} (vw) } \Ad_{Q_\x(v)} \cdot P_\x  
	& = \Ad_{Q_\y(w)} \Ad_{Q_{\x+\y}(vw)} P_\x\\
	& = \Ad_{Q_\y(w)} \sum_{k\geq 0} (vw)^k P_{k\x  + k\y +\x} \\
	& = \sum_{k, l\geq 0} (-1)^k \qbinom{-(k+1)}{l} v^k w^{k+l} P_{k\x + (k+l) \y +\x}
\end{align*}
Compare the coefficients of the term $ v^M w^N P_{M\x + N\y + \x} $, we need to show that
\begin{equation*}
\qbinom{N}{M} = (-1)^{N-M} \qbinom{-(N+1)}{N-M}, \ \ \ \text{for all}\  M, N \geq 0, 
\end{equation*}
which can be rewritten as
\begin{equation}
\label{eq:qbinomeq}
 (-1)^k	\qbinom{-(l+1)}{k} = \qbinom{k+l}{l} 
\end{equation}

Denote the generating function of the quantum binomials by
$$
 f^{(N)}(x) \coloneqq \sum_{k\ge 0} \qbinom{N}{k} x^k. 
$$
Then \eqref{eq:qbinomeq} is equivalent to 
\begin{align}
    f^{(-(l+1))} &= 
 \frac{1}{(1- s^l x ) (1 - s^{(l-2)} x) \cdots (1 - s^{-l}x) } \notag \\
 &= \sum_{k} \qbinom{k+l}{k} x^k  \label{generating function}
\end{align}
For $l = 0$ this is obvious; we prove We prove the general case by induction on $l$.  Assuming case the identity for $l$ consider it for $l+1$: 
\begin{align*}
	f^{(-(l+2))}(-x) 
	& = \frac{1}{1 - s^{-(l+1)}x}f^{(l+1)} (sx)\\
	&   \stackrel{\text{by induction hypothesis}}{=} (1 + s^{-(l+1) }x   + \cdots) \sum_k \qbinom{k+l}{k} s^k x^k \\
	& = \sum_N x^N \sum_{k=0}^{N} \qbinom{k+l}{k} s^{-(l+1)(N-k) + k}. 
\end{align*}
By repeatedly using the relation 
$$
\qbinom{d+1}{N} = s^N \qbinom{d}{N} + s^{N-d-1} \qbinom{d}{N-1},  
$$
we get that
$$
\qbinom{N+l+1}{N} = \sum_{k=0}^{N} \qbinom{k+l}{k} s^{-(l+1)(N-k) + k}
$$
Thus \eqref{generating function} is also true for $l+1$.

\underline{Case $P_{-\y}$}: This time, we have
\begin{align*}
	\Ad_{Q_\x(v)} \Ad_{Q_\y(w)} \cdot P_{-\y} 
	&= \Ad_{Q_\x(v)} \cdot P_{-\y}\\
	&= \sum_{k\geq 0} (-1)^k v^k P_{k \x -y}
\end{align*}
while
\begin{align*}
	\Ad_{Q_\y(w)} \Ad_{Q_{\x + \y }(vw) } \Ad_{Q_\x(v)} P_{-\y} 
	&= \Ad_{Q_\y(w)} \Ad_{Q_{\x + \y } (vw) } \sum_{k\geq 0} (-1)^k  v^k P_{k\x - \y} \\
	&= \Ad_{Q_\y(w)} \sum_{k, l\geq 0} (-1)^k \qbinom{-k-1}{l} v^{k+l} w^l P_{(k+l) \x + (l-1) \y }\\
	&= \sum_{k, l, t \geq 0} (-1)^k \qbinom{-1-k}{l} \qbinom{-(k+l)}{t}  v^{k+l} w^{l+t} P_{(k+l)\x + (l+t)\y -\y}
\end{align*}

Comparing the $ v^M w^NP_{M\x + N\y -\y} $ term, what must be shown is that:
\begin{equation}
	\sum_{0\geq l \geq M } (-1)^{M-l} \qbinom{(-M+l-1)}{l} \cdot \qbinom{-M}{N-l}  = 0, \ \ \ \text{for}\  M\geq 0, N> 0.
\end{equation}
Using \eqref{eq:qbinomeq}, note that $\qbinom{-M+l -1}{l} = (-1)^{l} \qbinom{M}{l}.$
Thus the above is equivalent to 
$$
 \sum_{0\leq l \leq M} \qbinom{M}{l}\qbinom{-M}{N-l} = 0, \ \ N>0.
$$
But this is immediate from the fact that $f^{(M)} \cdot f^{(-M)} = 1$.

%
%
%
%
%

\end{proof}

\begin{remark}
\label{rmk:pentagonrelation}
The pentagon relation for the adjoint action by the Baxter operators ensures consistent
scattering for the diagram of Example \ref{ex:basicdiagram}. 
\end{remark}

We can improve Proposition \ref{prop:adjoint-pentagon} after proving an intermediate result.
Note from the defining relation \eqref{eq:ms-skein-equation} that $\cA$ has a $\bZ^2$ grading.
We write $\cA_m$ for the $m$-graded piece, and likewise for any subset $C\subset \cA$ write $C_m \coloneqq C\cap \cA_m$.
Writing $Z(\cA)$ for the center of an algebra, $\cA$, we have the following proposition.

\begin{proposition}
\label{prop:trivial-center}
    The algebra $\cA$ has trivial center: $Z(\cA) = R.$
\end{proposition}

\begin{proof}


By \cite[Theorem 2]{MS17}, $\cA$ is the universal enveloping algebra of the Lie algebra generated by $P_\x, \ \x\in \bZ^2$. Thus we can use Poincar\'e-Birkhoff-Witt (PBW) theorem to construct an additive basis for $\cA$.

We give an order on $\bZ^2$. View $v\in \bZ^2$ as a vector in $\bR^2$ with slope in $ [0, 2\pi)$. We say $v_1 < v_2$ if the slope of $v_1$ is smaller than the slope of $v_2$. If they have the same slope, then $v_1 < v_2$ if the length of $v_1$ is smaller than the length of $v_2$. 

By the PBW theorem, 
$\cA$ has an additive basis consisting of increasing paths $\gamma$. Specifically, an increasing path $\gamma$ is an increasing sequence
    $(v_1, v_2, ... , v_n)$. 
For an increasing path $\gamma$ as above, write $P_\gamma := P_{v_1}\cdots P_{v_n}.$ 
 We call $n$ the length of such a path, and say that $\gamma$ is from $0$ to $m$, if $\sum_i v_i = m$. 
In particular, $\cA_m$ has an additive basis consisting of increasing paths $\gamma$ from $0$ to $m$.

It suffices to prove triviality of the $m$-graded piece $Z(\cA)_m$. Filter $\cA$ by the length of $\gamma$. For a contradiction, assume $c$ living in $gr(Z(A)_m)_n$ is not $0$. 

Let $c = \sum_i c_i P_{\gamma_i}$, 
and set $\gamma_i = (v_1^{(i)}, v_2^{(i)}, \cdots , v_n^{(i)} ) $. 
We can order the paths $\gamma_i$ by their first vectors. Namely, we can assume $v_1^{(1)} \le v_1^{(2)} \le v_1^{(3)} \le \cdots$, and if $v_1^{(1)} = v_1^{(2)}$, we ask $v_2^{(1)} \le v_2^{(2)}$, and so on. In particular $\gamma_1$ has the smallest first vector. 
By symmetry,  without loss of generality we can assume the slope of $v_1^{(1)}$ is in $(0, \pi/2]$ (if not, we can just rotate everything by a multiple of $\frac{\pi}{2}$). 
 
    Then let 
    $$
    K =  \sum_{\genfrac{}{}{0pt}{}{i\;\text{such that}}{v^{(i)}_1 = v^{(1)}_1}} c_i P_{\gamma_i} = P_{v_1^{(1)}}  K' 
    $$
    so $c = K + \text{terms with first vector greater than $v^{(1)}_1$}.$

    Since $c$ lives in the center we have (everything lives in $gr(\cA_m)_n $)
    \begin{align}
    0 &=    [P_{(1, 0)} , c] \notag \\ 
     &= \lambda P_{v_1^{(1)} + (1, 0)}  K' + P_{v_1^{(1)}} [P_{(1, 0)}, K'] + \cdots \label{eq: sum}
     \end{align}
    where $\lambda = \{\omega ( (1, 0), v_1^{(1)} ) \} $ is a nonzero constant. But note that by assumption, $v_1^{(1)}$ has slope in $(0, \pi/2]$.
    Thus now $v_1^{(1)} + (1, 0)$ becomes the smallest vector among all vectors appearing in the sum \eqref{eq: sum}. Since $\lambda P_{v_1^{(1)} + (1, 0)}K'$ contains all terms starting with $P_{v_1^{(1)} + (1, 0)}$, it has to be $0$. We get a contradiction. 
\end{proof}

Putting together Propositions \ref{prop:adjoint-pentagon} and \ref{prop:trivial-center}, we have the following,
due originally to Zenkevich \cite{Ze}, who extended the proof of Garsia-Mellit \cite{GM} by using a result of
Schiffman-Vasserot \cite[Corollary 1.5]{SV11}.

\begin{theorem}
\label{thm:baxter-pentagon}
    Let $\x$ and $\y$ as in Proposition \ref{prop:adjoint-pentagon}.  Then
    $$Q_\y(w) \cdot Q_\x(v) = Q_\x(v) \cdot Q_{\x+\y}(vw)\cdot Q_\y(w)$$
    in $\widehat{\cA} \cong \widehat{\Sk}(T^2).$
\end{theorem}
\begin{proof}
    The proof in $\widehat{\cA}$ is then immediate from Propositions \ref{prop:adjoint-pentagon} and \ref{prop:trivial-center},
    since in an algebra over $R$ with trivial center, if $Ad_a = Ad_b$ then $a=\lambda b,$ for some $\lambda \in R$. 
    The explicit form of the Baxter operators then fixes $\lambda = 1.$
\end{proof}

\subsubsection{The pentagon relation for the punctured torus}
\label{sec:skein-punctured-torus}

Recall that in skein algebra of the punctured torus $\Sk(T^2 - D)$ we still have a well-defined element $P_\x$ associated to each primitive vector $\x\in  \bZ^2$,  and a well-defined Baxter operator $Q_\x(v)$. 
\begin{conjecture}\label{conj: pentagon relation for punctured torus}
    For $\x , \y \in \bZ^2$, $\omega(\x, \y) =1$, the pentagon relation:
    \begin{equation}\label{eq: pentagon for punctured}
        Q_\x(v) \cdot Q_\y(w) = Q_\y(w) \cdot Q_{\x+\y} (vw) \cdot Q_\x(v)
    \end{equation}
    holds in $\widehat{\Sk}(T^2 - D)$. 
\end{conjecture}

On a general surface $\Sigma$, if two simple closed curves $l_1$ and $l_2$ intersect transversely at a single point, their neighborhood is homeomorphic to a punctured torus. Thus we have:
\begin{corollary}[of the conjecture]
    Let $\Sigma$, $l_1$ and $l_2$ be as above. Then the pentagon relation
    \begin{equation*}
        Q_{l_1} (v) \cdot Q_{l_2} (w) = Q_{l_2} (w) \cdot Q_{l_1 + l_2} (wv)\cdot  Q_{l_1} (v)
    \end{equation*}
    holds in $\widehat{\Sk}(\Sigma)$. 
\end{corollary}

Recall that
\cite{MS21} constructed a surjection from
    $\Sk(T^2 - D)^+ $ to (the positive part of) the elliptic Hall algebra $\cE_{\sigma, \bar{\sigma}}^+$, assuming a conjecture. Thus, Conjecture \ref{conj: pentagon relation for punctured torus} generalizes the pentagon relation for the elliptic hall algebra in \cite{Ze}, and hence the pentagon relation in \cite{GM}. In particular, it implies Theorem \ref{thm:baxter-pentagon}. 

    As is pointed in \cite{MS17}, $\Sk(T^2 - D)$ is indeed much larger than the elliptic Hall algebra, in the sense that 
    as a free R-module, the latter is isomorphic to a polynomial algebra with generators indexed by $\bZ^2$, while the former is isomorphic to a polynomial algebra with generators given by all conjugacy classes in $\pi_1(T^2 -D)$ \cite{P1}.

    We currently do not know how to prove Conjecture~\ref{conj: pentagon relation for punctured torus} for the full HOMFLYPT skein, although in Section~\ref{sec:higher-rank} we show that it holds in the reduction to any finite rank skein.
    
    However, in the full HOMFLYPT skein we can do some checks in low degrees, as we shall now describe. 

    Let's consider the coefficient of the term $v^m w^n$ in  \eqref{eq: pentagon for punctured}. When $n = 1$, the two sides matching up is equivalent to:
    $$
    \Ad_{Q_\x} P_\y = P_y + P_{\x + \y}. 
    $$
    This is true since if we go back to proof of Lemma \ref{lem:Ad-action}, in this special case we only need to use the relation $[P_{m\x}, P_\y] = \{m \} P_{m\x + \y} $, which holds in the skein algebra of punctured torus. 

    When $m, n \ge 2$, things become more complicated. Denote $u_{d \x_0} \coloneqq \frac{1}{ \{d\} } P_{\x_0}$ for a primitive $\x_0 \in \bZ^2$. 
    If we compute the $v^2w^2$ term, the relation we need is:
    $$
    [u_{2\x}, u_{2\y}] = [u_{\x}, u_{\x + 2\y}] - 2 u_{2\x + 2\y},
    $$
    which is already nonobvious.  We don't know any reason for it to hold a priori, while one can verify it by a direct but tedious skein diagram computation. 
    
    For the $v^m w^2$ ($m\ge 3$) terms, it turns out that the needed relation is 
    \begin{equation*}
        [u_{m\x} , u_{2\y}] - [u_{(m-1)\x}, u_{\x + 2\y}] + [u_{(m-2)x}, u_{2\x +2\y}] = 0.
    \end{equation*}

\section{Skein-Valued Wavefunctions}
\label{sec:wavefunctions}

In this section, we generalize the wavefunction constructions of \cite{ES19,SSZ} to include skeins and framings, 
and compute some examples.

\subsection{The Wavefunction}
\label{subsec:wavefunction}

Let us recall the skein-valued additive face relations of Equation \eqref{eq:skein-face-relation}.
Fix a cubic graph $\Gamma$ and zero-chains $\mathfrak{b}_\pm$ and $\mathfrak{p}_\pm$ on
the surface $S_\Gamma$.  As explained in Section \ref{sec:skeins with b and f},
these zero-chains define base points and framing points for the relative skein
$\Sk (\Lambda_\Gamma, \mathfrak{b}_\pm; \mathfrak{p}_\pm).$
Now fix a face $f$ of the graph with vertices $v_i$.
This data defines a face relation 
$$	A_{\Gamma, f} = \frame{l_1} + \frame{l_2} + \cdots + \frame{l_n} \in \Sk (S_\Gamma, \mathfrak{b}_\pm; \mathfrak{p}_\pm),$$
where $l_i$ is a path from $\fb_-$ to $\fb_+$ (or more precisely the parts of those 0-chains within $f$)
passing through the vertex $v_i$ (see Figure \eqref{fig: face relation}).

Given a geometric framed seed $ \underline{\mathbf{i}} = (\Gamma, \cT, \mathbf{f}, \gamma)$,  the framing path $\gamma$ induces 
an isomorphism of $R$-modules $\Sk(S_\Gamma, \fb_\pm ; \fp_\pm) \xrightarrow{\gamma} \Sk(S_\Gamma, \fb_\pm) $ \eqref{iso of skeins}.  Hence $A_{\Gamma, f}$ acts on $\widehat{\Sk}(L_\cT)$ by \eqref{action}. 

Take the ``necklace graph'' as in Figure \ref{subfig: neck}, with the framing and framing paths shown in the picture. According to \cite[Example 5.1]{SSZ}, it has an exact filling, which defines a geometric framed seed $\seed{i}_{\mathrm{neck}}$.

In this section, we consider all seeds ``connected'' to $\seed{i}_{\mathrm{neck}}$. To be more precise, let us use  $\mathbb{G}_{\mathrm{ad}}(\seed{i}_{\mathrm{neck}})$ to denote the connected component of the groupoid $\mathbb{G}_{\mathrm{ad}}$ (defined in Section \ref{sec:skein-cluster-groupoid}) containing $\seed{i}_{\mathrm{neck}}$. 
For each $\seed{i} \in \mathbb{G}_{\mathrm{ad}}(\seed{i}_{\mathrm{neck}})$, we define a wavefunction $\Psi_{\seed{i}}$, which will be the unique solution to the face relations:
$$
A_{\Gamma, f} \cdot \Psi_{\underline{\mathbf{i}}} = 0. 
$$

We start from the necklace seed $\underline{\mathbf{i}}_\mathrm{neck}$. Since its filling is exact, we expect $\Psi_\mathrm{neck} = 1 \in \widehat{\Sk}(\cH_g)$ . 
We will show in Corollary \ref{cor:uniq sol} that this is in fact the unique solution to the face relations.
Again as in \cite{SSZ}, consider an admissible mutation:
\begin{align*}
     \seed{i} \xlongrightarrow{k} \seed{i}' 
\end{align*}
given by mutating the edge $k$ with sign $\epsilon$. It induces an isomorphism: 
\begin{align*}
    \nu_k^\epsilon: \Sk(S_{\Gamma'}) \longrightarrow \Sk(S_{\Gamma}). 
\end{align*}
The framing $\mathbf{f}$ gives us an isomorphism:
$\Sk(S_\Gamma) \xrightarrow{\mathbf{f}} \Sk(\Sigma_g)$. Hence a mutation of framed seed gives rise to an automorphism $\mu_k^\cD \coloneqq \mathbf{f} \circ \mu_k \circ \mathbf{f}' $ of $\widehat{\Sk}(\Sigma_g)$, which factors as:
\begin{align*}
    \mu_k^\cD = \Ad_{Q^\epsilon_{\epsilon \mathbf{f}(k)}} \circ (\mathbf{f} \circ \nu_k \circ \mathbf{f}'^{-1} )  
\end{align*}
We can also characterize $\mu_k^\cD$ by the commutative diagram:
\begin{align*}
\xymatrix{
\Sk(S_{\Gamma'}) \ar[r]^-{\nu^\epsilon}  \ar[d]^-{\mathbf{f}} 
&\Sk(S_\Gamma) \ar[r]^-{\Ad_{Q^\epsilon_{\epsilon \mathbf{f}(k)}}}  
&\Sk(S_\Gamma) \ar[d]^-{\mathbf{f}'}  \\
\Sk(\cH_g) \ar[rr]^-{\mu_k^\cD} 
& &\Sk(\cH_g)
}
\end{align*}
We define $\Phi_k \coloneqq Q^\epsilon_{\epsilon \mathbf{f}(k)}$. If we have a path $\Vec{a}: \underline{\mathbf{i}} \rightarrow \underline{\mathbf{i}}_0 $ in $\mathbb{G}_\mathrm{ad}$, i.e. a sequence of admissible mutations
$$\Vec{a} :  \seed{i} = \seed{i}_0 \xrightarrow{k_1} \seed{i}_1 \xrightarrow{k_2} \cdots \xrightarrow{k_n} \seed{i}_n = \seed{i}_n =\seed{i}' 
$$
 then we define  
 $$\Phi_{\Vec{a}} \coloneqq \Phi_{k_n } \circ \cdots  \circ \Phi_{k_1} \in \widehat{\Sk}(\cH_g)
 $$
According to Theorem \ref{thm:face relation}, it is easy to see that if $\Psi'$ satisfies all the face relations $A_{\Gamma', f'}$ for $\seed{i}'$, 
then $\Phi_{\Vec{a}}\cdot \Psi'$ would satisfy the face relations $A_{\Gamma, f}$ for $\seed{i} $. Therefore, if $\Vec{a} $ is path from seed $\seed{i}$ to the necklace seed $\seed{i}_\mathrm{neck}$, we obtain the following candidate for the wave function of $\seed{i}$:
\begin{equation*}
    \Psi_{\underline{\mathbf{i}}} \coloneqq \Phi_{\Vec{a}} \cdot \Psi_\mathrm{neck}
\end{equation*}

The main goal of this section is to prove the following.

\begin{theorem}\label{thm:wavefunction}
    The above map:
\begin{equation*}
    \Psi : \ \mathrm{Ob} (\mathbb{G}_\mathrm{ad} (\underline{\mathbf{i}}_\mathrm{neck})) \longrightarrow \widehat{\Sk} (\cH_g), \quad \underline{\mathbf{i}} \longmapsto \Psi_{\underline{\mathbf{i}}}
\end{equation*}
    is well-defined, i.e.~is independent of the choice of path $\Vec{a}: \underline{\mathbf{i}} \rightarrow\underline{\mathbf{i}}_\mathrm{neck}$. 
    Moreover, each $\Psi_{\underline{\mathbf{i}}}$ is the unique solution to the face relations of the corresponding graph. 
\end{theorem}

We first show the uniqueness of the solution for the necklace graph. In order to do this we need the following lemma:
\begin{lemma}\label{lem: sliding relation}
    Let $M$ be a $3$-manifold with boundary $\Sigma$, and $\gamma$ a loop on $\Sigma$ bounding a disk $D$ in $M$. If an element  $\Psi \in {\Sk}(M)$ satisfies the ``sliding relation'':
    \begin{equation}\label{eq: sliding relation}
        (\skein{\gamma} - \bigcirc) \Psi = 0,
    \end{equation}
    then  $\Psi$ lives in the image of ${\Sk} (M \backslash D)\rightarrow {\Sk}(M).$ 
\end{lemma}
\begin{proof}
        Without loss of generality, we can work over the field $\mathrm{Frac} (R)$. 
        By $\cL (M)$ we denote the free module generated by the isotopy class of framed links in $M$. Then $\Sk(M) = \mathcal{L}(M ) / S$, where $S$ is the submodule generated by HOMFLY skein relations. 
		
		We define a filtration $F_\bullet$ on $\mathcal{L} (M)$ whose $n$-th piece $F_n \mathcal{L} (M)$ is spanned by links that can be isotopied to transversely intersect $D$ at no more than $n$ points. This induces a filtration on $\Sk (M)$. We want to show that  $\Psi \in F_0 \Sk(M)$.
		
		Assume that  $\Psi \in F_n\Sk(M)$ with $n \ge 1  $ satisfies the sliding relation \eqref{eq: sliding relation}. 
        We will show it also lives in $F_{n-1}\Sk(M)$. 
		
		In addition to the filtration defined above, the module $\mathcal{L} ({M})$ has grading given by the intersection number with $D$. Since $R$ is a homogeneous submodule, this also gives a grading on $\Sk(M)$. The ``sliding operator'' $(\skein{\gamma} - \bigcirc)$ preserves this grading. Thus each homogeneous component of  $\Psi$ is annihilated by $\skein{\gamma} - \bigcirc$. So we can assume $\Psi$ is homogeneous, i.e.~it is a sum of links having the same intersection number with $D$. 
  
	Consider the following map:
		\begin{equation} \label{eq:glue}
			\Sk_{a, n-a}(D\times I) \otimes \Sk_{n-a, a}(M \backslash D) \longrightarrow F_n\Sk(M).
		\end{equation}		
        Here $\Sk_{a, n-a}$ means the skein module with $n$ base points on $D\times \{0,1\}$, $a$ oriented positively and $n-a$ oriented negatively --- 
        see Figure \eqref{fig: decomposition}. The map is the obvious one, i.e.~if we have two tangles in the two parts we can just glue them together to get a link in the handlebody.

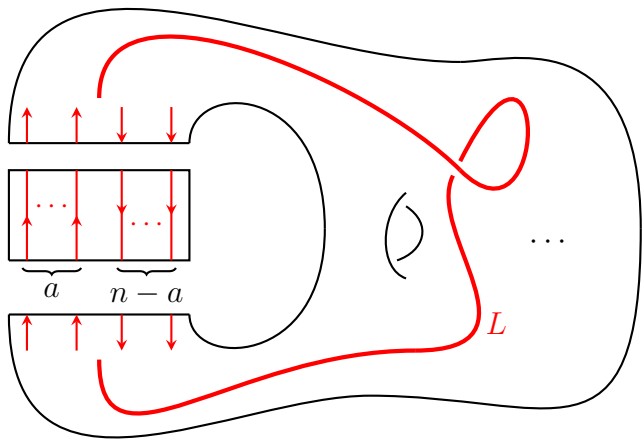
\begin{figure}[htpb]
	\centering
    \vspace{-1cm}
	\begin{tikzpicture}[scale=1.2]
	\draw[thick] (0, 0) -- (2, 0) -- (2, 1) -- (0, 1) -- (0,0);
	\draw[thick] (0, 1.3) -- (2, 1.3);
	\draw[thick] (0, -.6) -- (2, -.6); 
	
	\draw[thick, red, -stealth]   (.2, 0) -- (.2 , .5);
	\draw[thick, red] (.2, .4) -- (.2, 1);
	
	\draw[thick, red, -stealth]   (.75, 0) -- (.75 , .5);
	\draw[thick, red] (.75, .4) -- (.75, 1);
	
	\node at (.5, .6) {\color{red} \small$\cdots$ };
	
	\draw[thick, red, -stealth]   (1.8 , 1) -- (1.8 , .5);
	\draw[thick, red] (1.8, .6) -- (1.8, 0);
	
	\draw[thick, red, -stealth]   (1.25 , 1) -- (1.25 , .5);
	\draw[thick, red] (1.25, .6) -- (1.25, 0);
	
	\node at (1.55, .4) {\color{red} \small $\cdots$};
	
	\draw[thick, red, -stealth]   (.2, 1.3) -- (.2 , 1.7);
	\draw[thick, red, -stealth]   (.75, 1.3) -- (.75 , 1.7);
	\draw[thick, red, -stealth]   (1.8 , 1.7) -- (1.8 , 1.3);
	\draw[thick, red, -stealth]   (1.25 , 1.7) -- (1.25 , 1.3);
	
	\draw[thick, red, -stealth]   (.2, -1) -- (.2 , -.6);
	\draw[thick, red, -stealth]   (.75, -1) -- (.75 , -.6);
	\draw[thick, red, -stealth]   (1.8 , -.6) -- (1.8 , -1);
	\draw[thick, red, -stealth]   (1.25 , -.6) -- (1.25 , -1);
	
\draw [
thick,
decoration={
	brace,
	mirror,
	raise=0.1cm
},
decorate
] (.15, 0) -- (.8, 0) node [pos=0.5,anchor=north,yshift=-0.15cm] {$a$}; 

\draw [
thick,
decoration={
	brace,
	mirror,
	raise=0.1cm
},
decorate
] (1.2, 0) -- (1.85, 0) node [pos=0.5,anchor=north,yshift=-0.15cm] {$n-a$};

\draw[thick] (2, -.6) .. controls (2, -1.2) and (3.5, -1.2) .. (3.5, .35);
\draw[thick] (2, 1.3) .. controls (2, 2) and (3.5, 2) .. (3.5, .35);

\draw[thick] (0, 1.3) .. controls (0, 4) and (3, 2.2).. (5,2.2);

\draw[thick] (5, 2.2) ..controls (5.5, 2.2) and (7, 2.7)..  (7, .35);

\draw[thick] (0, -.6) .. controls (0, -3) and (2.5, -1.5).. (4, -1.5);

\draw[thick] (4, -1.5) .. controls (6, -1.5)and (7, -2.4) .. (7, .35); 

\begin{scope}[xshift= 3cm, yshift=3.3cm]
	\draw[thick] (1.4, -2.55) .. controls (1.1, -2.8) and (1.1, -3.4).. (1.4, -3.5);
	\draw[thick] (1.4, -2.7) .. controls (1.7, -2.9) and (1.6, -3.2) .. (1.3, -3.3);
\end{scope}

\node at (6, 0.2) {$\cdots$};

\draw[red, ultra thick] (1, 1.8) .. controls (1, 3.3) and (4, 2)..  (5, 1);
\draw[red, ultra thick] (5, 1) .. controls (6, 0) and (6, 3.1).. (5, 1.1); 
\draw[red, ultra thick] (4.92, .94) ..controls (4.6, .3) and (6, -1) .. (4.5, -1);
\draw[red, ultra thick] (4.5, -1) .. controls (2.5, -1) and (1, -2.5) .. (1, -1.1);

\node at (5.4, -.7) {\color{red} $L$};

	\end{tikzpicture}
	\vspace{-1cm}
	\caption{Decomposition of a skein $L \in \Sk(M)$ into factors in $\Sk_{a,n-a}(D\times I)\otimes \Sk_{n-a,a}(M\setminus D)$}
	\label{fig: decomposition}

\end{figure}

		There is an induced map:
		\begin{equation}
		 G:\ \ 	F_n/F_{n-1} \Sk_{a, n-a}(D \times I) \otimes \Sk_{n-a , a} (M \backslash D) \longrightarrow F_n/F_{n-1} \Sk (M)
		\end{equation}
		We now view $\Psi$ as an element in $F_n/F_{n-1} \Sk (M)$. Since each summand of $\Psi$ has the same intersection number with $D$, we can assume they all intersect $D$ at $a$ positive points and $n-a$ negative points. Thus $\Psi$ is the image of a decomposable element $T \otimes\Psi'$ under $G$, where $T$ is the trivial tangle in $\Sk_{a, n-a} (D\times I)$.   
		
		We need to show that  $\Psi = 0$ in $ F_n/F_{n-1} \Sk (M)$, which is equivalent to $T \in \ker (G(- \otimes\Psi'))$. Since $\Psi = G(T\otimes \Phi')$
        satisfies the sliding relation, i.e.
    \begin{align*}
        (P_\gamma - \bigcirc) G(T\otimes \Psi') &= G( ((P_\gamma - \bigcirc)T) \otimes \Psi') \\
        &= 0
    \end{align*}
        we know that $(\skein{\gamma} - \bigcirc) T$ lives in $\ker (G (- \otimes\Psi'))$. Because $\skein{\gamma} - \bigcirc$ preserves the grading $F_\bullet$, it also preserves the submodule $\ker (G (- \otimes\Psi'))$.  In \cite[Section 2]{gilmer2001homflypt}, it is shown that the operator $\skein{\gamma} - \bigcirc$ acting on  $F_n/F_{n-1} \Sk_{a, n-a}(D \times I)$ has an eigenbasis with nonzero eigenvalues. In particular, it is invertible (since we are working over $\mathrm{Frac}(R)$). Combining these facts we get that $T$ also lives in $\ker (G (- \otimes\Psi'))$, whence $\Psi\in F_{n-1} \Sk (M)$ as claimed. This proves the Lemma.
\end{proof}

\begin{corollary}\label{cor:uniq sol}
    The only solutions to the face relations of the necklace graph are the scalars. 
\end{corollary}
\begin{proof}
First, it is not hard to see that scalars are  solutions. We only need to prove the uniqueness. 
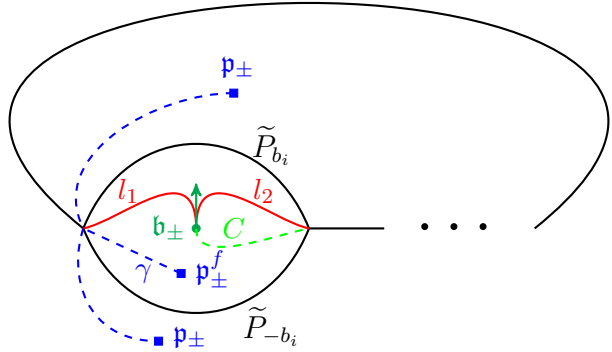
\begin{figure}[htpb] 
\centering
\vspace{-1cm}
\begin{tikzpicture}[scale= 1]
		
    \draw[thick] (0, 0) .. controls (0.6, 1.5) and (2.4, 1.5).. (3, 0);
    \draw[thick] (0, 0) .. controls (0.6, -1.5) and (2.4, -1.5).. (3, 0);
    \draw[thick] (4, 0) -- (3, 0) ;
    \draw[thick] (0, 0) .. controls (-5, 4)  and (11, 4) .. (6, 0);
    \node[thick] at (5, 0) {\huge$\cdots$};
    \filldraw[color=Green] (1.5, 0) circle (1.5pt) node[anchor = east]{\small$\mathfrak{b}_\pm$} ;

    \draw[thick, red] (1.5, 0) .. controls (1.5, 1) and (2.5, 0.1) .. (3, 0);
    \draw[thick, red] (1.5, 0) .. controls (1.5, 1) and (0.5, 0.1) .. (0, 0);

    \draw[thick, green, dashed] (1.5, 0) .. controls (1.5, -0.5) and (2.5, -0.1) .. (3, 0); 
    \node at (2, -0) {\small \color{green} $C$};
 				
    \filldraw[blue] (1.3, -0.6) +(-1.5pt, -1.5pt) rectangle +(1.5pt, 1.5pt) node[anchor = west]{\small$\fp_\pm^f$}; 
    
    \draw[thick, blue, dashed] (1.3, -0.6) -- (0, 0); 
    \node at (.8, -0.6){\small \color{blue} $\gamma$};

    \node at (0.6, 0.5) {\small \color{red} $l_1$};
    \node at (2.4, 0.5) {\small \color{red} $l_2$};
    
    \draw[thick, Green, ->, > = stealth'] (1.5, 0) -- (1.5, 0.6);

    \node at (2.5, 1.1) {$\frame{b_i}$}; 
    \node at (2.5, -1.3){$\frame{-b_i}$}; 

    \filldraw[blue] (2, 1.8) +(-1.5pt, -1.5pt) rectangle + (1.5pt, 1.5pt) node[anchor = south]{\small$\fp_\pm$};  
    \draw[dashed, blue, thick] (2, 1.8) .. controls (1, 1.8) and (-.5, 1) .. (0,0); 

    \filldraw[blue](1, -1.5)  +(-1.5pt, -1.5pt) rectangle +(1.5pt, 1.5pt) node[anchor = west]{\small$\fp_\pm$};
    \draw[dashed, blue, thick] (0, 0) ..controls (-.2, -.6) and (0, -1.5) .. (1, -1.5);
    
    \end{tikzpicture}
    \caption{One bead of a necklace}\label{fig: necklace}
\end{figure}

    After choosing a capping path $C$ and a framing path $\gamma$, as in Figure \ref{fig: necklace}, the face relations become:
    \begin{equation*}
        (\skein{b_i} - \bigcirc ) \Psi_{\mathrm{neck}}= 0. 
    \end{equation*}
    Each $b_i$ bounds a disk $D_i$. By lemma \ref{lem: sliding relation} we can assume $\Psi_\mathrm{neck}$ lives in $ \cH_g \backslash (\cup D_i) $, which is contractible. Hence it has to be a scalar. 
\end{proof}

\begin{proof}[Proof of Theorem \ref{thm:wavefunction}]
By Theorem \ref{thm:face relation}, we know that $\Psi_{\underline{\mathbf{i}}}$ is a solution of the corresponding face relations. We only need to show that for each seed living in $\mathbb{G}_{ad} (\underline{\mathbf{i}}_\mathrm{neck})$, there is a unique solution to the face relations. But again by Theorem \ref{thm:face relation}, if we have two solutions, we can mutate back to the necklace graph and get two solutions. Thus by Corollary \ref{cor:uniq sol}, the two solutions must be the same.

\end{proof}

Recall that \cite{ES19} defined the wavefunction of a Lagrangian $L$ by counting certain open $J$-holomorphic curves $C$ with boundary in $L$:
\begin{equation}\label{eq: ES wavefunction}
    \Psi_L = \sum_{u \in \cM } \omega(u) \cdot z^{ -\chi(u)} \cdot a^{ u \link L} \cdot [\partial u] \in \widehat{\Sk} (L) ,
\end{equation}
where $\omega(u)$ is a rational number,  $\chi(u)$ is the Euler characteristic of the domain of $u$, $u\link L$ is a linking number, and $\partial u$ is considered as an element in $\Sk(L)$ --- see \cite[Definition 5.1]{ES19}. 

  In \cite[Theorem 7.1]{SS}, Scharitzer and Shende showed that this wavefunction satisfies the same face relation as in
  Equation \eqref{eq:skein-face-relation}, if all the choices there are made properly 
    \footnote{They actually showed the Ekholm-Shende wavefuntion is annihilated by \eqref{eq:skein-face-relation} after choosing capping paths. But all the arguments there should still work without choosing capping paths.}.
     We've already seen that if a seed $\seed{i}$ lives in $\mathbb{G}_{\mathrm{ad}}(\seed{i}_\mathrm{neck})$, then the corresponding face relations have a unique solution. Hence our wavefunction $\Psi_{\seed{i}}$ coincides with the wavefunction of \cite{ES19}.

As is pointed out in \cite[Remark 3.2]{SS}, their arguments should work for any skein-valued curve counting setup compatible with SFT stretching (\cite{ES19,ekholm2022ghost} show that there exists one), i.e.~that any wavefunction in such setup should satisfy the face relations. Thus it will coincide with our wavefunction defined above. 

\begin{remark}
    Here is a heuristic explanation of how wavefunctions change under mutation. If two seeds $\seed{i} = (\Gamma, L, \cdots) $ and $\seed{i}' = (\Gamma', L', \cdots)$ are related by a mutation at an edge $e$, then as explained in Section \ref{sec:groupoid}, $L'$ can be obtained from $L$ by gluing a concordance $M_{\Gamma, e}$, isomorphic to a Harvey-Lawson Lagrangian $HL_\pm$. Thus by an SFT stretching argument, morally one has:
    $$
    \{ \text{counting in } L'\} = \{\text{counting in } M_{\Gamma, e} \} \cdot \{\text{counting in } L \}  . 
    $$
    which should give rise to
    \begin{align*}
        \Psi_{L'} &= \Psi_{M_{\Gamma, e}} \Psi_L \\
        & = Q_e \cdot \Psi_L  
    \end{align*}
\end{remark}

\begin{example}[Wavefunction of the canoe] \label{eg: wavefuntion of the canoe}

\begin{figure}[htpb]
	\centering
        \begin{subfigure}{1\textwidth}
        \centering
	\begin{tikzpicture}[scale = .8]
		\draw[thick] (0, 0) .. controls (0.6, 1.5) and (2.4, 1.5).. (3, 0);
		\draw[thick] (0, 0) .. controls (0.6, -1.5) and (2.4, -1.5).. (3, 0);
		\draw[thick] (4, 0) -- (3, 0);
		
		\draw[thick] (4, 0) .. controls (4.6, 1.5) and (6.4, 1.5).. (7, 0);
		\draw[thick] (4, 0) .. controls (4.6, -1.5) and (6.4, -1.5).. (7, 0);
		\draw[thick] (7, 0) -- (8, 0);
		
		\draw[thick] (0, 0) .. controls (-5, 4)  and (16, 4) .. (11, 0);
		\node[thick] at (9.5, 0) {\huge $\cdots$};
		\filldraw[color=Green] (1.5, 0) circle (1.5pt) ;
		\filldraw[color=Green] (5.5, 0) circle (1.5pt) ;
		
		\filldraw[blue] (1.3, -.2) +(-1.5pt, -1.5pt) rectangle +(1.5pt, 1.5pt) node[anchor = north] {\small$\fp_\pm^{f_1}$};
		\filldraw[blue] (5.3, -.2) +(-1.5pt, -1.5pt) rectangle +(1.5pt, 1.5pt) node[anchor = north] {\small$\fp_\pm^{f_2}$};
		
		\draw[blue, dashed] (1.3, -.2) -- (0, 0);
		\draw[blue, dashed] (5.3, -.2) -- (4, 0);
		
		\node at (2.5, 1.2) {\small$\frame{b_1}$}; 
		\node at (2.5, -1.4){\small$\frame{-b_1}$}; 
		
		\node at (6.5, 1.2) {\small$\frame{b_2}$}; 
		\node at (6.5, -1.4){\small$\frame{-b_2}$}; 
		
		\node at (3.5, .35) {\small $\frame{-a_1}$};
		\node at (7.5, .35) {\small $\frame{-a_2}$};
    \filldraw[blue] (2, 1.8) +(-1.5pt, -1.5pt) rectangle + (1.5pt, 1.5pt) node[anchor = south]{\small$\fp_\pm$};  
    \draw[dashed, blue, thick] (2, 1.8) .. controls (1, 1.8) and (-.5, 1) .. (0,0); 
    \filldraw[blue](1.5, -2)  +(-1.5pt, -1.5pt) rectangle +(1.5pt, 1.5pt) node[anchor = north]{\small$\fp_\pm$};
    \draw[dashed, blue, thick] (0, 0) ..controls (-1, -1) and (1, -1.8) .. (1.5, -2);
    
    \node at (1, .5) {\small $f_1$};
    \node at (5, .5) {\small $f_2$};
    
  \end{tikzpicture}
  \caption{$g$-necklace framed seed} \label{subfig: neck}
  \end{subfigure}\vspace{-1cm}

    \begin{subfigure}{1\textwidth}
    \centering
	\begin{tikzpicture}[scale = .8]
    \node at (5.5, 6) { $\xymatrix{ {}\\ \ar[d] \\{} }$};
		\draw[thick] (0, 0) .. controls (-5, 5)  and (16, 5) .. (11, 0);
		\draw[thick] (0, 0) -- (2.5,1.5) -- (2.5, -1.5) -- (0, 0);
		\draw[thick] (2.5, 1.5) -- (5.5, 1.5) -- (5.5, -1.5) -- (2.5, -1.5);
		\draw[thick] (5.5, 1.5) -- (7, 1.5);
		\draw[thick] (5.5, -1.5) -- (7, -1.5);
		\node at (8, 0) {\huge$\cdots$}  ;
		
		\filldraw[Green] (1.6, 0) circle(1.5pt);
		\filldraw[Green] (4, 0) circle(1.5pt);
		
		\filldraw[blue] (1.5, -.3) +(-1.5pt, -1.5pt) rectangle +(1.5pt, 1.5pt) node[anchor = west]  {\small \color{blue} $\fp_\pm^{f_1} $};
		\draw[thick, Blue, dashed] (0, 0) -- (1.5, -.3);
		
		\filldraw[blue] (3.7, -.3) +(-1.5pt, -1.5pt) rectangle +(1.5pt, 1.5pt) node[anchor = north]  {\small\color{blue} $\fp_\pm^{f_2} $};
		\draw[thick, Blue, dashed] (2.5, 1.5) -- (3.7, -.3);
		
		\node[rotate around = {32:(0,0)}] at (.9, 1) {\small$\frame{-a_1 +b_1}$}; 
		\node[rotate around = {-30:(0,0)}] at (1, -1.1) {\small$\frame{-b_1}$}; 
		\node at (2.9, 0) {\small$\frame{a_1}$}; 
		\node at (5.1, 0) {\small$\frame{a_2}$}; 
		\node at (4, 1.9) {\small$\frame{b_2-a_2}$}; 
		\node at (4, -1.9) {\small$\frame{-b_2+a_2}$};
   \filldraw[blue] (2, 2.3) +(-1.5pt, -1.5pt) rectangle + (1.5pt, 1.5pt) node[anchor = south]{\small$\fp_\pm$};  
    \draw[dashed, blue, thick] (2, 2.3) .. controls (1, 2.5) and (-.6, 1.9) .. (0,0); 
    \filldraw[blue](2, -2.5)  +(-1.5pt, -1.5pt) rectangle +(1.5pt, 1.5pt) node[anchor = north]{\small$\fp_\pm$};
    \draw[dashed, blue, thick] (0, 0) ..controls (-1, -1) and (1, -2.5) .. (2, -2.5);
    
    \node at (2.1, .7) {\small $f_1$};
    \node at (4.3, .8) {\small $f_2$};
    
    \end{tikzpicture}
    \caption{$g$-canoe framed seed}\label{subfig: canoe}
    \end{subfigure}
    \caption{Mutation from necklace to canoe}
    \label{fig: neck to canoe}
    
\end{figure}

Consider framed seed given by the $g$-canoe graph, which can be obtained from the necklace graph by mutating the edges labelled by $\frame{-a_i}$ --- see Figure \ref{fig: neck to canoe}. Thus by definition we have:
\begin{align}\label{eq: g-canoe}
    \Psi_{g\mathrm{-canoe}} = \prod_i Q_{a_i} .
\end{align}

\end{example}

\subsection{From the Necklace to the Canoe}

As shown in Corollary \ref{cor:uniq sol} above, the face relations of the genus-$g$ necklace graph have constant solutions,
which we normalize to get $$\Psi_{\Gamma^{\rm{neck}}_g} \equiv 1.$$

We consider the case $g=1$.  Note $S_{\Gamma^{\rm{neck}}_1}$ is a torus and we take a strand to
correspond to the $(0,1)$ cycle.
The tetrahedron graph $\Gamma^{\rm{canoe}}_1$ is the mutation along a necklace strand.
The surface is the boundary of the Aganagic-Vafa (AV) brane in $\bC^3$, described by the Harvey-Lawson cone ---
see \cite[Sections 5.1.1 and 5.1.3]{SSZ}.
By Example \ref{eg: wavefuntion of the canoe}
we have
$$\Psi_{AV} = Q_{(0,1)}\cdot 1.$$
As follows from Theorem \ref{thm:face relation}, this wavefunction satisfies the face relations of the $g=1$ canoe (tetrahedron). Let's consider face $f_1$ in Figure \ref{subfig: canoe}. If we choose the capping path to go from the middle to the upper-right vertex, then the corresponding face relation becomes:
\begin{align}
\label{eq:g1-canoe}
\left(a^{-1} \bigcirc + (-a)^{-1} P_{b_1} + P_{a_1}\right) \Psi_{AV} = 0,
\end{align}
which is the same equation as in \cite[Theorem 1.1]{ES20}, with the same (unique) solution.

We must resolve the question of framings.
There is no canonical basepoint for framings, and it turns out
that we have landed in framing $-1$, in the conventions of \cite{AKMV}. 
Changes of framing are canonical, however:  they are mapping classes of the handlebody filling.
When $g=1,$ this group is $\bZ,$ and
a change of framing $p\in \bZ$ sends $(0,1)\mapsto (p,1)$.
The change of framing also affects the framing path $\gamma,$ which changes
the isomorphism \eqref{iso of skeins}.  As a result, $P_{(0,1)}\mapsto a^{-p}P_{(p,1)},$
so $Q_{(0,1)} \mapsto Q_{(p,1)}(a^{-p})$.
Therefore, the framing-$(p-1)$ wavefunction is given by
$$\Psi^{(p-1)}_{AV} = Q_{(p,1)}(a^{-p})\cdot 1.$$
In fact the full open Gromov-Witten generating function is given by the one-leg topological vertex \cite{AKMV}:
$\Psi_{AV}^{(p-1)} = \mathrm{TV}_{\bullet,\varnothing,\varnothing}^{(p-1)}$
We should then expect
$$Q_{(p,1)}(a^{-p})\cdot 1 = \mathrm{TV}_{\bullet,\varnothing,\varnothing}^{(p-1)}$$

We  prove this in the next section, first for $p=0$ and then for general $p$.

\subsection{Relating the Canoe Wavefunction
to the Topological Vertex}

Recall that $\Psi_{AV}$ lives in the skein module of the solid torus, which is isomorphic to the ring of symmetric functions. In the this section we identify $P_{0, n} \in \Sk(\cH_1)$ with $p_n(x) \in R[x_1, x_2, \cdots]^{\mathrm{Sym}} $, and $W_\lambda \in \Sk(\cH_1) $ with $s_\lambda(x) \in R[x_1, x_2, \cdots]^{\mathrm{Sym}} $.

We recall from \cite{Z1} (see also \cite{DZ}) the one-leg topological vertex in framing $p$,
in a form that will be convenient for us:
$$
\Psi^{(p)} = \sum_\lambda q^{p\kappa_\lambda {/2}} 
s_\lambda(q^\rho) s_\lambda (x) ,
$$
where $q^\rho = (q^{-1/2}, q^{-3/2}, q^{-5/2}, ... )$, $x = (x_1, x_2, ...)$ and $s_\lambda$ are the Schur polynomials.

We can rewrite this using Cauchy's identity
$$
\sum_\lambda s_\lambda(x) s_\lambda(y) = \prod_{i, j \geq 0} \frac{1}{1-x_i y_j} = \e \left(\sum_d \frac{1}{d} p_d (x) p_d (y)\right).
$$
where $x = (x_1, x_2, ...), y_j = (y_1, y_2, ...)$ and $p_d$ is the dth power sum. 

In \cite[Proposition 4.1]{Z1}, Zhou proves
$$
   s_\lambda(q^{-\rho}) = (-1)^{|\lambda|}q^{-\frac{1}{2}\kappa_\lambda} s_\lambda (q^\rho), 
$$
where $\kappa_\lambda = \sum_i \lambda_i(\lambda_i - 2i +1)$. We can therefore write the framing $-1$ topological vertex as
\begin{align}
    \Psi^{(-1)} &= \sum_\lambda q^{-\frac{1}{2} \kappa_\lambda} s_\lambda(q^\rho) s_\lambda(x) \notag\\
    &= \sum_\lambda s_\lambda(q^{-\rho}) s_\lambda(-x) \notag\\
    &= \e(\sum_m \frac{1}{m} \frac{1}{q^{-m/2} - q^{m/2}} p_m(-x) ) \notag \\				
    &= Q_{(0,1)} (1) \cdot 1 \label{-1} 
\end{align}
We now treat the case of a general framing.  Framing $(p-1)$ is obtained by applying $p$ Dehn twists to the framing $-1$ solution.
This transforms the homology class $(0,1)$ to $(p,1)$, but as described above
$Q_{(0,1)}(1)$ is transformed to $Q_{(p,1)}(a^{-p}),$
so the framing-change morphism should be effected by the operator $Q_{p,1}(a^{-p})$ on wavefunctions. 
This is proved in the following.

\begin{prop}
\label{prop:TV}
$$
Q_{p,1} (a^{-p}) \cdot 1 = \Psi^{(p-1)}
$$
\end{prop}

\begin{proof}
Let $q^{\frac{p}{2}\kappa}$ be the operator which acts on $W_\lambda$ by multiplication by  $q^{\frac{p}{2}\kappa_\lambda}$.  Then
\begin{align*}
	\Psi^{(p-1)}  &= q^{\frac{p}{2} \kappa } \cdot \Psi^{(-1)} \\
        & = \Ad_{q^{p\kappa/2}} Q_{0, 1} \cdot 1 \\
	&= \Ad_{q^{\frac{p}{2} \kappa }} \e \left( \sum_n \frac{1}{n}\frac{(-1)^{n+1}}{ \{n\} } P_{0,n} \right) \cdot 1\\
	&= \e \left(\sum_n \frac{1}{n} \frac{(-1)^{n+1}}{\{n\}} \Ad_{q^{p\kappa/2}} (P_{0, n})\right) \cdot 1
\end{align*}
while 
\begin{align*}
    Q_{p,1} (a^{-p}) \cdot 1 = \e\left(\sum_n \frac{(-1)^{n+1} a^{-pn}}{n\{n\} } P_{pn, n} \right) \cdot 1 
\end{align*}

Thus we only need to show 
\begin{equation}\label{eq: ad of kappa}
\Ad_{q^{p\kappa/2}} P_{(0,n)} = a^{-pn} P_{(pn, n)}
\end{equation}
as operators acting on the skein of solid torus. 

According to \cite[Theorem 4.6]{MS17}
, 
the RHS of \eqref{eq: ad of kappa} acts by:
\begin{align*}
 a^{-pn} P_{(pn, n)} \cdot W_\lambda &=  a^{-pn} (a^{pn} \frac{\{pn\}}{ \{pn^2\} } (\sum_{\alpha \in \lambda + n}  C_{\alpha - \lambda} (q^{pn})(-1)^{ht(\alpha - \lambda) } W_\alpha  ))\\
 &= \frac{\{pn\}}{ \{pn^2\} } (\sum_{\alpha \in \lambda + n}  C_{\alpha - \lambda} (q^{pn})(-1)^{ht(\alpha - \lambda) } W_\alpha  ). 
\end{align*}
while the left-hand side of \eqref{eq: ad of kappa} acts by
\begin{align*}
 \Ad_{q^{p\kappa/2}} P_{(0,n)} \cdot W_\lambda = \sum_{\alpha \in \lambda + n} q^{\frac{p}{2}(\kappa_\alpha - \kappa_\lambda)} (-1)^{ht(\alpha - \lambda)} W_\alpha  ,
\end{align*}

Recall $\alpha \in \lambda + n$ means that $\alpha - \lambda$ is a ``border strip'' of size $n$, i.e. it is connected and contains no $2 \times 2$ squares, $\{n\} = q^{n/2} - q^{-n/2}$, and $C_{\alpha - \lambda} (t) = \sum_{\box \in  \alpha - \lambda } t^{c(\Box)}  $.

It remains to show 
\begin{align*}
 C_{\alpha - \lambda}  (q^{pn}) = \frac{\{pn^2\}}{\{pn\}} q^{\frac{p}{2} (\kappa_\alpha - \kappa_\lambda)}.
\end{align*}
We can assume $p = 1$, and use $\sum_{\Box \in \lambda} c(\Box) = \frac{1}{2} \kappa_\lambda$. The above is equivalent to
\begin{align}
\sum_{\Box \in \alpha - \lambda} q^{n c(\Box)} 
= (q^{\frac{n}{2} (n-1 )} + q^{\frac{n}{2}(n-3) } + \cdots q^{- \frac{n}{2} (n-1) } ) q^{ \sum_{\Box \in \alpha - \lambda} c(\Box)} \label{border}
\end{align}
Since $\alpha - \lambda $ is a ``border strip'' of size $n$, if we denote the very lower-left box of this border strip by $\boxtimes$, we have
$$\sum_{\Box \in \alpha - \lambda} q^{n c(\Box)}  = q^{n c(\boxtimes)} + q^{n (c(\boxtimes) +1) } + q^{n (c(\boxtimes) +2) } + \cdots q^{n (c(\boxtimes) + n-1) }
$$
and
 $$
 \sum_{\Box \in \alpha - \lambda} c(\Box) = n c(\boxtimes) + \frac{n(n-1)}{2}. 
 $$ 
  Thus the two sides of \eqref{border} match up. 
\end{proof}

Observe that
\begin{align*}
    Q_{0,1}(-1)^{-1} \cdot 1 
    = \sum_\lambda s_\lambda( q^\rho ) s_\lambda(x) 
    = \Psi^{(0)}_{AV} = Q_{1, 1}(a^{-1}) \cdot 1 
\end{align*}
One can think this as starting from the genus $1$ necklace graph and mutate a negative edge to get the $1$-canoe graph. 

Indeed, we have:
\begin{prop} \label{prop: inverse of Q}
We have an identity
\begin{align*}
    Q_{p, 1}(-t)^{-1} \cdot 1 = Q_{p+1, 1}(t a^{-1}) \cdot 1. 
\end{align*}
In particular, 
\begin{align*}
    Q_{p, 1}(-a^{-p})^{-1} \cdot 1 = \Psi^{(p)}_{AV}. 
\end{align*}
\end{prop}
\begin{proof}
This is a straightforward consequence of ~\eqref{eq: ad of kappa}: we compute
   \begin{align*}
        Q_{p,1}(-t a^{-p})^{-1} \cdot 1 
        &= \e (\sum_n \frac{1}{n\{n\}} a^{-pn} t^n P_{pn, n}) \cdot 1\\
        &\stackrel{\eqref{eq: ad of kappa}}{=} \Ad_{q^{p\kappa /2 }} \cdot \e (\sum_n \frac{1}{n\{n\} } t^n  P_{0, n}) \cdot 1\\
        &= q^{p\kappa /2} \cdot \e (\sum_n \frac{1}{n} p_n(q^\rho) p_n(tx)) \\
        &= \sum_\lambda q^{p\kappa_\lambda/2} s_{\lambda}(q^\rho) s_\lambda(tx) \\
        &= \Psi_{AV}^{(p)} (tx) \\
        & = Q_{p+1, 1} (a^{-p-1} t) \cdot 1.
   \end{align*} 
   Replacing $t$ by $a^p t$ gives the desired formula.
\end{proof}

\begin{example}[Conormal of the unknot.]
\label{eg:unknot-conormal}
    Let $L$ be the conormal of the unknot in $T^*S^3$. 
    \cite{ES20} showed that its wavefunction should live in $\Sk(\cH_1) \otimes \Sk(S^3) $ and satisfy the equation (choosing an orientation of $S^3$):
    \begin{align}\label{eq: unknot}
        (\bigcirc - P_{1,0} -\gamma a_L a^{-1} P_{0,1} + \gamma a P_{1, 1} ) \Psi_L = 0
    \end{align}
  Here $a_L$ and $a$ are variables in $\Sk(L)$ and $\Sk(S^3)$ respectively, and $\gamma$ is a monomial in $a$ and $a_L$.  In \cite{ES20} this recursion was solved in a pure combinatorial way.
    
    Although this Lagrangian doesn't come from our geometric seeds, we can still use the skein-cluster method to solve the equation~\eqref{eq: unknot}. 
    Indeed by Lemma \ref{lem:Ad-action}, we have:
    \begin{align*}
        \bigcirc - P_{1,0} -\gamma a_L a^{-1} P_{0,1} + \gamma a P_{1, 1} &= \Ad_{Q_{0, 1}(-\gamma a)^{-1}} (\bigcirc - \gamma a_L a^{-1} P_{0,1} -  P_{1,0} ) \\
        &= \Ad_{Q_{0, 1}(-\gamma a)^{-1}} \Ad_{Q_{-1, 1} (\gamma a_L a^{-1} )^{-1}} (\bigcirc - P_{1,0}) 
    \end{align*}
    Thus we get:
    \begin{align*}
        (\bigcirc - P_{1,0}) \cdot Q_{-1, 1} (\gamma a_L a^{-1} ) \cdot Q_{0, 1}(-\gamma a) \cdot \Psi_L = 0
    \end{align*}
    which implies that the unique solution to \eqref{eq: unknot} is:
    \begin{align*}
        \Psi_L &= Q_{0, 1}(-\gamma a)^{-1} \cdot Q_{-1, 1} (\gamma a_L a^{-1} )^{-1} \cdot 1\\
        & \stackrel {\text{Prop \ref{prop: inverse of Q}}}{=} \frac{Q_{0, 1} (-\gamma a^{-1})}{Q_{0, 1}(-\gamma a)}  . 
    \end{align*}
\end{example}

\section{Reduction to finite rank skeins}
\label{sec:higher-rank}

{In this section, we investigate the reduction of our HOMFLYPT skein constructions to finite-rank skein algebras, which are quantizations of the coordinate ring of $GL_N$-character varieties for surfaces. We focus on the fundamental case of the punctured torus, which captures the local geometry of two transversely intersecting curves in a general surface.  Here the quantized coordinate ring of the character variety embeds as a subalgebra in that of the cluster Poisson moduli space $\cX_{GL_N,\Sigma_{1,1}}$  of framed local systems on the punctured torus. 
This framework allows one to express the image of Baxter operator in the finite rank skein as a sequence of standard $q$-dilogarithm mutations. We then use cluster theory to prove that the image of the pentagon identity from Conjecture~\ref{conj: pentagon relation for punctured torus} holds in the rank-$N$ punctured torus skein for all $N$. This provides strong evidence for Conjecture~\ref{conj: pentagon relation for punctured torus}, and could perhaps, with more work, be upgraded to give a proof.

From the standard cluster perspective, the case of rank $N$ skeins on the closed torus actually turns out to be more involved that for its punctured counterpart, and is related to the cluster structure on the $\mathfrak{gl}_N$ spherical DAHA -- a topic of work in progress of the second-named author in collaboration with di Francesco, Kedem, Shapiro~\cite{daha}. For rank $N\leq 2$, however, the relevant moduli space is simply a symplectic leaf in the familiar Fock-Goncharov moduli space $\cX_{GL_N,\Sigma_{1,1}}$. For this reason our discussion of the closed torus will be mostly limited to the case $N=2$, leaving the general case for future work.}


\subsection{Finite rank and the spherical DAHA}
As just mentioned, in the setting of a closed genus 1 Legendrian, the reduction to the quantized microlocal rank $N$ moduli space is effected by passing from the skein algebra $\Sk(T^2)$ to its quotient, the $\mathfrak{gl}_N$ spherical DAHA  $ SH_{\sigma = q,\overline\sigma = q}(\mathfrak{gl}_N)$.  We also specialize the framing parameter to $a  = q^{\frac{N}{2}}$, so that the unknot $\bigcirc= [N]_q$ acts by multiplication by the quantum dimension of the vector representation of $U_q(\mathfrak{gl}_N)$. 
In the interest of lightening the notation, we will often abuse notation and denote the elements $P_{\mathbf x}$ in $\Sk(T^2)$ and their images $P_{\mathbf x}^{(N)}$ in $SH_{q, q}(\mathfrak{gl}_N)$ by the same symbol. 

To motivate what follows, we begin by discussing some aspects of the geometry of the classical, rank 2 chromatic Lagrangian. First,  recall that the spherical DAHA for $GL(2)$ is a quantization of the coordinate ring of the moduli space of $GL(2)$-local systems on the punctured torus; under this identification, the operators representing $P_{(1,0)}$ and $P_{(0,1)}$ correspond to quantum traces in the defining representation of the holonomies around the corresponding cycles on the torus. Specializing the parameter $t=-1$ cuts out the closed torus character variety inside that for the punctured torus. 

If $\mathbf{x}$ is a primitive vector, we can use the Newton formula to define operators $E_{2\mathbf{x}}$ corresponding to the trace in the second fundamental (so for $N=2$, the determinant character) of the holonomy around $\mathbf{x}$:  
$$
E_{2\mathbf{x}} = \frac{1}{2}\left(P_{\mathbf{x}}^2-P_{2\mathbf{x}}\right).
$$
For $N=2$, one can check that these operators satisfy 
\begin{align*}
q^{\frac{1}{2}}E_{(0,2)}P_{(1,0)}^2 + q^{-\frac{1}{2}}P_{(1,1)}^2 + q^{-\frac{1}{2}}P_{(0,1)}^2 E_{(2,0)} - q^{-1}P_{(1,0)}P_{(0,1)}P_{(1,1)} - 2(q^{\frac{1}{2}}+q^{-\frac{1}{2}})E_{(0,2)}E_{(2,0)}=0
\end{align*}
When $q=1$, this becomes the defining equation cutting out the 4-dimensional $GL(2)$ character variety $\mathrm{Ch}_{GL_2}(T^2)$ as a hypersurface inside $\mathbb{C}^3\times (\mathbb{C}^*)^2$:
\begin{align*}
E_{(0,2)}P_{(1,0)}^2 + P_{(1,1)}^2 + P_{(0,1)}^2 E_{(2,0)} -P_{(1,0)}P_{(0,1)}P_{(1,1)} - 4E_{(0,2)}E_{(2,0)}=0
\end{align*}
The variety $\mathrm{Ch}_{GL_2}(T^2)$ has a 2-dimensional singular locus, corresponding to the reducible rank 2-local systems on $T^2$.

The polynomial representation $\mathcal{F}^{(N)}$ of $SH_{q, q}(\mathfrak{gl}_N)$ is the space of symmetric Laurent polynomials in $N$ variables $x_1,\ldots, x_N$. The vertical generator  $P_{(0,1)}$ acts on  $\mathcal{F}^{(N)}$ by the operator
\begin{align}
    P_{(0,1)}\cdot f(x_1,\ldots, x_N) = (x_1+\cdots+x_N)f(x_1,\ldots, x_N)
\end{align}
of multiplication by the first elementary symmetric function $e_1(x_1\ldots,x_N)$ in the $x_i$.
The horizontal generator $P_{(1,0)}$ acts by a $q$-difference operator: writing
$$
Y_i\cdot f(x_1,\ldots,x_N) = f\left(x_1,\ldots,  qx_i,\ldots,x_N\right),
$$
the action of $P_{(1,0)}$ on $\mathcal{F}^{(N)}$ is given by
\begin{align}
P_{(1,0)}\longmapsto 
\sum_{i=1}^N\prod_{j\neq i} \frac{q^{\frac{1}{2}}x_i-q^{-\frac{1}{2}}x_j}{x_i-x_j}Y_i.
\end{align}
In other words, $P_{(1,0)}$ acts by $q^{\frac{1-N}{2}}M_1$, where $M_1$ is the $t=q$ specialization of the first \emph{Macdonald difference operator} in $N$ variables.
For example, when $N=2$ we have
\begin{align}
P_{(1,0)}^{(2)}\mapsto 
 \frac{q^{1/2}x_1-q^{-1/2}x_2}{x_1-x_2}Y_1 + \frac{q^{1/2}x_2-q^{-1/2}x_1}{x_2-x_1}Y_2,
\end{align}
The action of the operators $P_{(1,0)}, P_{(0,1)}$ preserves the subspace
\begin{align*}
\mathcal{F}^{(N)} &= \mathbb{C}(q)[x^\pm_1,\ldots, x^\pm_N]^{S_N}
\\
&=\bigoplus_{l(\lambda)\leq N}\mathbb{C}(q)\cdot W_\lambda
\end{align*}
which is spanned by the Schur polynomials $W_\lambda(x_1,\ldots x_N)$ labelled by all generalized partitions $\lambda=(\lambda_1\geq\lambda_2\geq\cdots\geq \lambda_N)$. In the running example $N=2$, the $W_\lambda$ for $\lambda=(\lambda_1\geq\lambda_2\geq0)$ are given explicitly by
\begin{align*}
W_\lambda(x_1,x_2) &= \frac{x_1^{\lambda_1+1}x_2^{\lambda_2} - x_1^{\lambda_2}x_2^{\lambda_1+1}}{x_1-x_2}\\
&= x_1^{\lambda_1}x_2^{\lambda_2} + x_1^{\lambda_1-1}x_2^{\lambda_2+1}+\cdots + x_1^{\lambda_2}x_2^{\lambda_1}.
\end{align*}
The Macdonald operators act diagonally in the basis of Macdonald polynomials, whose $t=q$ specializations are just the Schur polynomials:
$$
M_1 \cdot W_\lambda = \left(\sum_{k=1}^Nq^{\lambda_k+N-k}\right)W_\lambda .
$$
Using the combinatorial identity~\eqref{eq:content-sum} and recalling our specialization $a\mapsto q^{\frac{N}{2}}$, we can therefore rewrite the action of $P_{(1,0)}$ in the form
\begin{align*}
P_{(1,0)}\cdot W_{\lambda} 
&=q^{\frac{1-N}{2}}M_1\cdot W_\lambda\\
&=\left(\sum_{k=1}^Nq^{\lambda_k-k+\frac{N+1}{2}}\right)W_\lambda \\
&= \left(\bigcirc + a(q^{\frac{1}{2}}-q^{-\frac{1}{2}})c_\lambda(q)\right)W_{\lambda}
\end{align*}
which recovers the description~\eqref{eq:P-eigenvalues} of the eigenvalues. 
Similarly, by the Pieri rule for Schur polynomials we have
$$
P_{(1,0)}\cdot W_{\lambda}(x_1,\ldots, x_N) =\sum_{\mu=\lambda+\square}W_\mu(x_1,\ldots,x_N),
$$
where now the only nonzero terms in the sum are partitions with at most $N$ rows.

As a difference operator on the space of $N$-variable symmetric functions, the face relation for the framing $p=-1$  wavefunction $\Psi(x_1,\ldots,x_N)$ reads
\begin{align}
\left([N]_{q} - q^{\frac{1-N}{2}}M_1 + q^{\frac{N}{2}}e_1(x_1,\ldots, x_N) \right)\cdot \Psi(x_1,\ldots,x_N) = 0,
\end{align}
which has the unique solution
$$
\Psi^{(-1)}(x_1,\ldots,x_N) = \prod_{k=1}^N\Phi(X_k).
$$
where $\Phi$ is the standard quantum dilogarithm function from~\eqref{eq:qdl-def}.
We can identify the piece of the rank 2 chromatic Lagrangian living over the canoe graph as the locus in $\mathrm{Ch}_{GL_2}(T^2)$ cut out by the specialization of the annihilator of the framing $(-1)$-wavefunction $\Psi^{(-1)}=\Phi(X_1)\Phi(X_2)$. 
In terms of the elements $X_e\in GL_N/\mathrm{Ad}_{GL_N}$ described in item~\eqref{item-edge-fn} of Section~\ref{sec:chromatic}, the classical limit of the canoe face relation~\eqref{eq:g1-canoe}
reads
$$
\mathbf{1}_{GL_2} = X_{(1,0)} - X_{(0,1)} \in GL_N/\mathrm{Ad}_{GL_N}.
$$
Taking the trace of this relation in the two fundamental representations of $GL_2$ yields two polynomial functions
$$
A_{\omega_1} = 2 - P_{(1,0)} + P_{(0,1)}, \quad A_{\omega_2} = E_{(0,2)} +P_{(1,0)} - 1 - E_{(2,0)} 
$$
which vanish on the rank 2 chromatic Lagrangian. However, these elements alone do not suffice to generate the local chromatic ideal $\mathcal{I}_{canoe}$. Indeed, the ideal $\langle A_{\omega_1},A_{\omega_2}\rangle$ has two associated primes, one of which is the  ideal
\begin{align}
\label{eq:skein-ideal}
\mathcal{I}_{canoe} = \langle A_{\omega_1},~ &A_{\omega_2},P_{(1,1)} -(E_{(2,0)}-E_{(0,2)})^2 +E_{(2,0)}+E_{(0,2)} \rangle,
\end{align}
given by the classical limit of the annihilator of the {framing $(-1)$-}canoe wavefunction $\Psi^{(-1)}=\Phi(X_1)\Phi(X_2)$.
We note that no two of these generators suffice to generate  the ideal $\mathcal{I}_{canoe}$ in the classical coordinate ring of the character variety. However in view of the Poisson brackets 
$$
\{P_{(1,0)},P_{(0,1)}\} = P_{(1,1)}, \quad \{E_{(2,0)},P_{(0,1)}\} = E_{(2,0)}P_{(0,1)}, \quad \{E_{(0,2)},P_{(1,0)}\} = -E_{(0,2)}P_{(1,0)},
$$
we see that $A_{\omega_1},A_{\omega_2}$ do nonetheless generate $\mathcal{I}_{canoe}$ as a Poisson ideal.

The choice of Schur basis $\{W_\lambda\}$ for $\mathcal{F}^{(N)}$ identifies the latter with a subspace of the space of functions on the set $\mathrm{Par}_N$ of partitions with $l(\lambda)\leq N$:
\begin{align*}
&\{ \phi:\mathrm{Par}_N\rightarrow \mathbb{C}(q)\}\simeq \mathcal{F}^{(N)}\\
&\phi \longmapsto \sum_{\lambda}\phi(\lambda)W_\lambda ,
\end{align*}
so that a function $\phi$ encodes the expansion coefficients of a symmetric function into the Schur basis. In this way, we can embed $\mathcal{F}^{(N)}$ into the vector space
$$
\mathcal{H}^{(N)}:=\{ \phi:\mathbb{Z}^N\rightarrow \mathbb{C}(q)\},
$$
as the subspace $\mathcal{F}^{(N)}\subset\mathcal{H}^{(N)}$ consisting of all functions that vanish at all but finitely many points, and which vanish everywhere outside of the cone of generalized partitions $\mathrm{Par}_N\subset\mathbb{Z}^N$ in the lattice $\mathbb{Z}^N$.

The polynomial representation of $SH_{q,q}$ delivers a kind of `abelianization' of the skein algebra and its module defined by the solid torus, although one somewhat different from those familiar from cluster theory: the generators of the algebra act by not by Laurent polynomials in a quantum torus, but rather by elements in a localization of such an algebra, owing to the denominators in the Macdonald difference operators. 

We will now explain how to construct an alternative abelianization in which the generators do indeed correspond to genuine quantum Laurent polynomials, and the Baxter operator associated to the flip of triangulation is realized by the action of a product of quantum dilogarithms. 

Here we will do this explicitly only for $N=2$ where we can appeal to the results of Fock and Goncharov \cite{FG1} on the cluster structure of $\cX_{S,PGL_2}$ where $S=T^2\setminus D^2$ is the punctured torus. Extending to the higher rank case involves significantly more complicated cluster combinatorics, {the latter being the subject of the work in progress~\cite{daha}.} 

%

It is convenient to work with the singular character variety $\mathrm{Ch}_{GL_2}(T^2)$ and its quantization by means of the framed moduli space $\mathscr{X}_{GL_2,S}$ of local systems with Borel decoration at the puncture -- indeed, the latter is a rational variety which admits a cluster parametrization.
As explained in Section 2.4 of~\cite{SSZ}, the moduli space $\mathscr{X}_{PGL_2,S}$ is a genus 1 analog of the Springer resolution $T^*\mathbb{P}^1\rightarrow \mathcal{U}$ resolving the singular variety of unipotent elements in the group $PGL_2(\mathbb{C})$. 

Thus if one has a sheaf on the classical moduli space $\mathrm{Ch}_{GL_2}(T^2)$ (for example, the classical limit of the skein module for the solid torus), one can pull back that sheaf to the symplectic tori given by intersecting the cluster charts on $\mathscr{X}_{GL_2,S}$ with the appropriate symplectic leaf -- see Proposition 2.10 of \cite{SSZ} for details.  For the quantum moduli space, such localizations will produce a compatible system of modules over the quantum tori associated to different charts in the cluster atlas and define an object in the category of descent data introduced in~\cite{SSZ},  Remark 4.4. The practical upshot is that one can study the skein module $\Sk(D^2\times S^1)$ by means of its \emph{localizations} to modules over the different quantum tori in the cluster atlas, which will have the advantage that the generators of $\Sk(T^2)$ act by \emph{universally Laurent} elements of the quantum tori, with the intertwiners between charts being given by standard quantum cluster transformations.

%

\subsection{\texorpdfstring{$q$}{q}-Whittaker basis and abelian clusters for the skein algebra}
At the quantum level, we construct the cluster abelianization by passing from the Schur basis $\{W_\lambda\}$ to another basis $\{R_\lambda\}$ for the ring of symmetric functions, the $q$-Whittaker functions. They are characterized by the  initialization $R_{(0,0)}=1$ and their Pieri rules
\begin{align}
\label{eq:qW-pieri}
(x_1+x_2)R_{\lambda_1,\lambda_2} &=  R_{{\lambda_1+1,\lambda_2}}+ (1-q^{\lambda_1-\lambda_2})R_{{\lambda_1,\lambda_2+1}}, \\ \nonumber (x_1x_2)R_{\lambda_1,\lambda_2} &=  R_{\lambda_1+1,\lambda_2+1} ,
\end{align}
Note that every term in the Pieri rule~\eqref{eq:qW-pieri} indeed corresponds to a partition: the only time we cannot add a box to row 2 while maintaining the partition inequality is when $\lambda_1=\lambda_2$, and in this case the coefficient $(1-q^{\lambda_1-\lambda_2})$ of $R_{{\lambda_1,\lambda_2+1}}$ in~\eqref{eq:qW-pieri} vanishes. Hence the Pieri rule for the $R_\lambda$ takes the form of a $q$-weighted sum over partitions $\mu$ with at most 2 parts obtained from $\lambda$ by adding a single box. One checks that in the $q$-Whittaker basis, the canoe wavefunction 
$$
\Psi(X_1,X_2)=\Phi(X_1)\Phi(X_2)
$$ is expanded as
\begin{align}
\label{eq:whittaker-wavefunction}
\Psi = \sum_{n\geq0} \frac{(-q)^{\frac{n}{2}}}{\prod_{k=1}^n(1-q^{k})}R_{(n,0)},
\end{align}
with the sum being reduced to one over partitions with a single row. 

The basis $\{R_\lambda\}$ for $\mathcal{F}^{(2)}$ again identifies the latter with a space of functions $\phi:\mathrm{Par}_N\rightarrow\mathbb{C}(q)$, which as before we embed into $\mathcal{H}^{(2)}=\{\phi:\mathbb{Z}^N\rightarrow\mathbb{C}(q)\}$ as the subspace of all functions vanishing outside the cone of partitions. Let us introduce the following difference operators on $\mathcal{H}^{(2)}$:
$$
V_1\cdot \phi(\lambda_1,\lambda_2) = \phi(\lambda_1-1,\lambda_2), \quad V_2\cdot \phi(\lambda_1,\lambda_2) = \phi(\lambda_1,\lambda_2-1),
$$
$$
 U_i\cdot\phi(\lambda) = q^{\lambda_i}\phi(\lambda),
$$
These operators generate a standard quantum torus
$$
U_iV_j=q^{\delta_{i,j}}V_iU_j.
$$
In terms of this quantum torus, the Pieri rules take the form
\begin{align}
\label{eq:toda}
(x_1+x_2)\cdot \phi &= \left(V_1 + V_2-qU_1U_2^{-1}V_2\right)\cdot\phi\\
\nonumber x_1x_2\cdot\phi &= V_1V_2\cdot\phi,
\end{align}
so that as an operator we have
$$
P_{(0,1)} \mapsto V_1 + V_2-qU_1U_2^{-1}V_2.
$$
The operators~\eqref{eq:toda} are the Hamiltonians $H_1,H_2$ of the $U_q(\mathfrak{gl}_2)$ open Toda chain, which arise as the quantum traces of holonomies around the $(0,1)$-curve on the torus. 

The Macdonald operator $M_1$, on the other hand, acts in the $R_\lambda$ basis as follows:
$$
M_1\cdot R_{\lambda_1,\lambda_2} = \left(  q^{\lambda_1+1} + q^{\lambda_2}\right)R_{\lambda_1,\lambda_2}  -q^{\lambda_2+1}(1-q^{\lambda_1-\lambda_2})(1-q^{\lambda_1-\lambda_2-1})R_{\lambda_1-1,\lambda_2+1}.
$$
Note that again each term on the right-hand-side either vanishes or corresponds to a genuine partition. So at the level of the expansion coefficients $\phi:\mathrm{Par}_N\rightarrow\mathbb{C}(q)$, 
the element $P_{(1,0)}$ acts by
$$
P_{(1,0)}\cdot\phi = q^{-\frac{1}{2}}\left\{qU_1 + U_2 - (1-qU_1U_2^{-1}) (1-q^2U_1U_2^{-1})U_2V_1^{-1}V_2\right\}\cdot\phi
$$
Thus by passing to the $q$-Whittaker basis we have constructed an embedding of the skein algebra into the quantum torus $\mathcal{T}^q_{U,V}$ generated by $U_i,V_i$. The annihilator in $\mathcal{T}^q_{U,V}$ of the wavefunction $\Psi^{(-1)}=\Phi(X_1)\Phi(X_2)$ regarded as an element of $\mathcal{H}^{(2)}$ is simply the left ideal 
\begin{align}
\mathcal{I}_{UV}:=\label{eq:ideal-uv}
\mathcal{T}^q_{U,V}\cdot \langle U_2-1, ~1-U_1 + q^{\frac{1}{2}}V_1\rangle \subset \mathcal{T}^q_{U,V},
\end{align}
in contrast to its more complicated annihilator in $\mathrm{Sk}(T^2)$ whose classical limit is ~\eqref{eq:skein-ideal}.

In fact, the image of $\Sk(T^2)$ is universally Laurent with respect to the standard Fock-Goncharov cluster structure on $\mathcal{X}_{PGL_2,\Sigma_{1,1}}$. To illustrate this, let us obtain a more symmetric realization of the skeins associated to our homology basis elements by performing the two `negative mutations' given by  conjugation by the product of dilogarithms $\Phi(m_1)\Phi(m_2)$, where
$$
m_1 = -q^{\frac{1}{2}}U_1V_1V_2^{-1}U_2^{-1}, \qquad m_2 = -q^{\frac{1}{2}}U_1U_2^{-1}.
$$
One checks that in this new coordinate system we have
\begin{align}
\label{eq:symm-embed}
    \mathrm{Ad}_{\Phi(m_1)\Phi(m_2)}\cdot P_{(1,0)} &= -q^{\frac{1}{2}}V_2^{-1}V_1U_1 + q^{-\frac{1}{2}}U_2 - q^{\frac{1}{2}}V_1^{-1}V_2U_2,\\
       \nonumber \mathrm{Ad}_{\Phi(m_1)\Phi(m_2)}\cdot P_{(0,1)} &= V_1+V_2 -U_1V_1^2V_2^{-1}U_2^{-1}\\
              \nonumber \mathrm{Ad}_{\Phi(m_1)\Phi(m_2)}\cdot P_{(1,1)} &= U_1V_1 + V_2U_2 - q V_1^{-1}V_2^2U_2.
\end{align}

Figure~\ref{fig:loc-quiver} illustrates  the corresponding quiver that describes the `local' cluster structure on $\mathcal{X}_{S,PGL_2}$ for the cluster torus in which the embedding~\eqref{eq:symm-embed} takes its values:

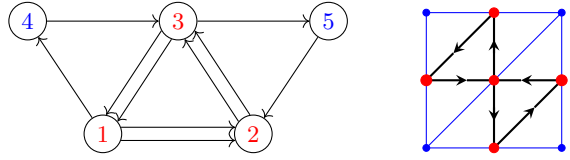
\begin{figure}[ht]

\begin{tikzpicture}[every node/.style={inner sep=0, minimum size=0.5cm, circle, draw, fill=white}, x=1cm, y=0.75cm]

\node (1) at (2,2) {\color{blue}\scriptsize 5};
\node (2) at (1,0) {\color{red}\scriptsize 2};
\node (3) at (0,2) {\color{red}\scriptsize 3};
\node (4) at (-1,0) {\color{red}\scriptsize 1};
\node (5) at (-2,2) {\color{blue}\scriptsize 4};

\draw [->] (2.100) -- (3.-30);
\draw [->] (2.140) -- (3.-70);
\draw [->] (3.-110) -- (4.40);
\draw [->] (3.-150) -- (4.80);
\draw [->] (4.-20) -- (2.-160);
\draw [->] (4.20) -- (2.160);
\draw [->] (1) -- (2);
\draw [->] (3) -- (1);
\draw [->] (5) -- (3);
\draw [->] (4) -- (5);
\end{tikzpicture}
\qquad
\begin{tikzpicture}[rotate=90,xscale=-.9,yscale=.9]
    \draw[blue] (-1,-1)--(1,-1)--(1,1)--(-1,1)--(-1,-1)--(1,1);
    \draw[fill,blue] (-1,-1) circle (.05cm);
    \draw[fill,blue] (-1,1) circle (.05cm);
    \draw[fill,blue] (1,-1) circle (.05cm);
    \draw[fill,blue] (1,1) circle (.05cm);

    \pgfmathsetmacro{\a}{.6}
    \pgfmathsetmacro{\b}{1-\a}
    
    \draw[black,thick] (\b,0)--(1,0);
    \draw[black,thick,-stealth] (0,0)--(\a,0);
    
    \draw[black,thick,-stealth] (0,0)--(-\a,0);
    \draw[black,thick] (-\b,0)--(-1,0);
    
    \draw[black,thick,-stealth] (0,-1)--(0,-\b);
    \draw[black,thick] (0,-\b)--(0,0);
    
    \draw[black,thick,-stealth] (0,1)--(0,\b);
    \draw[black,thick] (0,\b)--(0,0);
    
    \draw[black,thick,-stealth] (-1,0)--(-\b,\a);
    \draw[black,thick] (-\b,\a)--(0,1);
    
    \draw[black,thick,-stealth] (1,0)--(\b,-\a);
    \draw[black,thick] (\b,-\a)--(0,0-1);

    \draw[red,fill] (0,0) circle (.07cm);
    \draw[red,fill] (1,0) circle (.07cm);
    \draw[red,fill] (-1,0) circle (.07cm);
    \draw[red,fill] (0,-1) circle (.08cm);
    \draw[red,fill] (0,1) circle (.08cm);
    
\end{tikzpicture}
\caption{Left:  Quiver describing the local cluster structure at the canoe graph for the rank 2 chromatic Lagrangian; the red vertices are regarded as mutable and the blue as frozen (with respect to the canoe).  Right: corresponding ideal triangulation (in blue) of once-punctured torus.}
\label{fig:loc-quiver}
\end{figure}
Let us equip this quiver with the structure of a framed seed in the sense of~\cite{SSZ} as follows:
\begin{align*}
&X_{e_1} \mapsto -q^{-\frac{1}{2}}V_1V_2^{-1} , \quad X_{e_2} \mapsto -q^{-\frac{1}{2}}U_1^{-1}V_1^{-2}V_2^2U_2 \\
& X_{e_3}= -q^{\frac{1}{2}} V_2^{-1}U_2^{-1}V_1U_1, \quad X_{e_4} \mapsto -q^{\frac{1}{2}}V_{1}^{-1}V_2U_2 , \quad X_{e_5}\mapsto V_1.
\end{align*}
Then we have
\begin{align}
\label{eq:cluster-embed1}
P_{(1,0)} &\longmapsto X_{e_4}+ X_{e_4+e_1} + X_{e_4+e_1+e_3},\\
\nonumber P_{(0,1)} &\longmapsto X_{e_5}+ X_{e_5+e_3} + X_{e_5+e_3+e_2},
\end{align}
and we can now identify the two negative mutations we performed to reach this cluster as being done first at vertex 2, and then at vertex 3.

We conclude this discussion by noting that the solution in $\mathcal{H}^{(2)}$ to the difference equations~\eqref{eq:ideal-uv} can be determined using the cluster structure on $\mathcal{X}_{PGL_2,\Sigma_{1,1}}$ (i.e. without recourse to the Baxter operator) using the familiar strategy of mutating back to a simpler ideal. Indeed, let us choose an equivalent second generator
$$
q^{-\frac{1}{2}}V_1^{-1}(1-U_1U_2^{-1}) + 1
$$
 for $\mathcal{I}_{UV}$ obtained by multiplying the original from the left by $q^{-\frac{1}{2}}V_1^{-1}$ and using on the right the relation $U_2=1$. Then conjugating this new generator by $\Phi(m_2)$ we arrive at the element $1+q^{-\frac{1}{2}}V_1^{-1}$. Hence the ideals \begin{align*}
   \mathcal{I}_{UV}:= &=\mathcal{T}^q_{U,V}\cdot \langle U_2-1, ~1-U_1 + q^{\frac{1}{2}}V_1\rangle\\
\mathcal{I}'_{UV}:=&=\mathcal{T}^q_{U,V}\cdot \langle U_2-1, ~1+q^{-\frac{1}{2}}V_1^{-1}\rangle
\end{align*}
glue together along the mutation at vertex 2.

At the level of the representation $\mathcal{F}^{(2)}$, the intertwiner corresponding to this mutation is realized by the multiplication operator 
$$
\phi(\lambda_1,\lambda_2)\mapsto \frac{\phi(\lambda_1,\lambda_2)}{\prod_{k=1}^{\lambda_1-\lambda_2}(1-q^{k})}.
$$
Hence the difference equations coming from the generators $U_2-1,1+q^{-\frac{1}{2}}V_1^{-1}$ for the mutated ideal tell us that the Fourier coefficients of the wavefunction with respect to the renormalized $q$-Whittaker basis
$$
\widetilde{R}_\lambda = \frac{R_{\lambda}}{\prod_{k=1}^{\lambda_1-\lambda_2}(1-q^{k})}
$$
are simply $\phi(\lambda_1,\lambda_2) = \delta_{\lambda_2,0}(-q)^{\frac{\lambda_1}{2}}$, and we recover the solution~\eqref{eq:whittaker-wavefunction}.

\subsection{Abelianization of the Baxter operator}
We now explain how to abelianize the Baxter operator and thereby realize the flip from the genus 1 necklace to the canoe graph as a composite of cluster transformations in the traditional (i.e. non skein-theoretic) sense. Doing so will involve performing mutations at the blue vertices in the quiver of Figure~\ref{fig:loc-quiver}, which were regarded as frozen with respect to the local cluster structure associated to the canoe chart.

Let us first recall, following~\cite{KN11} and~\cite{Kel}, a general scheme for deducing dilogarithm identities from relations in the cluster modular groupoid. Fix an initial cluster seed $\bc_0$ corresponding to a basis $\Pi_0=\{e_1,\ldots, e_d\}$ for lattice $\Lambda$. Let us say that a vector $\lambda=\sum\lambda_ie_i\in\Lambda$ is \emph{tropically positive} (resp. tropically negative) with respect to the basis $\Pi_0=\{e_i\}$ if each of its coordinates $\lambda_i$  is positive (resp. negative). 
\begin{proposition}[Tropical sign coherence]
\label{prop:sgn-coh}
Suppose that $\vec\mu = \mu_{i_l}\circ\cdots\circ\mu_{i_1}$ is a sequence of mutations taking seed $\bc_0$ to some other seed $\bc=\vec{\mu}(\bc_0)$. Then there is a unique sequence of signs $(\epsilon_1,\ldots, \epsilon_l)$ such that for all $1\leq k\leq l$, the vector
$$
f_k = \nu^{\epsilon_{k-1}}_{i_{k-1}}\circ\cdots \nu^{\epsilon_{1}}_{i_{1}}(e_k)
$$
is tropically positive with respect to $\Pi_0$.
\end{proposition}
The elements of the basis 
$$
\vec\mu(\Pi_0):=\{\vec\mu(e_j)\}_{k=1}^d, \quad \vec\mu(e_j)=\nu^{\epsilon_{l}}_{i_{l}}\circ\cdots \nu^{\epsilon_{1}}_{i_{1}}(e_j)
$$
for $\Lambda$ are called the \emph{tropical $\cX$-variables} (or \emph{c-vectors}) associated to the seed $\bc$ (with respect to basepoint $\bc_0$).  By Proposition~\ref{prop:sgn-coh}, each tropical $\cX$-variable is either tropically positive or tropically negative with respect to $\Pi_0$. 

The sign-coherent automorphism part of the cluster transformation  $\vec\mu$ is given by
$$
\vec\mu = \Ad_{\Phi_{\vec\mu}}, \quad \Phi_{\vec\mu}=\Phi(X_{f_1})^{\epsilon_1}\circ \cdots\circ\Phi(X_{f_l})^{\epsilon_l}.
$$
The following important result says that the tropical $\cX$-variables suffice to detect when two sequences of mutations induce the same cluster transformation.
\begin{proposition}\cite[Theorem 3.5]{KN11}
\label{trop-criterion}
Let $\vec{\mu}_1,\vec{\mu}_2$  be two sequences of mutations based at $\Pi_0$. Then $\Phi_{\vec\mu_1}= \Phi_{\vec\mu_2}$ if and only if the (non-ordered) sets of tropical $\cX$-variables $\vec\mu_1(\Pi_0)$ and $\vec\mu_2(\Pi_0)$ are identical.
\end{proposition}
We now turn to the problem of realizing the Baxter operator $Q_{(0,1)}(z)$ cluster-theoretically, which amounts to determining its image under the algebra homomorphism~\eqref{eq:symm-embed}.  By its definition, the Baxter operator in the non-reduced algebra $\mathcal{E}_{q,q}$ is the solution of the difference equation
\begin{align}
Q_{(0,1)}(qz) = \Gamma_+(z)Q_{(0,1)}(z),
\end{align}
where $\Gamma_+(z)$ is the vertex operator
$$
\Gamma_+(z) = \exp\left(-\sum_{n>0} \frac{\left(-q^{\frac{1}{2}}z\right)^n}{n}P_{(0,n)}\right).
$$
Now recall that in the rank-$N$ quotient, the generators $P_{(0,k)}$ act on the ring of symmetric polynomials in $x_1,\ldots, x_N$ by multiplication by the $k$-th power sum symmetric function. Hence the image $\Gamma^{(N)}_+(z)$ of the vertex operator $\Gamma_+(z)$ in this quotient becomes a degree $N$ polynomial given by the generating function of the elementary symmetric polynomials:
$$
\Gamma^{(N)}_+(z) = \prod_{k=1}^N\left(1+q^{\frac{1}{2}}zx_k\right)=:\sum_{j=0}^Ne_{(0,k)}q^{\frac{j}{2}}z^j.
$$
On the other hand, we know by the Pieri rule for Schur functions that the element $e_{(k,0)}$ acts on the space of Taylor coefficients $\mathcal{F}^{(N)}=\{\phi:\mathrm{Par}_N\rightarrow \mathbb{C}(q)\}$ of symmetric functions with respect to the Schur basis by the $k$-th Hamiltonian of the $\mathfrak{gl}_N$ $q$-difference open Toda chain.

This difference equation has been solved explicitly in~\cite{SS18}, where the solution was shown to be given by a product of quantum dilogarithms corresponding to a sequence of cluster transformations. Rewriting the operator from from Definition 4.2 of ~\cite{SS18} with respect to our present conventions, the specialization of this solution for $N=2$
reads
$$
Q^{(2)}_{(0,1)}(1) = \Phi(V_1)\Phi(-U_1V_1^2V_2^{-1}U_2^{-1})\Phi(V_2).
$$
Returning to the quiver in Figure~\ref{fig:loc-quiver}, we see that $Q^{(2)}_{(0,1)}$ is identified with the operator defining the automorphism part of the sequence of signed mutations $\mu_{2,+}\circ\mu_{3,+}\circ\mu_{5,+}$. 

The Baxter operators $Q_{(1,0)}$ and $Q_{(1,1)}$ can be similarly determined:
\begin{align*}
Q^{(2)}_{(1,0)}(1)  &= \Phi(-q^{\frac{1}{2}}V_{1}^{-1}V_2U_2)\Phi(q^{-\frac{1}{2}}U_2)\Phi(-q^{\frac{1}{2}}V_2^{-1}V_1U_1),\\
Q^{(2)}_{(1,1)}(1) &= \Phi(U_1V_1)\Phi(- q V_1^{-1}V_2^2U_2)\Phi(V_2U_2).
\end{align*}
Using these formulas, one sees that the composite $Q_{(1,0)}Q_{(0,1)}$ corresponds to the sign-coherent automorphism part of the sequence mutations $\vec\mu_1:=\mu_{2}\circ\mu_{1}\circ\mu_{5}\circ \mu_{3}\circ \mu_{1}\circ\mu_{4}$. On the other hand, the composite $Q_{(0,1)}Q_{(1,1)}Q_{(1,0)}$ is realized as the sign-coherent automorphism part of the mutation sequence $\vec\mu_2:=\mu_{5}\circ\mu_{3}\circ\mu_{2}\circ\mu_{1}\circ\mu_{3}\circ\mu_{4}\circ \mu_{2}\circ\mu_{3}\circ\mu_{5}$.

One easily checks that these two sequences lead to clusters with identical (up to permutation) sets of tropical $\cX$-variables
$$
\{e_1,e_2,e_3, -(e_1+e_2+e_4), -(e_2+e_3+e_5)\}
$$ 
Thus Proposition~\ref{trop-criterion} implies that the elements $Q_{(0,1)}Q_{(1,1)}Q_{(1,0)}$ and $Q_{(1,0)}Q_{(0,1)}$ defining the sign-coherent automorphism parts of the two sequences are also equal. We conclude that the rank 2 specialization of the pentagon identity from Theorem~\ref{thm:baxter-pentagon}
is in fact an identity in the standard cluster modular group of the quiver in Figure~\ref{fig:loc-quiver}.

\subsection{The pentagon identity in arbitary finite rank} \label{subsec: pentagon in finite rank}
\begin{figure}[ht]
\begin{tikzpicture}[scale=.6]

	\pgfmathsetmacro{\N}{5}
	\pgfmathsetmacro{\k}{\N-1}
	\pgfmathsetmacro{\Nmt}{\N-2}
    \pgfmathsetmacro{\dx}{4}
    \pgfmathsetmacro{\dy}{4}
    \pgfmathsetmacro{\eps}{0.2}
    \pgfmathsetmacro{\varoffset}{0.05}
     \pgfmathsetmacro{\baxoffset}{3}


\node[draw,circle] (1_1) at (1*\dx,1*\dy ) { \scriptsize${x_{1,1}}$};
\node[draw,circle] (2_2) at (2*\dx,2*\dy ) {\scriptsize ${x_{2,2}}$};
\node[draw,circle] (3_3) at (3*\dx,3*\dy ) { \scriptsize${x_{3,3}}$};
\node[draw,circle] (4_4) at (4*\dx,4*\dy ) { \scriptsize${x_{4,4}}$};
\node[draw,circle] (2_1) at (2*\dx,1*\dy ) {\scriptsize ${x_{2,1}}$};
\node[draw,circle] (3_2) at (3*\dx,2*\dy ) {\scriptsize ${x_{3,2}}$};
\node[draw,circle] (4_3) at (4*\dx,3*\dy ) { \scriptsize${x_{4,3}}$};
\node[draw,circle] (0_4) at (0*\dx,4*\dy ) {\scriptsize ${x_{0,4}}$};
\node[draw,circle] (3_1) at (3*\dx,1*\dy ) {\scriptsize ${x_{3,1}}$};
\node[draw,circle] (4_2) at (4*\dx,2*\dy ) {\scriptsize ${x_{4,2}}$};
\node[draw,circle] (0_3) at (0*\dx,3*\dy ) { \scriptsize${x_{0,3}}$};
\node[draw,circle] (1_4) at (1*\dx,4*\dy ) {\scriptsize ${x_{1,4}}$};
\node[draw,circle] (4_1) at (4*\dx,1*\dy ) {\scriptsize ${x_{4,1}}$};
\node[draw,circle] (0_2) at (0*\dx,2*\dy ) {\scriptsize ${x_{0,2}}$};
\node[draw,circle] (1_3) at (1*\dx,3*\dy ) { \scriptsize${x_{1,3}}$};
\node[draw,circle] (2_4) at (2*\dx,4*\dy ) {\scriptsize ${x_{2,4}}$};
\node[draw,circle] (0_1) at (0*\dx,1*\dy ) {\scriptsize ${x_{0,1}}$};
\node[draw,circle] (1_2) at (1*\dx,2*\dy ) {\scriptsize ${x_{1,2}}$};
\node[draw,circle] (2_3) at (2*\dx,3*\dy ) { \scriptsize${x_{2,3}}$};
\node[draw,circle] (3_4) at (3*\dx,4*\dy ) {\scriptsize ${x_{3,4}}$};
\node[draw,circle] (0_0) at (0*\dx,0*\dy ) {\scriptsize ${x_{0,0}}$};
\node[draw,circle] (1_0) at (1*\dx,0*\dy ) {\scriptsize ${x_{1,0}}$};
\node[draw,circle] (2_0) at (2*\dx,0*\dy ) {\scriptsize ${x_{2,0}}$};
\node[draw,circle] (3_0) at (3*\dx,0*\dy ) {\scriptsize ${x_{3,0}}$};
\node[draw,circle] (0_0) at (0*\dx,0*\dy ) {\scriptsize ${x_{0,0}}$};
\node[draw,circle] (1_0) at (1*\dx,0*\dy ) {\scriptsize ${x_{1,0}}$};
\node[draw,circle] (2_0) at (2*\dx,0*\dy ) {\scriptsize ${x_{2,0}}$};
\node[draw,circle] (3_0) at (3*\dx,0*\dy ) {\scriptsize ${x_{3,0}}$};


    \foreach \j in {2,...,\k}{
        \pgfmathtruncatemacro{\UpperLimit}{\j-1}
        \foreach \i in {1,...,\UpperLimit}{
            \pgfmathtruncatemacro{\Nexti}{\i-1}
            \draw[->] (\Nexti_\j) -- (\i_\j);
        }
    }
    
    \foreach \j in {1,...,\k}{
        \pgfmathtruncatemacro{\LowerLimit}{\j-1}
        \foreach \i in {\LowerLimit,...,\Nmt}{
            \pgfmathtruncatemacro{\Nexti}{\i+1}
            \draw[<-] (\i_\j) -- (\Nexti_\j);
        }
    }

    \foreach \i in {1,...,\Nmt}{
        \pgfmathtruncatemacro{\LowerLimit}{\i}
        \foreach \j in {\LowerLimit,...,\Nmt}{
            \pgfmathtruncatemacro{\Nextj}{\j+1}
            \pgfmathtruncatemacro{\Newi} {\i-1}
            \draw[->] (\Newi_\j) -- (\Newi_\Nextj);
        }
    }
    
    \foreach \i in {1,...,\k}{
        \pgfmathtruncatemacro{\UpperLimit}{\i-1}
        \foreach \j in {0,...,\UpperLimit}{
            \pgfmathtruncatemacro{\Nextj}{\j+1}
            \pgfmathtruncatemacro{\Newi} {\i-1}
            \draw[<-] (\Newi_\j) -- (\Newi_\Nextj);
        }
    }

    \foreach \i in {2,...,\Nmt}{
    \pgfmathtruncatemacro{\UpperLimit}{\i-1}
        \foreach \m in {1,...,\UpperLimit}{
        	\pgfmathtruncatemacro{\Newi} {\i-1}
            \pgfmathtruncatemacro{\Currx}{\Newi-\m+1}
            \pgfmathtruncatemacro{\Curry}{\N-\m}
            \pgfmathtruncatemacro{\Nextx}{\Newi-\m}
            \pgfmathtruncatemacro{\Nexty}{\N-\m-1}
            \draw[->] (\Currx_\Curry) -- (\Nextx_\Nexty);
        }
    }
    
    \foreach \i in {2,...,\N}{
    \pgfmathtruncatemacro{\UpperLimit}{\N-\i+1}
        \foreach \m in {1,...,\UpperLimit}{
        	\pgfmathtruncatemacro{\Newi} {\i-1}
            \pgfmathtruncatemacro{\Currx}{\Newi+\m-2}
            \pgfmathtruncatemacro{\Curry}{\m-1}
            \pgfmathtruncatemacro{\Nextx}{\Newi+\m-1}
            \pgfmathtruncatemacro{\Nexty}{\m}
            \draw[->] (\Currx_\Curry) -- (\Nextx_\Nexty);
        }
    }
        

\foreach \j in {1,...,\k}{
\pgfmathtruncatemacro{\lm}{\k}

\draw[<-] (0_\j) to [out=160,in=20] (\lm_\j);
}

\foreach \j in {2,...,\k}{
\pgfmathtruncatemacro{\Prevj}{\j-1}
\pgfmathtruncatemacro{\lm}{\k}
\draw[->] (0_\j) to  (\lm_\Prevj);
}

\foreach \i in {1,...,\k}{
        	\pgfmathtruncatemacro{\Newi} {\i-1}
\draw[->] (\Newi_\k) to [out=75,in=-75] (\Newi_0);
}

\foreach \i in {1,...,\Nmt}{
\pgfmathtruncatemacro{\Nexti}{\i+1}
   \pgfmathtruncatemacro{\Newi} {\i-1}
\draw[<-] (\Newi_\k) to (\i_0);
}

\pgfmathtruncatemacro{\lm}{\k}
\draw[->] (0_0) to [out=90,in=-110] (\lm_\lm);    


\node[draw,circle,thick,ForestGreen] (A) at (\k*\dx-\dx/2,\N*\dy ) { $x_A$};

\node[draw,circle,thick, blue] (B) at (-\dx,\dy/2 ) {$x_B$};


\pgfmathtruncatemacro{\lol}{\N} 
\pgfmathtruncatemacro{\wut}{\k}   
\draw[<-]  (4_4) -- (A);
\draw[->]  (3_4) -- (A);

\draw[->]  (0_0) to  (B);
\draw[->]  (B) to (0_1);

\end{tikzpicture}
\caption{The initial seed $\mathcal{Q}$ and its tropical variables}
\label{fig:baxquiv-init}
\end{figure}
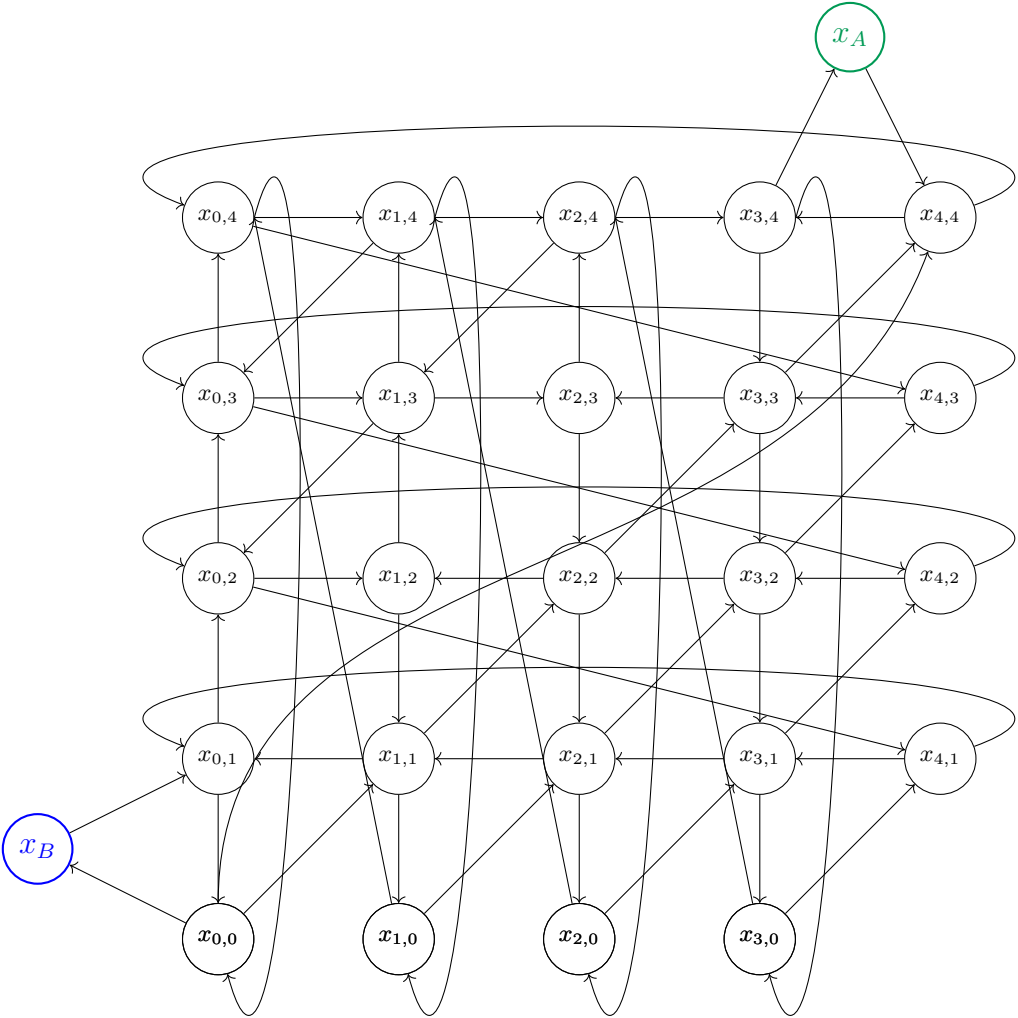

We now explain how this argument can be extended to establish the rank $N$ specialization of the pentagon identity of Conjecture~\ref{conj: pentagon relation for punctured torus} for the Baxter operators in the skein of the punctured torus $T^2\setminus D$. In this case, the relevant moduli space is the Fock-Goncharov cluster variety $\cX_{GL_N,\Sigma_{1,1}}$ parametrizing $G=GL_N$-local systems on a punctured torus, with the additional data of a flat section of the associated $G/B$-bundle in a neighborhood the puncture. As explained in~\cite{FG1}, each ideal triangulation of $\Sigma_{1,1}$ determines a cluster chart on  $\cX_{GL_N,\Sigma_{1,1}}$, whose quiver is obtained by subdividing each triangle into an $N$-triangulation, and then amalgamating the corresponding quivers along the edges of the ideal triangulation. Strictly speaking, the construction from~\cite{FG1} applies to the group $G=PGL_N$; the passage to $GL_N$ is achieved by introducing two frozen $\cX$-variables $x_A,x_B$ which we use to parametrize the determinant of the holonomy around the $(0,1)$ and $(1,0)$-curves on the torus.
The tokamak-shaped quiver obtained by applying this construction to  $\Sigma_{1,1}$ for $N=5$ is illustrated in Figure~\ref{fig:baxquiv-init}.

We now show how to realize the Baxter operators $Q_{(1,0)},Q_{(0,1)},Q_{(1,1)}$ as sequences of mutations in the cluster algebra by `defrosting' the variables $x_A,x_B$.  Just as for $N=2$, this requires understanding how to express the $k$-th fundamental quantum trace $e_{k\mathbf{x}}$ of the monodromy around along the fundamental cycle of such annulus in cluster coordinates. 

\begin{figure}[ht]

\qquad
\begin{tikzpicture}[scale=.4]

	\pgfmathsetmacro{\N}{4}
	\pgfmathsetmacro{\k}{\N-1}
	\pgfmathsetmacro{\Nmt}{\N-2}
    \pgfmathsetmacro{\dx}{2.5}
    \pgfmathsetmacro{\dy}{2.5}
    \pgfmathsetmacro{\eps}{0.2}
    \pgfmathsetmacro{\varoffset}{0.05}
     \pgfmathsetmacro{\baxoffset}{3}


\node[] (s_1) at (0*\dx,1*\dy ) {};
\node[] (s_2) at (0*\dx,2*\dy ) {};

\node[] (t_1) at (6*\dx,1*\dy ) {};
\node[] (t_2) at (6*\dx,2*\dy ) {};

\node[draw,circle,fill=white, inner sep = .08cm] (w_1) at (2*\dx,1*\dy ) {};

\node[draw,circle,fill=white, inner sep = .08cm] (w_2) at (4*\dx,2*\dy ) {};

\node[draw,circle,fill=black, inner sep = .08cm] (b_1) at (4*\dx,1*\dy ) {};

\node[draw,circle,fill=black, inner sep = .08cm] (b_2) at (2*\dx,2*\dy ) {};


\node[] at (3*\dx,2.5*\dy) {$X_{e_5}$};

\node[] at (\dx,1.5*\dy) {$X_{e_2}$};

\node[] at (3*\dx,1.5*\dy) {$X_{e_3}$};

\node[] at (5*\dx,1.5*\dy) {$X_{e_2}$};

\draw[->] (s_1) -- (w_1);

\draw[->] (s_2) -- (b_2);

\draw[->] (b_1) -- (t_1);

\draw[->] (w_2) -- (t_2);

\draw[->] (w_1) -- (b_1);

\draw[->] (b_2) -- (w_2);


\draw[->] (b_2) -- (w_1);


\draw[->] (b_1) -- (w_2);

\end{tikzpicture}
\caption{
Directed graph computing the transport along $(0,1)$-curve for $N=2$; the left and right boundaries in the graph are identified.}
\label{fig:trinetworkN2}
\end{figure}
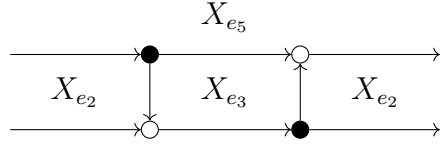

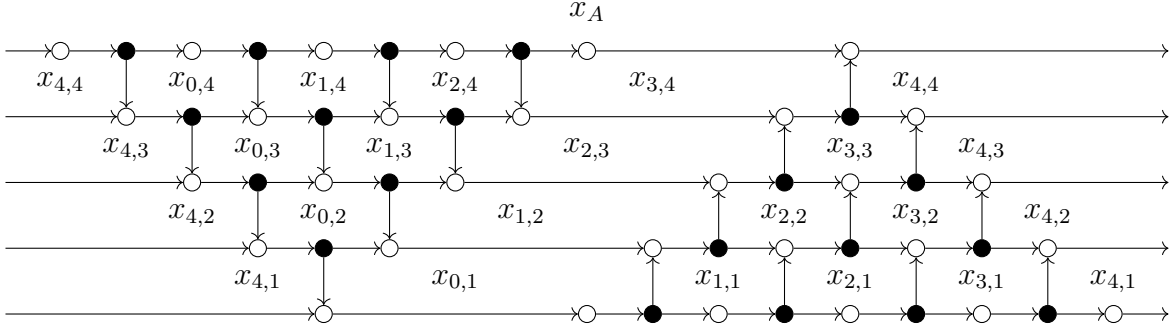
\begin{figure}[ht]

\qquad
\begin{tikzpicture}[scale=.35]

	\pgfmathsetmacro{\N}{5}
	\pgfmathsetmacro{\k}{\N-1}
	\pgfmathsetmacro{\Nmt}{\N-2}
    \pgfmathsetmacro{\dx}{2.5}
    \pgfmathsetmacro{\dy}{2.5}
    \pgfmathsetmacro{\eps}{0.2}
    \pgfmathsetmacro{\varoffset}{0.05}
     \pgfmathsetmacro{\baxoffset}{3}

\node[] at (5*\dx,5.6*\dy) {$x_{A}$};

    \foreach \j in {1,...,\N}{
        \foreach \i in {1,...,\j}{
			\node[] (s_\j) at (-\N*\dx+\dx,\j*\dy ) {};
			\node[] (t_\j) at (2*\N*\dx+2*\dx+2*\dx,\j*\dy ) {};
        	}
        }

    \foreach \j in {1,...,\N}{
        \foreach \i in {1,...,\j}{
			\node[draw,circle,fill=white, inner sep = .08cm] (w_\j_\i) at (-\j*\dx+2*\i*\dx,\j*\dy ) {};
        	}
        }

    \foreach \j in {3,...,\N}{
        \pgfmathtruncatemacro{\limit}{\j-1}

        \foreach \i in {2,...,\limit}{
                \pgfmathtruncatemacro{\im}{\i-2}
			\node (dl_\j_\i) at (-\j*\dx+2*\i*\dx,\j*\dy-.5*\dy ) {$x_{\im,\limit}$};
        	}
        }        

\foreach \j in {2,...,\N}{
 \pgfmathtruncatemacro{\cc}{\N-1}
  \pgfmathtruncatemacro{\cd}{\j-1}
\node at (2*\dx-\j*\dx,\j*\dy-.5*\dy) {$x_{\cc,\cd}$};
}


\foreach \j in {2,...,\N}{
 \pgfmathtruncatemacro{\cc}{\j-1}
  \pgfmathtruncatemacro{\cd}{\j-2}

\node (ml_\j) at (\j*\dx+\dx,\j*\dy-.5*\dy) {$x_{\cd,\cc}$};
}

    \foreach \j in {2,...,\N}{
        \pgfmathtruncatemacro{\Upperlimit}{\j-1}
        \foreach \i in {1,...,\Upperlimit}{
			\node[draw,circle,fill=black,inner sep=.08cm] (b_\j_\i) at (-\j*\dx+\dx+2*\i*\dx,\j*\dy ) {};
        	}
        }


    \foreach \j in {1,...,\N}{
        \pgfmathtruncatemacro{\Limit}{\N-\j+1}
        \foreach \i in {1,...,\Limit}{
			\node[draw,circle,fill=white, inner sep = .08cm] (wu_\j_\i) at (-3*\dx+\N*\dx+2*\i*\dx+\j*\dx,\j*\dy ) {};
        	}
        }
        
    \foreach \j in {1,...,\k}{
        \pgfmathtruncatemacro{\Limit}{\N-\j}
        \foreach \i in {1,...,\Limit}{
			\node[draw,circle,fill=black, inner sep = .08cm] (bu_\j_\i) at (-2*\dx+\N*\dx+2*\i*\dx+\j*\dx,\j*\dy ) {};
        	}
        }

    \foreach \j in {1,...,\k}{
        \pgfmathtruncatemacro{\Limit}{\N-\j+1}
        \foreach \i in {2,...,\Limit}{
        \pgfmathtruncatemacro{\im}{\i+\j-2}
			\node (ul_\j_\i) at (-3*\dx+\N*\dx+2*\i*\dx+\j*\dx,\j*\dy+.5*\dy ) {$x_{\im,\j}$};
        	}
        }

\draw[->] (s_1) -- (w_1_1);

\foreach \j in {1,...,\N}{
\pgfmathtruncatemacro{\compj}{\N+1-\j}
\draw[->] (wu_\j_\compj) -- (t_\j);
}

\foreach \j in {1,...,\N}{
\draw[->] (w_\j_\j) -- (wu_\j_1);
}

    \foreach \j in {2,...,\N}{
    	\draw[->] (s_\j) -- (w_\j_1);
        \pgfmathtruncatemacro{\Limit}{\j-1}
        \foreach \i in {1,...,\Limit}{
            \pgfmathtruncatemacro{\Nexti}{\i+1}
            \draw[->] (w_\j_\i) -- (b_\j_\i);
            \draw[->] (b_\j_\i) -- (w_\j_\Nexti);
            \draw[->] (b_\j_\i) -- (w_\Limit_\i);
        }
    }

    \foreach \j in {1,...,\k}{
        \pgfmathtruncatemacro{\Limit}{\N-\j}
        \foreach \i in {1,...,\Limit}{
            \pgfmathtruncatemacro{\Nexti}{\i+1}
            \pgfmathtruncatemacro{\Nextj}{\j+1}
            \draw[->] (wu_\j_\i) -- (bu_\j_\i);
            \draw[->] (bu_\j_\i) -- (wu_\j_\Nexti);
            \draw[->] (bu_\j_\i) -- (wu_\Nextj_\i);
        }
    }


\end{tikzpicture}
\caption{
Initial black-white graph for computing transport along the $(1,0)$-curve for $N=5$.}
\label{fig:trinetworkN5}
\end{figure}

Note that in the ideal triangulation in Figure~\ref{fig:baxquiv-init}, each of the $(1,0),(1,0)$ and $(1,1)$-curves is isolated inside an ideal annulus obtained by identifying two opposite sides of a triangulated ideal quadrilateral. Following~\cite{SS17}, the $k$-th fundamental quantum trace $e_{k\mathbf{x}}$ of the monodromy around along the fundamental cycle of such annulus can be computed as a partition function of non-intersecting $k$-multiloops in a directed black-white graph dual to the quiver. Each such multiloop $\mathbf{\ell}=(\ell_1,\ldots, \ell_k)$ contributes a cluster monomial $X_{\sum_i\sum_{f\subset \ell_i}}$, where the inner sum is taken over all faces `inside' (or `above', as drawn in Figures~\ref{fig:trinetworkN2} and~\ref{fig:trinetworkN5}) the path $\ell_i$.

To illustrate the definition, and compare with the discussion in the previous section, we draw this directed black-white graph for $N=2$ in Figure~\ref{fig:trinetworkN2}. 

Identifying left and right pairs of terminal vertices, there are exactly three closed loops in this graph (two based at the top left terminal vertex, and one based at the bottom left terminal vertex), corresponding to the three summands in the formula for $P_{(0,1)}=e_{(0,1)}$ in~\ref{eq:cluster-embed1}.  As in the previous section the face $X_{e_5}$ `at infinity' is regarded as a frozen $\cX$-variable with respect to the cluster structure on $\cX_{GL_2,\Sigma_{1,1}}$. As there is only one non-intersecting 2-path in Figure~\ref{fig:trinetworkN2}, the determinant of the holonomy around the $(0,1)$-curve is given by $X_{e_2+e_3 + 2e_5}$.

For general $N$, the graph computing the transport along the $(0,1)$-curve is dual to the $N$-triangulation of the punctured torus, and takes the form shown in Figure~\ref{fig:trinetworkN5} for $N=5$. 

By Proposition {4.2} of~\cite{SS17}, this combinatorial formula for $e_{k\mathbf{x}}$ is invariant under mutations corresponding to square moves on the black-white graph. Thus to compute (for example) $Q_{(0,1)}$, we have the freedom to pass to a square-move-equivalent cluster in which the expression for $e_{(0,k)}$ in cluster coordinates becomes as simple as possible. 

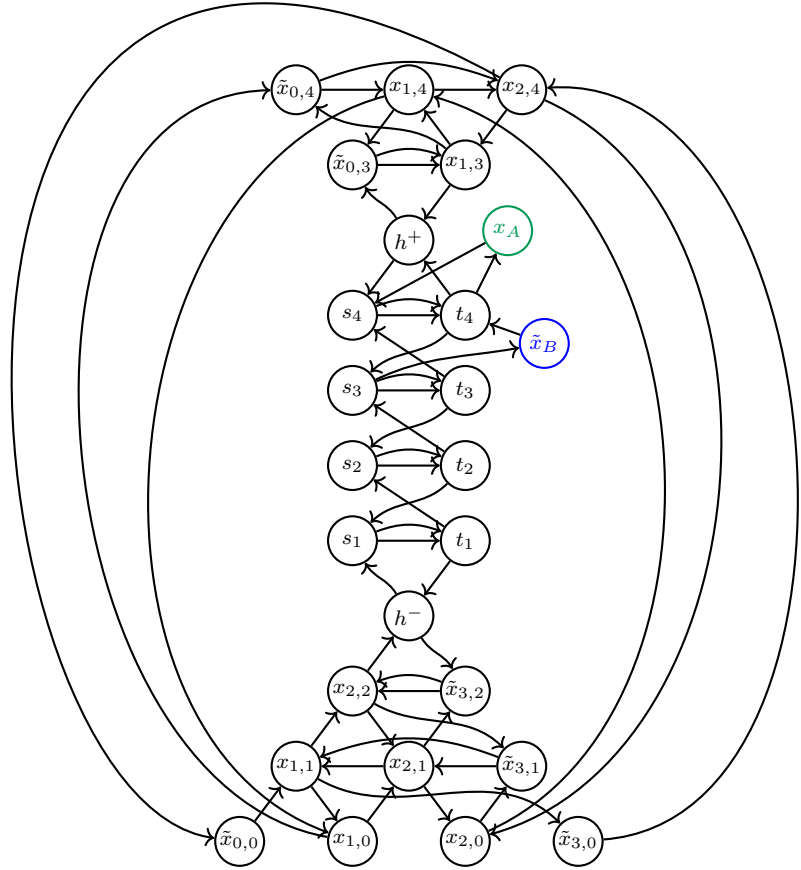
\begin{figure}[ht]
\begin{tikzpicture}[every node/.style={inner sep=0, minimum size=0.65cm, thick, draw, circle}, x=0.75cm, y=0.5cm]

\node (1) at (-2,0) {\tiny{$x_{1,1}$}};
\node (2) at (0,0) {\tiny{$x_{2,1}$}};
\node (3) at (2,0) {\tiny{$\tilde x_{3,1}$}};
\node (4) at (-1,2) {\tiny{$x_{2,2}$}};
\node (5) at (1,2) {\tiny{$\tilde x_{3,2}$}};
\node (6) at (0,4) {\tiny{$h^-$}};
\node (7) at (-1,6) {\tiny{$s_1$}};
\node (8) at (1,6) {\tiny{$t_1$}};
\node (9) at (-1,8) {\tiny{$s_2$}};
\node (10) at (1,8) {\tiny{$t_2$}};
\node (11) at (-1,10) {\tiny{$s_3$}};
\node (12) at (1,10) {\tiny{$t_3$}};
\node (13) at (-1,12) {\tiny{$s_4$}};
\node (14) at (1,12) {\tiny{$t_4$}};
\node (15) at (0,14) {\tiny{$h^+$}};
\node (16) at (-1,16) {\tiny{$\tilde x_{0,3}$}};
\node (17) at (1,16) {\tiny{$x_{1,3}$}};
\node (18) at (-2,18) {\tiny{$\tilde x_{0,4}$}};
\node (19) at (0,18) {\tiny{$x_{1,4}$}};
\node (20) at (2,18) {\tiny{$x_{2,4}$}};

\node (21) at (-3,-2) {\tiny{$\tilde x_{0,0}$}};
\node (22) at (-1,-2) {\tiny{$ x_{1,0}$}};
\node (23) at (1,-2) {\tiny{$x_{2,0}$}};
\node (24) at (3,-2) {\tiny{$\tilde x_{3,0}$}};

\node[ForestGreen] (A) at (1.75, 14.25) {\tiny $x_A$};
\node[blue] (B) at (2.4, 11.25) {\tiny $\tilde x_B$};
\draw [->, thick] (14) -- (A);
\draw [->, thick] (A) -- (13);
\draw [->, thick] (11) to [out=25,in=190] (B);
\draw [->, thick] (B) -- (14);

\draw [->, thick] (21) -- (1);
\draw [->, thick] (1) -- (22);
\draw [->, thick] (22) -- (2);
\draw [->, thick] (2) -- (23);
\draw [->, thick] (23) -- (3);
\draw [->, thick] (1) to [out = -30, in=135] (24);
\draw [->, thick] (24) to [bend right=90] (20);
\draw [->, thick] (20) .. controls +(-12, 10) and +(-5, 1) .. (21);
\draw [->, thick] (19) to [out = -165 , in=160] (22);
\draw [->, thick] (23) to [bend right=70] (19);
\draw [->, thick] (20) to [bend left=70] (23);
\draw [->, thick] (22) to [out = 170 , in=-179] (18);

\draw [->, thick] (3) -- (2);
\draw [->, thick] (2) -- (1);
\draw [->, thick] (1) -- (4);
\draw [->, thick] (4) -- (2);
\draw [->, thick] (2) -- (5);
\draw [->, thick] (5) -- (4);
\draw [->, thick] (4) -- (6);

\draw [->, thick] (6) to [out = -60, in=120] (5);
\draw [->, thick] (5) to [bend right = 21] (4);
\draw [->, thick] (4) to [out = -30, in=135] (3);
\draw [->, thick] (3) to [bend right = 21] (1);

\draw [->, thick] (8) -- (6);
\draw [->, thick] (6) to [out = 120, in=-60] (7);
\draw [->, thick] (7) -- (8);
\draw [->, thick] (8) -- (9);
\draw [->, thick] (9) -- (10);
\draw [->, thick] (10) -- (11);
\draw [->, thick] (11) -- (12);
\draw [->, thick] (12) -- (13);
\draw [->, thick] (13) -- (14);
\draw [->, thick] (14) -- (15);
\draw [->, thick] (15) -- (13);
\draw [->, thick] (13) to [bend left = 21] (14);
\draw [->, thick] (14) to [out = -135, in=45] (11);
\draw [->, thick] (11) to [bend left = 21] (12);
\draw [->, thick] (12) to [out = -135, in=45] (9);
\draw [->, thick] (9) to [bend left = 21] (10);
\draw [->, thick] (10) to [out = -135, in=45] (7);
\draw [->, thick] (7) to [bend left = 21] (8);

\draw [->, thick] (18) -- (19);
\draw [->, thick] (19) -- (20);
\draw [->, thick] (20) -- (17);
\draw [->, thick] (17) -- (19);
\draw [->, thick] (19) -- (16);
\draw [->, thick] (16) -- (17);
\draw [->, thick] (17) -- (15);

\draw [->, thick] (15) to [out = 120, in=-60] (16);
\draw [->, thick] (16) to [bend left = 21] (17);
\draw [->, thick] (17) to [out = 140, in=-40] (18);
\draw [->, thick] (18) to [bend left = 21] (20);

\end{tikzpicture}
\caption{Quiver $\sigma_{(0,1)}(\mathcal{Q})$ obtained by applying the mutation sequence $\sigma_{(0,1)}$ to the initial quiver $\mathcal{Q}$.}
\label{fig:coxeter-quiver}
\end{figure}

\begin{figure}[ht]

\qquad
\begin{tikzpicture}[scale=.4]

	\pgfmathsetmacro{\N}{4}
	\pgfmathsetmacro{\k}{\N-1}
	\pgfmathsetmacro{\Nmt}{\N-2}
    \pgfmathsetmacro{\dx}{2.5}
    \pgfmathsetmacro{\dy}{2.5}
    \pgfmathsetmacro{\eps}{0.2}
    \pgfmathsetmacro{\varoffset}{0.05}
     \pgfmathsetmacro{\baxoffset}{3}


\node[] at (3*\dx,5.5*\dy) {$x_{A}$};

\node[] (s_1) at (0*\dx,1*\dy ) {};
\node[] (s_2) at (0*\dx,2*\dy ) {};
\node[] (s_3) at (0*\dx,3*\dy ) {};
\node[] (s_4) at (0*\dx,4*\dy ) {};
\node[] (s_5) at (0*\dx,5*\dy ) {};

\node[] (t_1) at (6*\dx,1*\dy ) {};
\node[] (t_2) at (6*\dx,2*\dy ) {};
\node[] (t_3) at (6*\dx,3*\dy ) {};
\node[] (t_4) at (6*\dx,4*\dy ) {};
\node[] (t_5) at (6*\dx,5*\dy ) {};

\node[draw,circle,fill=white, inner sep = .08cm] (w_1) at (2*\dx,1*\dy ) {};
\node[draw,circle,fill=white, inner sep = .08cm] (w_3) at (2*\dx,3*\dy ) {};
\node[draw,circle,fill=white, inner sep = .08cm] (w_5) at (2*\dx,5*\dy ) {};

\node[draw,circle,fill=white, inner sep = .08cm] (w_2) at (4*\dx,2*\dy ) {};
\node[draw,circle,fill=white, inner sep = .08cm] (w_4) at (4*\dx,4*\dy ) {};

\node[draw,circle,fill=black, inner sep = .08cm] (b_1) at (4*\dx,1*\dy ) {};
\node[draw,circle,fill=black, inner sep = .08cm] (b_3) at (4*\dx,3*\dy ) {};
\node[draw,circle,fill=black, inner sep = .08cm] (b_5) at (4*\dx,5*\dy ) {};

\node[draw,circle,fill=black, inner sep = .08cm] (b_2) at (2*\dx,2*\dy ) {};
\node[draw,circle,fill=black, inner sep = .08cm] (b_4) at (2*\dx,4*\dy ) {};


\node[] at (\dx,1.5*\dy) {$s_{1}$};
\node[] at (\dx,2.5*\dy) {$t_{2}$};
\node[] at (\dx,3.5*\dy) {$s_3$};
\node[] at (\dx,4.5*\dy) {$t_{4}$};

\node[] at (3*\dx,1.5*\dy) {$t_{1}$};
\node[] at (3*\dx,2.5*\dy) {$s_2$};
\node[] at (3*\dx,3.5*\dy) {$t_{3}$};
\node[] at (3*\dx,4.5*\dy) {$s_4$};

\node[] at (5*\dx,1.5*\dy) {$s_1$};
\node[] at (5*\dx,2.5*\dy) {$t_{2}$};
\node[] at (5*\dx,3.5*\dy) {$s_3$};
\node[] at (5*\dx,4.5*\dy) {$t_{4}$};

\draw[->] (s_1) -- (w_1);
\draw[->] (s_3) -- (w_3);
\draw[->] (s_5) -- (w_5);

\draw[->] (s_2) -- (b_2);
\draw[->] (s_4) -- (b_4);

\draw[->] (b_1) -- (t_1);
\draw[->] (b_3) -- (t_3);
\draw[->] (b_5) -- (t_5);

\draw[->] (w_2) -- (t_2);
\draw[->] (w_4) -- (t_4);

\draw[->] (w_1) -- (b_1);
\draw[->] (w_3) -- (b_3);
\draw[->] (w_5) -- (b_5);

\draw[->] (b_2) -- (w_2);
\draw[->] (b_4) -- (w_4);


\draw[->] (b_2) -- (w_1);
\draw[->] (b_2) -- (w_3);

\draw[->] (b_4) -- (w_5);
\draw[->] (b_4) -- (w_3);

\draw[->] (b_1) -- (w_2);
\draw[->] (b_3) -- (w_4);
\draw[->] (b_3) -- (w_2);
\draw[->] (b_5) -- (w_4);

\end{tikzpicture}
\caption{
Network computing the transport along $(0,1)$-curve in the cluster $\sigma_{(1,0)}(\Pi_0)$ given by the quiver in Figure~\ref{fig:coxeter-quiver}}
\label{fig:coxeter-network}
\end{figure}
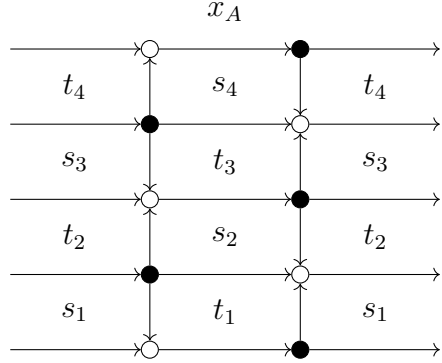

In particular, we can apply the mutation sequence denoted $\phi_2\circ\phi_1$ defined in Section 5.6 of~\cite{SS17}, after which the subquiver on all nodes except that $\tilde x_B$ is isomorphic to one shown in Figure~\ref{fig:coxeter-quiver}. In fact, it is convenient to post-compose this sequence with the further mutations $\tau_-^{N-2}$, where $\tau_-$ is the Dehn twist mutation sequence obtained by mutating at the sinks of all double arrows in the quiver~\ref{fig:coxeter-quiver}. We write $\sigma_{(0,1)}=\tau_-^{N-2}\circ\phi^{(0,1)}_2\circ\phi^{(0,1)}_1$ for the corresponding sequence of mutations. 
In the figure, we  label the vertices by their $c$-vectors/tropical $\mathscr{X}$-coordinates in the cluster $\sigma_{(0,1)}(\Pi_0)$, abbreviating
 \begin{align*}
h^+&= x_{0,2}\prod_{j=2}^{N-1}\prod_{i=0}^{N-1}x_{i,j}^{-1}, \quad h^- = x_{N-2,N-2}\prod_{j=1}^{N-2}\prod_{i=0}^{N-1}x_{i,j}^{-1}
 \end{align*}
 \begin{align*}
s_j = x_{N-1,j}, \quad t_j = \prod_{i=0}^{N-2}x_{i,j}
 \end{align*}
 \begin{align*}
 \tilde x_B = x_B \prod_{j=1}^{N-2}\prod_{i=0}^{j-1} x_{i,j}.
 \end{align*}
 Note that both $h^\pm$ are of negative tropical sign. Explicit formulas for the remaining tropical variables $\tilde x_{i,j}$ can easily be written. As we will shortly see, they will not be needed in the computations to follow and so we suppress them here.
The directed graph computing the $e_{(0,k)}$ in terms of the cluster variables from $\sigma_{(0,1)}(\Pi_0)$ is dual to the quiver $\sigma_{(0,1)}(\mathcal{Q})$, and is illustrated in Figure~\ref{fig:coxeter-network} for $N=5$.

By the results of~\cite{SS17}, it follows that in this cluster the quantum holonomy $e_{(0,k)}$ is identified with the $k$-th $\mathfrak{gl}_N$ $q$-difference Toda chain Hamiltonian.  
As shown in~\cite{SS18}, the solution to the difference equation~\ref{eq:bax-qde} is given by the automorphism part of the sequence of $2N-1$ mutations
\begin{align}
\label{eq:bax-muts}
Q^{\mathrm{isolated}}_{(0,1)} = \mu_{s_{1}}\circ\mu_{t_{1}}\circ \cdots\circ\mu_{s_{N-1}}\circ\mu_{t_{N-1}}\circ\mu_{A}.
\end{align}
The quiver $Q^{\mathrm{isolated}}_{(0,1)}\circ\sigma_{(0,1)}(\mathcal{Q})$ is shown in Figure~\ref{fig:coxeter-after-bax} for $N=5$.
{The key observation is that the tropical $\cX$-variables in $Q^{\mathrm{isolated}}_{(0,1)}\circ\sigma_{(0,1)}(\mathcal{Q})$ are \emph{identical} to those of $\sigma_{(0,1)}\circ(\mathcal{Q})$, except for the the determinantal variables $(x_A,\tilde x_B)$ and $(\tilde x_A,\widetilde{\widetilde x}_B)$}, where we have
$$
{\widetilde{x}}_A =x_A^{-1}\prod_{i=0}^{N-1}\prod_{j=1}^{N-1}x_{i,j}^{-1} ,\quad \widetilde{\widetilde{x}}_B = 
x_Ax_B\prod_{j=1}^{N-1}\prod_{i=0}^{j-1}x_{i,j}
$$
Hence setting
$$
\mu_{(0,1)} = \sigma_{(0,1)}^{-1}\circ Q^{\mathrm{isolated}}_{(0,1)}\circ \sigma_{(0,1)}
$$
it follows that the same is true of the sets of tropical variables for the seeds $\mathcal{Q}$ and $\mu_{(0,1)}(\mathcal{Q})$. 
The resulting seed $\mu_{(0,1)}(\mathcal{Q})$, and its tropical $\mathscr{X}$-variables/$c$-vectors, are shown in Figure~\ref{fig:quiver-after-10}, where we have abbreviated
\begin{align}
\label{eq:atilde}
x'_A=&x_A^{-1}\prod_{i=0}^{N-1}\prod_{j=1}^{N-1}x_{i,j}^{-1}\\
\nonumber x_{A+B}&=x_Ax_B\prod_{j=1}^{N-1}\prod_{i=0}^{j-1}x_{i,j}.
\end{align}
\begin{figure}[ht]
\begin{tikzpicture}[every node/.style={inner sep=0, minimum size=0.65cm, thick, draw, circle}, x=0.75cm, y=0.5cm]

\node (1) at (-2,0) {\tiny{$x_{1,1}$}};
\node (2) at (0,0) {\tiny{$x_{2,1}$}};
\node (3) at (2,0) {\tiny{$\tilde x_{3,1}$}};
\node (4) at (-1,2) {\tiny{$x_{2,2}$}};
\node (5) at (1,2) {\tiny{$\tilde x_{3,2}$}};
\node (6) at (0,4) {\tiny{$h^-$}};
\node (7) at (-1,6) {\tiny{$s_1$}};
\node (8) at (1,6) {\tiny{$t_1$}};
\node (9) at (-1,8) {\tiny{$s_2$}};
\node (10) at (1,8) {\tiny{$t_2$}};
\node (11) at (-1,10) {\tiny{$s_3$}};
\node (12) at (1,10) {\tiny{$t_3$}};
\node (13) at (-1,12) {\tiny{$s_4$}};
\node (14) at (1,12) {\tiny{$t_4$}};
\node (15) at (0,14) {\tiny{$h^+$}};
\node (16) at (-1,16) {\tiny{$\tilde x_{0,3}$}};
\node (17) at (1,16) {\tiny{$x_{1,3}$}};
\node (18) at (-2,18) {\tiny{$\tilde x_{0,4}$}};
\node (19) at (0,18) {\tiny{$x_{1,4}$}};
\node (20) at (2,18) {\tiny{$x_{2,4}$}};

\node (21) at (-3,-2) {\tiny{$\tilde x_{0,0}$}};
\node (22) at (-1,-2) {\tiny{$ x_{1,0}$}};
\node (23) at (1,-2) {\tiny{$x_{2,0}$}};
\node (24) at (3,-2) {\tiny{$\tilde x_{3,0}$}};

\node[ForestGreen] (A) at (2, 4.00) {\tiny ${\widetilde{x}}_A$};
\node[blue] (B) at (-2.3, 14.4) {\tiny $\widetilde{\widetilde{x}}_B$};
\draw [->, thick] (7) -- (A);
\draw [->, thick] (A) -- (8);
\draw [->, thick] (13) to  (B);
\draw [->, thick] (B) -- (15);

\draw [->, thick] (21) -- (1);
\draw [->, thick] (1) -- (22);
\draw [->, thick] (22) -- (2);
\draw [->, thick] (2) -- (23);
\draw [->, thick] (23) -- (3);
\draw [->, thick] (1) to [out = -30, in=135] (24);
\draw [->, thick] (24) to [bend right=90] (20);
\draw [->, thick] (20) .. controls +(-12, 10) and +(-5, 1) .. (21);
\draw [->, thick] (19) to [out = -165 , in=160] (22);
\draw [->, thick] (23) to [bend right=70] (19);
\draw [->, thick] (20) to [bend left=70] (23);
\draw [->, thick] (22) to [out = 170 , in=-179] (18);

\draw [->, thick] (3) -- (2);
\draw [->, thick] (2) -- (1);
\draw [->, thick] (1) -- (4);
\draw [->, thick] (4) -- (2);
\draw [->, thick] (2) -- (5);
\draw [->, thick] (5) -- (4);
\draw [->, thick] (4) -- (6);

\draw [->, thick] (6) to [out = -60, in=120] (5);
\draw [->, thick] (5) to [bend right = 21] (4);
\draw [->, thick] (4) to [out = -30, in=135] (3);
\draw [->, thick] (3) to [bend right = 21] (1);

\draw [->, thick] (8) -- (6);
\draw [->, thick] (6) to [out = 120, in=-60] (7);
\draw [->, thick] (7) -- (8);
\draw [->, thick] (8) -- (9);
\draw [->, thick] (9) -- (10);
\draw [->, thick] (10) -- (11);
\draw [->, thick] (11) -- (12);
\draw [->, thick] (12) -- (13);
\draw [->, thick] (13) -- (14);
\draw [->, thick] (14) -- (15);
\draw [->, thick] (15) -- (13);
\draw [->, thick] (13) to [bend left = 21] (14);
\draw [->, thick] (14) to [out = -135, in=45] (11);
\draw [->, thick] (11) to [bend left = 21] (12);
\draw [->, thick] (12) to [out = -135, in=45] (9);
\draw [->, thick] (9) to [bend left = 21] (10);
\draw [->, thick] (10) to [out = -135, in=45] (7);
\draw [->, thick] (7) to [bend left = 21] (8);

\draw [->, thick] (18) -- (19);
\draw [->, thick] (19) -- (20);
\draw [->, thick] (20) -- (17);
\draw [->, thick] (17) -- (19);
\draw [->, thick] (19) -- (16);
\draw [->, thick] (16) -- (17);
\draw [->, thick] (17) -- (15);

\draw [->, thick] (15) to [out = 120, in=-60] (16);
\draw [->, thick] (16) to [bend left = 21] (17);
\draw [->, thick] (17) to [out = 140, in=-40] (18);
\draw [->, thick] (18) to [bend left = 21] (20);

\end{tikzpicture}

\caption{The result of applying the Baxter operator  $Q^{\mathrm{isolated}}_{(0,1)}$. }
\label{fig:coxeter-after-bax}
\end{figure}
So with respect to the coordinates of this initial cluster, the Baxter operator  $Q_{(0,1)}$ coincides with the automorphism part of the sequence of sign-coherent mutations
$$
\mu_{(0,1)} = \sigma_{(0,1)}^{-1}\circ Q^{\mathrm{isolated}}_{(0,1)}\circ \sigma_{(0,1)}.
$$

The same considerations apply to the isolated $(1,0)$ and $(1,1)$--curves, and we obtain analogous realizations for the corresponding Baxter operators. 

\begin{figure}[ht]
\begin{tikzpicture}[scale=.6]

	\pgfmathsetmacro{\N}{5}
	\pgfmathsetmacro{\k}{\N-1}
	\pgfmathsetmacro{\Nmt}{\N-2}
    \pgfmathsetmacro{\dx}{4}
    \pgfmathsetmacro{\dy}{4}
    \pgfmathsetmacro{\eps}{0.2}
    \pgfmathsetmacro{\varoffset}{0.05}
     \pgfmathsetmacro{\baxoffset}{3}


\node[draw,circle] (1_1) at (1*\dx,1*\dy ) { \scriptsize${x_{1,1}}$};
\node[draw,circle] (2_2) at (2*\dx,2*\dy ) {\scriptsize ${x_{2,2}}$};
\node[draw,circle] (3_3) at (3*\dx,3*\dy ) { \scriptsize${x_{3,3}}$};
\node[draw,circle] (4_4) at (4*\dx,4*\dy ) { \scriptsize${x_{4,4}}$};
\node[draw,circle] (2_1) at (2*\dx,1*\dy ) {\scriptsize ${x_{2,1}}$};
\node[draw,circle] (3_2) at (3*\dx,2*\dy ) {\scriptsize ${x_{3,2}}$};
\node[draw,circle] (4_3) at (4*\dx,3*\dy ) { \scriptsize${x_{4,3}}$};
\node[draw,circle] (0_4) at (0*\dx,4*\dy ) {\scriptsize ${x_{0,4}}$};
\node[draw,circle] (3_1) at (3*\dx,1*\dy ) {\scriptsize ${x_{3,1}}$};
\node[draw,circle] (4_2) at (4*\dx,2*\dy ) {\scriptsize ${x_{4,2}}$};
\node[draw,circle] (0_3) at (0*\dx,3*\dy ) { \scriptsize${x_{0,3}}$};
\node[draw,circle] (1_4) at (1*\dx,4*\dy ) {\scriptsize ${x_{1,4}}$};
\node[draw,circle] (4_1) at (4*\dx,1*\dy ) {\scriptsize ${x_{4,1}}$};
\node[draw,circle] (0_2) at (0*\dx,2*\dy ) {\scriptsize ${x_{0,2}}$};
\node[draw,circle] (1_3) at (1*\dx,3*\dy ) { \scriptsize${x_{1,3}}$};
\node[draw,circle] (2_4) at (2*\dx,4*\dy ) {\scriptsize ${x_{2,4}}$};
\node[draw,circle] (0_1) at (0*\dx,1*\dy ) {\scriptsize ${x_{0,1}}$};
\node[draw,circle] (1_2) at (1*\dx,2*\dy ) {\scriptsize ${x_{1,2}}$};
\node[draw,circle] (2_3) at (2*\dx,3*\dy ) { \scriptsize${x_{2,3}}$};
\node[draw,circle] (3_4) at (3*\dx,4*\dy ) {\scriptsize ${x_{3,4}}$};
\node[draw,circle] (0_0) at (0*\dx,0*\dy ) {\scriptsize ${x_{0,0}}$};
\node[draw,circle] (1_0) at (1*\dx,0*\dy ) {\scriptsize ${x_{1,0}}$};
\node[draw,circle] (2_0) at (2*\dx,0*\dy ) {\scriptsize ${x_{2,0}}$};
\node[draw,circle] (3_0) at (3*\dx,0*\dy ) {\scriptsize ${x_{3,0}}$};
\node[draw,circle] (0_0) at (0*\dx,0*\dy ) {\scriptsize ${x_{0,0}}$};
\node[draw,circle] (1_0) at (1*\dx,0*\dy ) {\scriptsize ${x_{1,0}}$};
\node[draw,circle] (2_0) at (2*\dx,0*\dy ) {\scriptsize ${x_{2,0}}$};
\node[draw,circle] (3_0) at (3*\dx,0*\dy ) {\scriptsize ${x_{3,0}}$};


    \foreach \j in {2,...,\k}{
        \pgfmathtruncatemacro{\UpperLimit}{\j-1}
        \foreach \i in {1,...,\UpperLimit}{
            \pgfmathtruncatemacro{\Nexti}{\i-1}
            \draw[->] (\Nexti_\j) -- (\i_\j);
        }
    }
    
    \foreach \j in {1,...,\k}{
        \pgfmathtruncatemacro{\LowerLimit}{\j-1}
        \foreach \i in {\LowerLimit,...,\Nmt}{
            \pgfmathtruncatemacro{\Nexti}{\i+1}
            \draw[<-] (\i_\j) -- (\Nexti_\j);
        }
    }

    \foreach \i in {1,...,\Nmt}{
        \pgfmathtruncatemacro{\LowerLimit}{\i}
        \foreach \j in {\LowerLimit,...,\Nmt}{
            \pgfmathtruncatemacro{\Nextj}{\j+1}
            \pgfmathtruncatemacro{\Newi} {\i-1}
            \draw[->] (\Newi_\j) -- (\Newi_\Nextj);
        }
    }
    
    \foreach \i in {1,...,\k}{
        \pgfmathtruncatemacro{\UpperLimit}{\i-1}
        \foreach \j in {0,...,\UpperLimit}{
            \pgfmathtruncatemacro{\Nextj}{\j+1}
            \pgfmathtruncatemacro{\Newi} {\i-1}
            \draw[<-] (\Newi_\j) -- (\Newi_\Nextj);
        }
    }

    \foreach \i in {2,...,\Nmt}{
    \pgfmathtruncatemacro{\UpperLimit}{\i-1}
        \foreach \m in {1,...,\UpperLimit}{
        	\pgfmathtruncatemacro{\Newi} {\i-1}
            \pgfmathtruncatemacro{\Currx}{\Newi-\m+1}
            \pgfmathtruncatemacro{\Curry}{\N-\m}
            \pgfmathtruncatemacro{\Nextx}{\Newi-\m}
            \pgfmathtruncatemacro{\Nexty}{\N-\m-1}
            \draw[->] (\Currx_\Curry) -- (\Nextx_\Nexty);
        }
    }
    
    \foreach \i in {2,...,\N}{
    \pgfmathtruncatemacro{\UpperLimit}{\N-\i+1}
        \foreach \m in {1,...,\UpperLimit}{
        	\pgfmathtruncatemacro{\Newi} {\i-1}
            \pgfmathtruncatemacro{\Currx}{\Newi+\m-2}
            \pgfmathtruncatemacro{\Curry}{\m-1}
            \pgfmathtruncatemacro{\Nextx}{\Newi+\m-1}
            \pgfmathtruncatemacro{\Nexty}{\m}
            \draw[->] (\Currx_\Curry) -- (\Nextx_\Nexty);
        }
    }
        

\foreach \j in {1,...,\k}{
\pgfmathtruncatemacro{\lm}{\k}

\draw[<-] (0_\j) to [out=160,in=20] (\lm_\j);
}

\foreach \j in {2,...,\k}{
\pgfmathtruncatemacro{\Prevj}{\j-1}
\pgfmathtruncatemacro{\lm}{\k}
\draw[->] (0_\j) to  (\lm_\Prevj);
}

\foreach \i in {1,...,\k}{
        	\pgfmathtruncatemacro{\Newi} {\i-1}
\draw[->] (\Newi_\k) to [out=75,in=-75] (\Newi_0);
}

\foreach \i in {1,...,\Nmt}{
\pgfmathtruncatemacro{\Nexti}{\i+1}
   \pgfmathtruncatemacro{\Newi} {\i-1}
\draw[<-] (\Newi_\k) to (\i_0);
}

\pgfmathtruncatemacro{\lm}{\k}
\draw[->] (0_0) to [out=90,in=-110] (\lm_\lm);    


\node[draw,circle,minimum size=.75cm] (A) at (6,-4) {\tiny $x'_A$};

\node[draw,circle,minimum size=.75cm] (B) at (-1*\dx,5*\dy ) {\tiny$x_{A+B}$};


\pgfmathtruncatemacro{\lol}{\N} 
\pgfmathtruncatemacro{\wut}{\k}   
\draw[->]  (0_1)  to [out=-175,in=-145] (A);
\draw[->]  (A) to [out=-15,in=-80]  (4_1);

\draw[<-]  (0_0) to [out=160,in=-60]  (B);
\draw[->]  (4_4) to [out=90,in=0] (B) ;

\end{tikzpicture}

\caption{Tropical variables for the seed $\mu_{0,1}(\mathcal{Q})$}
\label{fig:quiver-after-10}
\end{figure}

\begin{figure}[ht]

\qquad
\begin{tikzpicture}[scale=.6]

	\pgfmathsetmacro{\N}{5}
	\pgfmathsetmacro{\k}{\N-1}
	\pgfmathsetmacro{\Nmt}{\N-2}
    \pgfmathsetmacro{\dx}{4}
    \pgfmathsetmacro{\dy}{4}
    \pgfmathsetmacro{\eps}{0.2}
    \pgfmathsetmacro{\varoffset}{0.05}
     \pgfmathsetmacro{\baxoffset}{3}


\node[draw,circle,fill=red] (1_1) at (1*\dx,1*\dy ) { \scriptsize${x_{1,1}}$};
\node[draw,circle,fill=red] (2_2) at (2*\dx,2*\dy ) {\scriptsize ${x_{2,2}}$};
\node[draw,circle,fill=red] (3_3) at (3*\dx,3*\dy ) { \scriptsize${x_{3,3}}$};
\node[draw,circle,fill=red] (4_4) at (4*\dx,4*\dy ) { \scriptsize${x_{4,4}}$};
\node[draw,circle,fill=yellow] (2_1) at (2*\dx,1*\dy ) {\scriptsize ${x_{2,1}}$};
\node[draw,circle,fill=yellow] (3_2) at (3*\dx,2*\dy ) {\scriptsize ${x_{3,2}}$};
\node[draw,circle,fill=yellow] (4_3) at (4*\dx,3*\dy ) { \scriptsize${x_{4,3}}$};
\node[draw,circle,fill=yellow] (0_4) at (0*\dx,4*\dy ) {\scriptsize ${x_{0,4}}$};
\node[draw,circle,fill=green] (3_1) at (3*\dx,1*\dy ) {\scriptsize ${x_{3,1}}$};
\node[draw,circle,,fill=green] (4_2) at (4*\dx,2*\dy ) {\scriptsize ${x_{4,2}}$};
\node[draw,circle,fill=green] (0_3) at (0*\dx,3*\dy ) { \scriptsize${x_{0,3}}$};
\node[draw,circle,fill=green] (1_4) at (1*\dx,4*\dy ) {\scriptsize ${x_{1,4}}$};
\node[draw,circle,fill=cyan] (4_1) at (4*\dx,1*\dy ) {\scriptsize ${x_{4,1}}$};
\node[draw,circle,fill=cyan] (0_2) at (0*\dx,2*\dy ) {\scriptsize ${x_{0,2}}$};
\node[draw,circle,fill=cyan] (1_3) at (1*\dx,3*\dy ) { \scriptsize${x_{1,3}}$};
\node[draw,circle,fill=cyan] (2_4) at (2*\dx,4*\dy ) {\scriptsize ${x_{2,4}}$};
\node[draw,circle] (0_1) at (0*\dx,1*\dy ) {\scriptsize ${x_{0,1}}$};
\node[draw,circle] (1_2) at (1*\dx,2*\dy ) {\scriptsize ${x_{1,2}}$};
\node[draw,circle] (2_3) at (2*\dx,3*\dy ) { \scriptsize${x_{2,3}}$};
\node[draw,circle] (3_4) at (3*\dx,4*\dy ) {\scriptsize ${x_{3,4}}$};
\node[draw,circle,fill=red] (0_0) at (0*\dx,0*\dy ) {\scriptsize ${x_{0,0}}$};
\node[draw,circle,fill=yellow] (1_0) at (1*\dx,0*\dy ) {\scriptsize ${x_{1,0}}$};
\node[draw,circle,fill=green] (2_0) at (2*\dx,0*\dy ) {\scriptsize ${x_{2,0}}$};
\node[draw,circle,fill=cyan] (3_0) at (3*\dx,0*\dy ) {\scriptsize ${x_{3,0}}$};
\node[draw,circle] (0_0) at (0*\dx,0*\dy ) {\scriptsize ${x_{0,0}}$};
\node[draw,circle] (1_0) at (1*\dx,0*\dy ) {\scriptsize ${x_{1,0}}$};
\node[draw,circle] (2_0) at (2*\dx,0*\dy ) {\scriptsize ${x_{2,0}}$};
\node[draw,circle] (3_0) at (3*\dx,0*\dy ) {\scriptsize ${x_{3,0}}$};


    \foreach \j in {2,...,\k}{
        \pgfmathtruncatemacro{\UpperLimit}{\j-1}
        \foreach \i in {1,...,\UpperLimit}{
            \pgfmathtruncatemacro{\Nexti}{\i-1}
            \draw[->] (\Nexti_\j) -- (\i_\j);
        }
    }
    
    \foreach \j in {1,...,\k}{
        \pgfmathtruncatemacro{\LowerLimit}{\j-1}
        \foreach \i in {\LowerLimit,...,\Nmt}{
            \pgfmathtruncatemacro{\Nexti}{\i+1}
            \draw[<-] (\i_\j) -- (\Nexti_\j);
        }
    }

    \foreach \i in {1,...,\Nmt}{
        \pgfmathtruncatemacro{\LowerLimit}{\i}
        \foreach \j in {\LowerLimit,...,\Nmt}{
            \pgfmathtruncatemacro{\Nextj}{\j+1}
            \pgfmathtruncatemacro{\Newi} {\i-1}
            \draw[->] (\Newi_\j) -- (\Newi_\Nextj);
        }
    }
    
    \foreach \i in {1,...,\k}{
        \pgfmathtruncatemacro{\UpperLimit}{\i-1}
        \foreach \j in {0,...,\UpperLimit}{
            \pgfmathtruncatemacro{\Nextj}{\j+1}
            \pgfmathtruncatemacro{\Newi} {\i-1}
            \draw[<-] (\Newi_\j) -- (\Newi_\Nextj);
        }
    }

    \foreach \i in {2,...,\Nmt}{
    \pgfmathtruncatemacro{\UpperLimit}{\i-1}
        \foreach \m in {1,...,\UpperLimit}{
        	\pgfmathtruncatemacro{\Newi} {\i-1}
            \pgfmathtruncatemacro{\Currx}{\Newi-\m+1}
            \pgfmathtruncatemacro{\Curry}{\N-\m}
            \pgfmathtruncatemacro{\Nextx}{\Newi-\m}
            \pgfmathtruncatemacro{\Nexty}{\N-\m-1}
            \draw[->] (\Currx_\Curry) -- (\Nextx_\Nexty);
        }
    }
    
    \foreach \i in {2,...,\N}{
    \pgfmathtruncatemacro{\UpperLimit}{\N-\i+1}
        \foreach \m in {1,...,\UpperLimit}{
        	\pgfmathtruncatemacro{\Newi} {\i-1}
            \pgfmathtruncatemacro{\Currx}{\Newi+\m-2}
            \pgfmathtruncatemacro{\Curry}{\m-1}
            \pgfmathtruncatemacro{\Nextx}{\Newi+\m-1}
            \pgfmathtruncatemacro{\Nexty}{\m}
            \draw[->] (\Currx_\Curry) -- (\Nextx_\Nexty);
        }
    }
        

\foreach \j in {1,...,\k}{
\pgfmathtruncatemacro{\lm}{\k}

\draw[<-] (0_\j) to [out=160,in=20] (\lm_\j);
}

\foreach \j in {2,...,\k}{
\pgfmathtruncatemacro{\Prevj}{\j-1}
\pgfmathtruncatemacro{\lm}{\k}
\draw[->] (0_\j) to  (\lm_\Prevj);
}

\foreach \i in {1,...,\k}{
        	\pgfmathtruncatemacro{\Newi} {\i-1}
\draw[->] (\Newi_\k) to [out=75,in=-75] (\Newi_0);
}

\foreach \i in {1,...,\Nmt}{
\pgfmathtruncatemacro{\Nexti}{\i+1}
   \pgfmathtruncatemacro{\Newi} {\i-1}
\draw[<-] (\Newi_\k) to (\i_0);
}

\pgfmathtruncatemacro{\lm}{\k}
\draw[->] (0_0) to [out=90,in=-110] (\lm_\lm);    


\node[draw,circle,minimum size=.75cm] (A) at (4.5*\dx,-.5*\dy) {\tiny $x'_{A+B}$};

\node[draw,circle,minimum size=.75cm] (B) at (-1*\dx,0.5*\dy ) {$\tiny x_B$};


\pgfmathtruncatemacro{\lol}{\N} 
\pgfmathtruncatemacro{\wut}{\k}   
\draw[->]  (4_1)  to [out=-25,in=75] (A);
\draw[->]  (A) to [out=-115,in=-100]  (3_0);

\draw[->] (B)  to [out=90,in=170]  (0_1);
\draw[->]  (0_0) to [out=-150,in=-90] (B) ;

\end{tikzpicture}
\caption{Tropical variables for the seed $\mu_{(1,1)}\mu_{(0,1)}(\mathcal{Q})$}
\label{fig:quiver-after-11-10}
\end{figure}
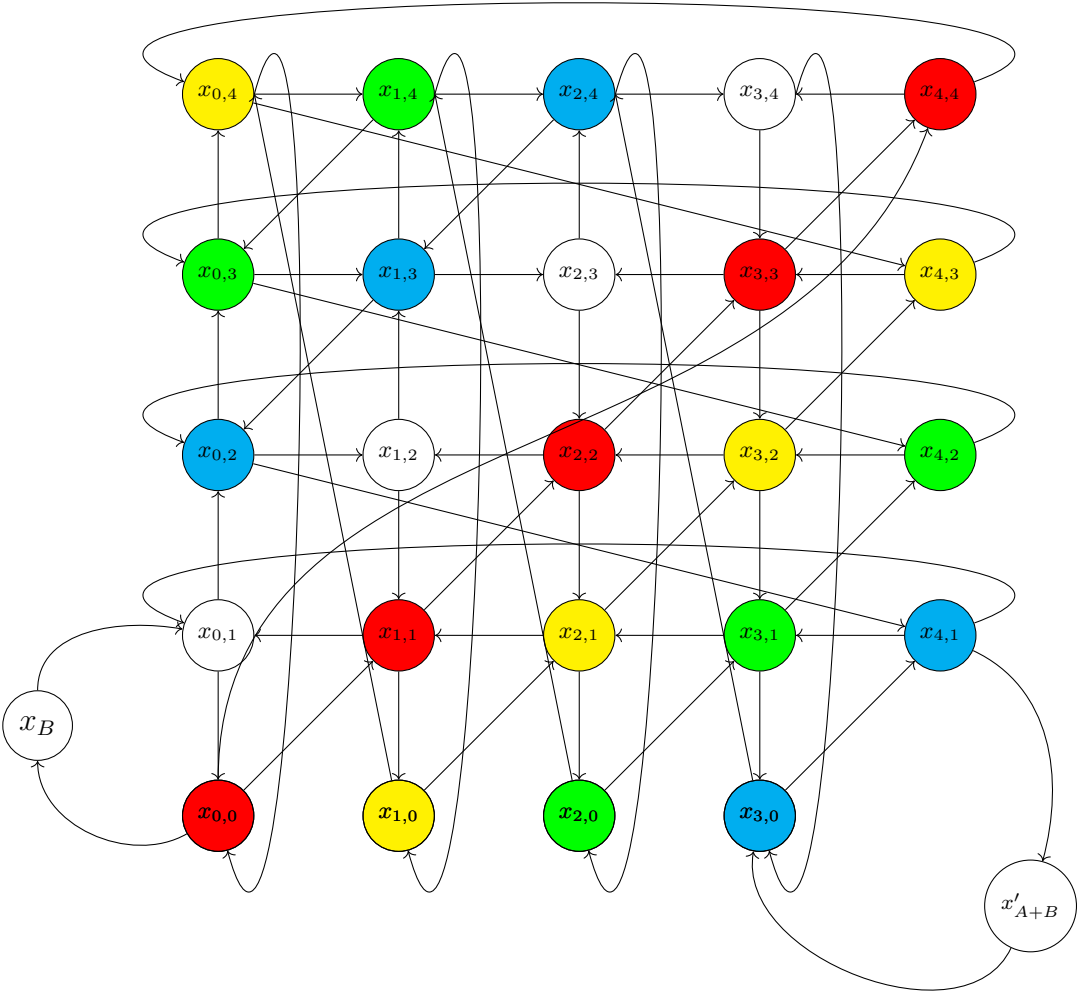

Indeed, note that in Figure~\ref{fig:quiver-after-10} the vertex with tropical variable $x_{A+B}$ is now attached to the $(1,1)$-curve, so that in this cluster the Baxter operator $Q_{(1,1)}$ is given by the automorphism part of the sign-coherent mutation sequence $\mu_{(1,1)}$ obtained in the same way as the sequence $\mu_{(0,1)}$ was in the initial cluster. The quiver $\mu_{(1,1)}\mu_{(0,1)}(\mathcal{Q})$ and its tropical variables are illustrated in Figure~\ref{fig:quiver-after-11-10} for $N=5$, where we have set
$$
x'_{A+B} = x_{A+B}\prod_{i=1}^{N-1}\prod_{k=0}^{N-1}x_{i+k-1,k},
$$
and we understand the index of $x_{i+k-1,k}$ in the product modulo $N$. The factors in the products over $k$ with fixed $i$ are indicated chromatically in Figure~\ref{fig:quiver-after-11-10}.

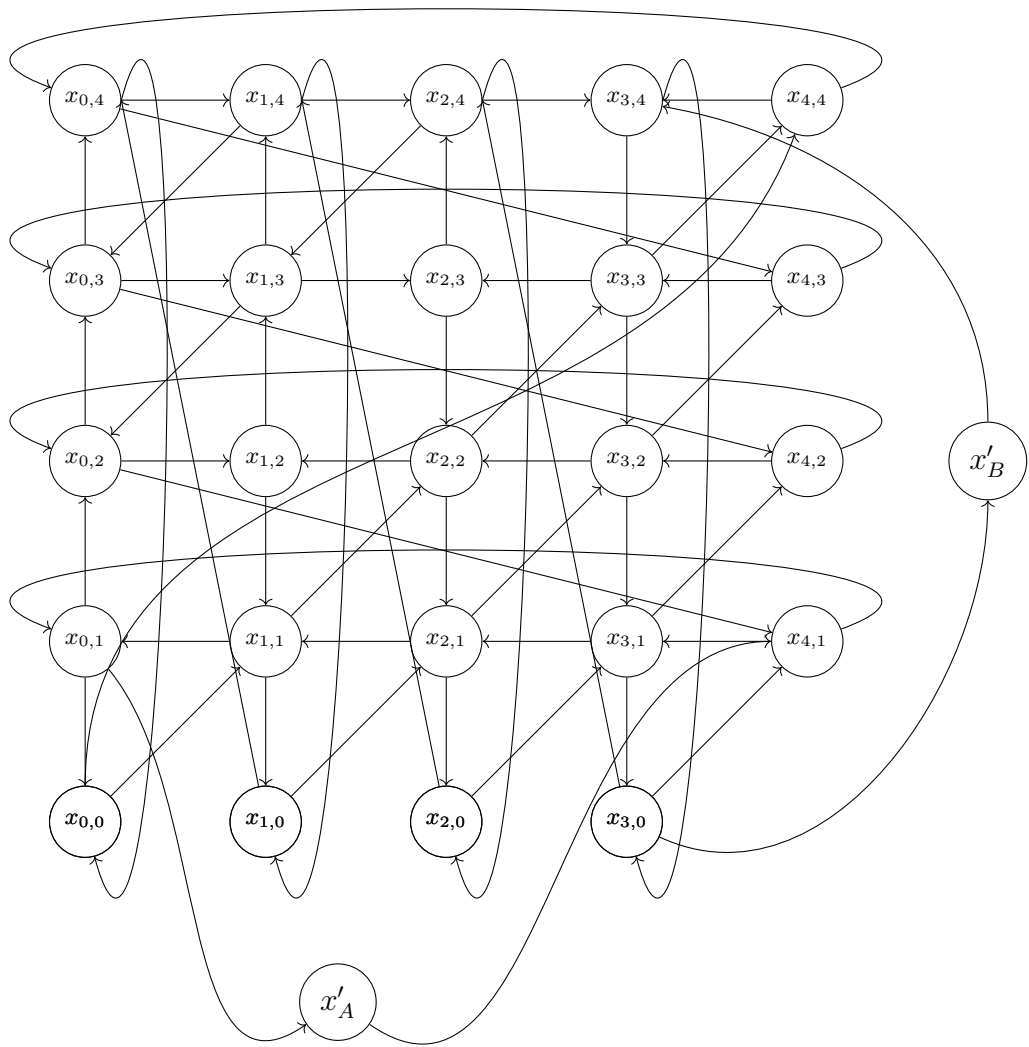
\begin{figure}[ht]

\qquad
\begin{tikzpicture}[scale=.6]

	\pgfmathsetmacro{\N}{5}
	\pgfmathsetmacro{\k}{\N-1}
	\pgfmathsetmacro{\Nmt}{\N-2}
    \pgfmathsetmacro{\dx}{4}
    \pgfmathsetmacro{\dy}{4}
    \pgfmathsetmacro{\eps}{0.2}
    \pgfmathsetmacro{\varoffset}{0.05}
     \pgfmathsetmacro{\baxoffset}{3}


\node[draw,circle] (1_1) at (1*\dx,1*\dy ) { \scriptsize${x_{1,1}}$};
\node[draw,circle] (2_2) at (2*\dx,2*\dy ) {\scriptsize ${x_{2,2}}$};
\node[draw,circle] (3_3) at (3*\dx,3*\dy ) { \scriptsize${x_{3,3}}$};
\node[draw,circle] (4_4) at (4*\dx,4*\dy ) { \scriptsize${x_{4,4}}$};
\node[draw,circle] (2_1) at (2*\dx,1*\dy ) {\scriptsize ${x_{2,1}}$};
\node[draw,circle] (3_2) at (3*\dx,2*\dy ) {\scriptsize ${x_{3,2}}$};
\node[draw,circle] (4_3) at (4*\dx,3*\dy ) { \scriptsize${x_{4,3}}$};
\node[draw,circle] (0_4) at (0*\dx,4*\dy ) {\scriptsize ${x_{0,4}}$};
\node[draw,circle] (3_1) at (3*\dx,1*\dy ) {\scriptsize ${x_{3,1}}$};
\node[draw,circle] (4_2) at (4*\dx,2*\dy ) {\scriptsize ${x_{4,2}}$};
\node[draw,circle] (0_3) at (0*\dx,3*\dy ) { \scriptsize${x_{0,3}}$};
\node[draw,circle] (1_4) at (1*\dx,4*\dy ) {\scriptsize ${x_{1,4}}$};
\node[draw,circle] (4_1) at (4*\dx,1*\dy ) {\scriptsize ${x_{4,1}}$};
\node[draw,circle] (0_2) at (0*\dx,2*\dy ) {\scriptsize ${x_{0,2}}$};
\node[draw,circle] (1_3) at (1*\dx,3*\dy ) { \scriptsize${x_{1,3}}$};
\node[draw,circle] (2_4) at (2*\dx,4*\dy ) {\scriptsize ${x_{2,4}}$};
\node[draw,circle] (0_1) at (0*\dx,1*\dy ) {\scriptsize ${x_{0,1}}$};
\node[draw,circle] (1_2) at (1*\dx,2*\dy ) {\scriptsize ${x_{1,2}}$};
\node[draw,circle] (2_3) at (2*\dx,3*\dy ) { \scriptsize${x_{2,3}}$};
\node[draw,circle] (3_4) at (3*\dx,4*\dy ) {\scriptsize ${x_{3,4}}$};
\node[draw,circle] (0_0) at (0*\dx,0*\dy ) {\scriptsize ${x_{0,0}}$};
\node[draw,circle] (1_0) at (1*\dx,0*\dy ) {\scriptsize ${x_{1,0}}$};
\node[draw,circle] (2_0) at (2*\dx,0*\dy ) {\scriptsize ${x_{2,0}}$};
\node[draw,circle] (3_0) at (3*\dx,0*\dy ) {\scriptsize ${x_{3,0}}$};
\node[draw,circle] (0_0) at (0*\dx,0*\dy ) {\scriptsize ${x_{0,0}}$};
\node[draw,circle] (1_0) at (1*\dx,0*\dy ) {\scriptsize ${x_{1,0}}$};
\node[draw,circle] (2_0) at (2*\dx,0*\dy ) {\scriptsize ${x_{2,0}}$};
\node[draw,circle] (3_0) at (3*\dx,0*\dy ) {\scriptsize ${x_{3,0}}$};


    \foreach \j in {2,...,\k}{
        \pgfmathtruncatemacro{\UpperLimit}{\j-1}
        \foreach \i in {1,...,\UpperLimit}{
            \pgfmathtruncatemacro{\Nexti}{\i-1}
            \draw[->] (\Nexti_\j) -- (\i_\j);
        }
    }
    
    \foreach \j in {1,...,\k}{
        \pgfmathtruncatemacro{\LowerLimit}{\j-1}
        \foreach \i in {\LowerLimit,...,\Nmt}{
            \pgfmathtruncatemacro{\Nexti}{\i+1}
            \draw[<-] (\i_\j) -- (\Nexti_\j);
        }
    }

    \foreach \i in {1,...,\Nmt}{
        \pgfmathtruncatemacro{\LowerLimit}{\i}
        \foreach \j in {\LowerLimit,...,\Nmt}{
            \pgfmathtruncatemacro{\Nextj}{\j+1}
            \pgfmathtruncatemacro{\Newi} {\i-1}
            \draw[->] (\Newi_\j) -- (\Newi_\Nextj);
        }
    }
    
    \foreach \i in {1,...,\k}{
        \pgfmathtruncatemacro{\UpperLimit}{\i-1}
        \foreach \j in {0,...,\UpperLimit}{
            \pgfmathtruncatemacro{\Nextj}{\j+1}
            \pgfmathtruncatemacro{\Newi} {\i-1}
            \draw[<-] (\Newi_\j) -- (\Newi_\Nextj);
        }
    }

    \foreach \i in {2,...,\Nmt}{
    \pgfmathtruncatemacro{\UpperLimit}{\i-1}
        \foreach \m in {1,...,\UpperLimit}{
        	\pgfmathtruncatemacro{\Newi} {\i-1}
            \pgfmathtruncatemacro{\Currx}{\Newi-\m+1}
            \pgfmathtruncatemacro{\Curry}{\N-\m}
            \pgfmathtruncatemacro{\Nextx}{\Newi-\m}
            \pgfmathtruncatemacro{\Nexty}{\N-\m-1}
            \draw[->] (\Currx_\Curry) -- (\Nextx_\Nexty);
        }
    }
    
    \foreach \i in {2,...,\N}{
    \pgfmathtruncatemacro{\UpperLimit}{\N-\i+1}
        \foreach \m in {1,...,\UpperLimit}{
        	\pgfmathtruncatemacro{\Newi} {\i-1}
            \pgfmathtruncatemacro{\Currx}{\Newi+\m-2}
            \pgfmathtruncatemacro{\Curry}{\m-1}
            \pgfmathtruncatemacro{\Nextx}{\Newi+\m-1}
            \pgfmathtruncatemacro{\Nexty}{\m}
            \draw[->] (\Currx_\Curry) -- (\Nextx_\Nexty);
        }
    }
        

\foreach \j in {1,...,\k}{
\pgfmathtruncatemacro{\lm}{\k}

\draw[<-] (0_\j) to [out=160,in=20] (\lm_\j);
}

\foreach \j in {2,...,\k}{
\pgfmathtruncatemacro{\Prevj}{\j-1}
\pgfmathtruncatemacro{\lm}{\k}
\draw[->] (0_\j) to  (\lm_\Prevj);
}

\foreach \i in {1,...,\k}{
        	\pgfmathtruncatemacro{\Newi} {\i-1}
\draw[->] (\Newi_\k) to [out=75,in=-75] (\Newi_0);
}

\foreach \i in {1,...,\Nmt}{
\pgfmathtruncatemacro{\Nexti}{\i+1}
   \pgfmathtruncatemacro{\Newi} {\i-1}
\draw[<-] (\Newi_\k) to (\i_0);
}

\pgfmathtruncatemacro{\lm}{\k}
\draw[->] (0_0) to [out=90,in=-110] (\lm_\lm);    


\node[draw,circle,minimum size=.75cm] (A) at (1.4*\dx,-1*\dy) {$\tiny x_A'$};

\node[draw,circle,minimum size=.75cm] (B) at (5*\dx,2*\dy ) {$x'_B$};




\pgfmathtruncatemacro{\lol}{\N} 
\pgfmathtruncatemacro{\wut}{\k}   
\draw[->]  (0_1)  to [out=-50,in=-145] (A);
\draw[->]  (A) to [out=-35,in=-180]  (4_1);

\draw[->]  (3_0) to [out=-25,in=-90]  (B);
\draw[->]  (B) to [out=90,in=-10] (3_4) ;

\end{tikzpicture}
\caption{Tropical variables for the seed $\mu_{(1,0)}\mu_{(1,1)}\mu_{(0,1)}(\mathcal{Q})=\mu_{(0,1)}\mu_{(1,0)}(\mathcal{Q})$. Here we have abbreviated 
$
x'_A =x_A^{-1}\prod_{i=0}^{N-1}\prod_{j=1}^{N-1}x_{i,j}^{-1} , \quad x'_B=x_B^{-1}\prod_{i=1}^{N-1}\prod_{j=0}^{N-1}x_{i-1,j}^{-1}
$.
}
\label{fig:quiver-final-longway}
\end{figure}

After applying the mutations $\mu_{(1,1)}$, the vertex carrying tropical $\cX$-coordinate $ x_B$ is now re-attached to the $(1,0)$-curve, and thus in this cluster we may realize $Q_{(1,0)}$ as the automorphism part of the mutation sequence $\mu_{(1,0)}$. The seed $\mu_{(1,0)}\mu_{(1,1)}\mu_{(0,1)}(\mathcal{Q})$ is illustrated in Figure~\ref{fig:quiver-final-longway}, where we have introduced another abbreviation
\begin{align}
\label{eq:btilde}
x'_B=x_B^{-1}\prod_{i=1}^{N-1}\prod_{j=0}^{N-1}x_{i-1,j}^{-1}.
\end{align}

\begin{figure}[ht]

\qquad
\begin{tikzpicture}[scale=.6]

	\pgfmathsetmacro{\N}{5}
	\pgfmathsetmacro{\k}{\N-1}
	\pgfmathsetmacro{\Nmt}{\N-2}
    \pgfmathsetmacro{\dx}{4}
    \pgfmathsetmacro{\dy}{4}
    \pgfmathsetmacro{\eps}{0.2}
    \pgfmathsetmacro{\varoffset}{0.05}
     \pgfmathsetmacro{\baxoffset}{3}


\node[draw,circle] (1_1) at (1*\dx,1*\dy ) { \scriptsize${x_{1,1}}$};
\node[draw,circle] (2_2) at (2*\dx,2*\dy ) {\scriptsize ${x_{2,2}}$};
\node[draw,circle] (3_3) at (3*\dx,3*\dy ) { \scriptsize${x_{3,3}}$};
\node[draw,circle] (4_4) at (4*\dx,4*\dy ) { \scriptsize${x_{4,4}}$};
\node[draw,circle] (2_1) at (2*\dx,1*\dy ) {\scriptsize ${x_{2,1}}$};
\node[draw,circle] (3_2) at (3*\dx,2*\dy ) {\scriptsize ${x_{3,2}}$};
\node[draw,circle] (4_3) at (4*\dx,3*\dy ) { \scriptsize${x_{4,3}}$};
\node[draw,circle] (0_4) at (0*\dx,4*\dy ) {\scriptsize ${x_{0,4}}$};
\node[draw,circle] (3_1) at (3*\dx,1*\dy ) {\scriptsize ${x_{3,1}}$};
\node[draw,circle] (4_2) at (4*\dx,2*\dy ) {\scriptsize ${x_{4,2}}$};
\node[draw,circle] (0_3) at (0*\dx,3*\dy ) { \scriptsize${x_{0,3}}$};
\node[draw,circle] (1_4) at (1*\dx,4*\dy ) {\scriptsize ${x_{1,4}}$};
\node[draw,circle] (4_1) at (4*\dx,1*\dy ) {\scriptsize ${x_{4,1}}$};
\node[draw,circle] (0_2) at (0*\dx,2*\dy ) {\scriptsize ${x_{0,2}}$};
\node[draw,circle] (1_3) at (1*\dx,3*\dy ) { \scriptsize${x_{1,3}}$};
\node[draw,circle] (2_4) at (2*\dx,4*\dy ) {\scriptsize ${x_{2,4}}$};
\node[draw,circle] (0_1) at (0*\dx,1*\dy ) {\scriptsize ${x_{0,1}}$};
\node[draw,circle] (1_2) at (1*\dx,2*\dy ) {\scriptsize ${x_{1,2}}$};
\node[draw,circle] (2_3) at (2*\dx,3*\dy ) { \scriptsize${x_{2,3}}$};
\node[draw,circle] (3_4) at (3*\dx,4*\dy ) {\scriptsize ${x_{3,4}}$};
\node[draw,circle] (0_0) at (0*\dx,0*\dy ) {\scriptsize ${x_{0,0}}$};
\node[draw,circle] (1_0) at (1*\dx,0*\dy ) {\scriptsize ${x_{1,0}}$};
\node[draw,circle] (2_0) at (2*\dx,0*\dy ) {\scriptsize ${x_{2,0}}$};
\node[draw,circle] (3_0) at (3*\dx,0*\dy ) {\scriptsize ${x_{3,0}}$};
\node[draw,circle] (0_0) at (0*\dx,0*\dy ) {\scriptsize ${x_{0,0}}$};
\node[draw,circle] (1_0) at (1*\dx,0*\dy ) {\scriptsize ${x_{1,0}}$};
\node[draw,circle] (2_0) at (2*\dx,0*\dy ) {\scriptsize ${x_{2,0}}$};
\node[draw,circle] (3_0) at (3*\dx,0*\dy ) {\scriptsize ${x_{3,0}}$};


    \foreach \j in {2,...,\k}{
        \pgfmathtruncatemacro{\UpperLimit}{\j-1}
        \foreach \i in {1,...,\UpperLimit}{
            \pgfmathtruncatemacro{\Nexti}{\i-1}
            \draw[->] (\Nexti_\j) -- (\i_\j);
        }
    }
    
    \foreach \j in {1,...,\k}{
        \pgfmathtruncatemacro{\LowerLimit}{\j-1}
        \foreach \i in {\LowerLimit,...,\Nmt}{
            \pgfmathtruncatemacro{\Nexti}{\i+1}
            \draw[<-] (\i_\j) -- (\Nexti_\j);
        }
    }

    \foreach \i in {1,...,\Nmt}{
        \pgfmathtruncatemacro{\LowerLimit}{\i}
        \foreach \j in {\LowerLimit,...,\Nmt}{
            \pgfmathtruncatemacro{\Nextj}{\j+1}
            \pgfmathtruncatemacro{\Newi} {\i-1}
            \draw[->] (\Newi_\j) -- (\Newi_\Nextj);
        }
    }
    
    \foreach \i in {1,...,\k}{
        \pgfmathtruncatemacro{\UpperLimit}{\i-1}
        \foreach \j in {0,...,\UpperLimit}{
            \pgfmathtruncatemacro{\Nextj}{\j+1}
            \pgfmathtruncatemacro{\Newi} {\i-1}
            \draw[<-] (\Newi_\j) -- (\Newi_\Nextj);
        }
    }

    \foreach \i in {2,...,\Nmt}{
    \pgfmathtruncatemacro{\UpperLimit}{\i-1}
        \foreach \m in {1,...,\UpperLimit}{
        	\pgfmathtruncatemacro{\Newi} {\i-1}
            \pgfmathtruncatemacro{\Currx}{\Newi-\m+1}
            \pgfmathtruncatemacro{\Curry}{\N-\m}
            \pgfmathtruncatemacro{\Nextx}{\Newi-\m}
            \pgfmathtruncatemacro{\Nexty}{\N-\m-1}
            \draw[->] (\Currx_\Curry) -- (\Nextx_\Nexty);
        }
    }
    
    \foreach \i in {2,...,\N}{
    \pgfmathtruncatemacro{\UpperLimit}{\N-\i+1}
        \foreach \m in {1,...,\UpperLimit}{
        	\pgfmathtruncatemacro{\Newi} {\i-1}
            \pgfmathtruncatemacro{\Currx}{\Newi+\m-2}
            \pgfmathtruncatemacro{\Curry}{\m-1}
            \pgfmathtruncatemacro{\Nextx}{\Newi+\m-1}
            \pgfmathtruncatemacro{\Nexty}{\m}
            \draw[->] (\Currx_\Curry) -- (\Nextx_\Nexty);
        }
    }
        

\foreach \j in {1,...,\k}{
\pgfmathtruncatemacro{\lm}{\k}

\draw[<-] (0_\j) to [out=160,in=20] (\lm_\j);
}

\foreach \j in {2,...,\k}{
\pgfmathtruncatemacro{\Prevj}{\j-1}
\pgfmathtruncatemacro{\lm}{\k}
\draw[->] (0_\j) to  (\lm_\Prevj);
}

\foreach \i in {1,...,\k}{
        	\pgfmathtruncatemacro{\Newi} {\i-1}
\draw[->] (\Newi_\k) to [out=75,in=-75] (\Newi_0);
}

\foreach \i in {1,...,\Nmt}{
\pgfmathtruncatemacro{\Nexti}{\i+1}
   \pgfmathtruncatemacro{\Newi} {\i-1}
\draw[<-] (\Newi_\k) to (\i_0);
}

\pgfmathtruncatemacro{\lm}{\k}
\draw[->] (0_0) to [out=90,in=-110] (\lm_\lm);    


\node[draw,circle,minimum size=.75cm] (A) at (\k*\dx-\dx/2,\N*\dy ) {$x_A$ };


\node[draw,circle,minimum size=.75cm] (B) at (5*\dx,2*\dy ) {$x'_B$};



\pgfmathtruncatemacro{\lol}{\N} 
\pgfmathtruncatemacro{\wut}{\k}   
\draw[->]  (3_4)  to [out=90,in=-180] (A);
\draw[->]  (A) to [out=0,in=90]  (4_4);

\draw[->]  (3_0) to [out=-25,in=-90]  (B);
\draw[->]  (B) to [out=90,in=-10] (3_4) ;


\draw[->]  (A) to [out=45,in=45]  (B);

\end{tikzpicture}
\caption{
Tropical variables for quiver $\mu_{(1,0)}(\mathcal{Q})$. }
\label{fig:quiver-after-01}
\end{figure}
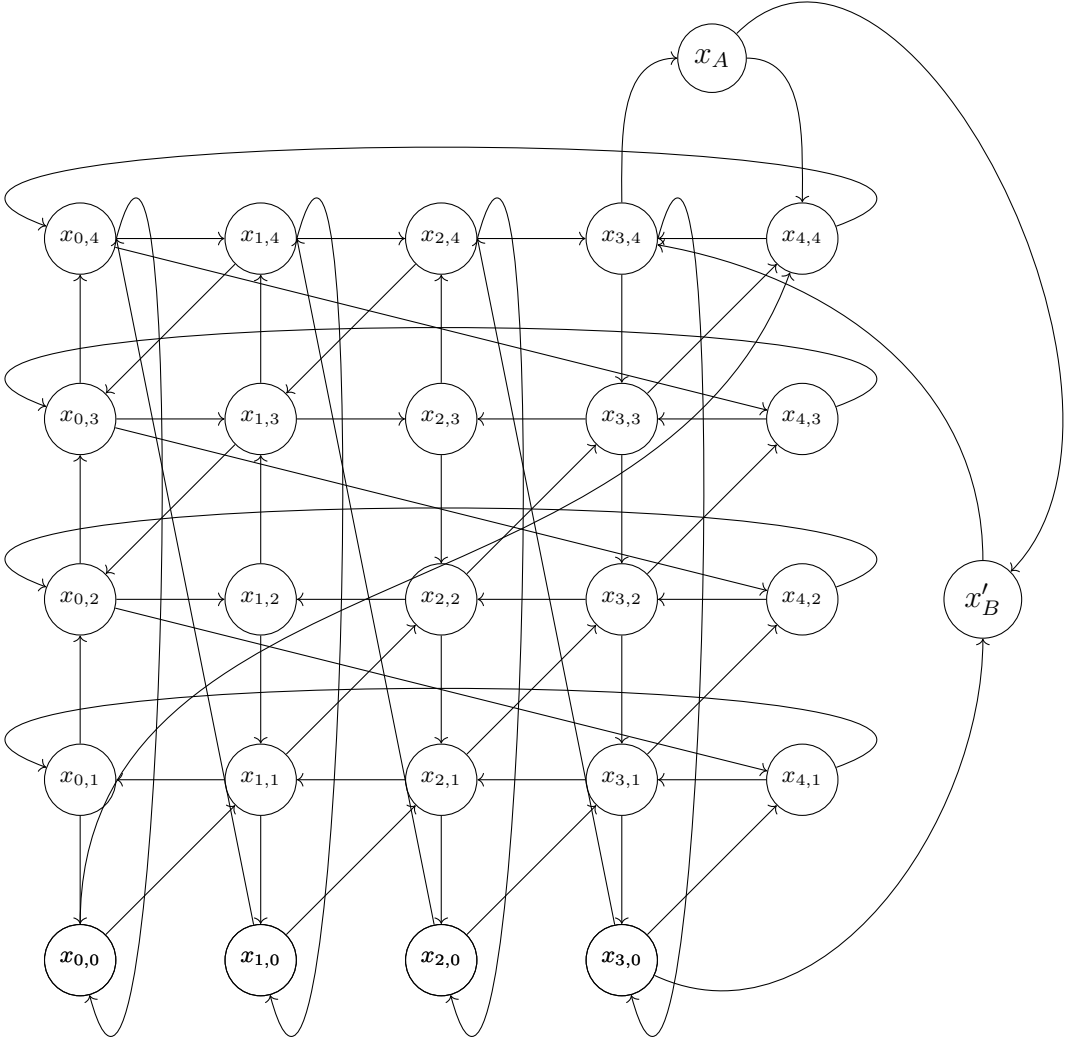

On the other hand, starting from the initial seed $\mathcal{Q}$ of Figure~\ref{fig:baxquiv-init}, the Baxter operator $Q_{(1,0)}$ is realized as the automorphism part of the mutation sequence $\mu_{(1,0)}$ associated to the vertex $x_B$ attached to the $(1,0)$-curve. This leads to the seed $\mu_{(1,0)}(\mathcal{Q})$ shown in Figure~\ref{fig:quiver-after-01}.

Finally, in this cluster we realize the Baxter operator $Q_{(0,1)}$ using the determinantal node $x_A$ attached to the $(0,1)$-curve, and observe that the set of tropical $\mathscr{X}$-variables in the resulting seed $\mu_{(0,1)}\mu_{(1,0)}(\mathcal{Q})$ is identical to that in $\mu_{(1,0)}\mu_{(1,1)}\mu_{(0,1)}(\mathcal{Q})$. Hence, as in the $N=2$ case, it follows from Theorem 3.5 of~\cite{KN11} that the non-commutative power series defining the automorphism parts of the corresponding quantum cluster transformations are equal, so that
$$
Q_{(1,0)}Q_{(0,1)} = Q_{(0,1)}Q_{(1,1)}Q_{(0,1)}
$$
and we thus obtain the pentagon identity for the Baxter operators in the rank $N$ skein.


.

\bibliographystyle{alpha}
\bibliography{refs.bib}

\end{document}